\pdfoutput=1
\RequirePackage{ifpdf}
\ifpdf % We are running pdfTeX in pdf mode
\documentclass[pdftex]{sigma}
\else
\documentclass{sigma}
\fi

\numberwithin{equation}{section}

\usepackage{kbordermatrix}
\usepackage{tikz}
\usetikzlibrary{calc,decorations.pathmorphing,decorations.markings,decorations.pathreplacing,matrix,fit,shapes.geometric}

\renewcommand\ss{\scriptstyle}

\newcommand{\ket}[1]{\left|#1\right\rangle}

\DeclareMathOperator{\codim}{codim}
\DeclareMathOperator{\Ker}{Ker}
\let\Im\relax\DeclareMathOperator{\Im}{Im}
\newcommand\ZZ{{\mathbb Z}}
\newcommand\QQ{{\mathbb Q}}
\newcommand\CC{{\mathbb C}}
\newcommand\PP{{\mathbb P}}
\newcommand\RR{{\mathbb R}}
\newcommand\Mat[1]{\operatorname{Mat}_{#1}}
\newcommand\cross{\mathit{cr}}
\newcommand\rk{\operatorname{rk}}
\newcommand\wt{\operatorname{wt}}
\newcommand\cl{{c\ell}}
\newcommand\TL{{\rm TL}}
\newcommand\CTL{\widetilde{{\rm TL}}}
\newcommand\LP{\mathcal L}
\newcommand\CLP{\widetilde{\mathcal L}}
\newcommand\CS{\mathit{CS}}
\newcommand\cyc[1]{\circlearrowright (#1)}

\newtheorem{thm}{Claim}[section]
\newtheorem{conj}[thm]{Conjecture}

\newtheorem{lem}[thm]{Lemma}
\newtheorem{cor}[thm]{Corollary}
\theoremstyle{definition}
\newtheorem{rmk}[thm]{Remark}
\newtheorem{ex}[thm]{Example}

\makeatletter
\newcommand{\gettikzxy}[3]{% returns coordinates
  \tikz@scan@one@point\pgfutil@firstofone#1\relax
\pgfmathsetmacro{#2}{\the\pgf@x/\linkpatternunit}
\pgfmathsetmacro{#3}{\the\pgf@y/\linkpatternunit}
}
\tikzset{label anchor/.code={%
    \let\tikz@auto@anchor=\pgfutil@empty
    \def\tikz@anchor{#1}
  },
  label anchor/.default=center
}
\makeatother
\tikzset{arrow/.style={postaction={decorate,thick,decoration={markings,mark = at position #1 with {\arrow{>}}}}},arrow/.default=0.5}
\tikzset{invarrow/.style={postaction={decorate,thick,decoration={markings,mark = at position #1 with {\arrow{<}}}}},invarrow/.default=0.5}
\newdimen\linkpatternunit%
\newcount\linkpatternsize%
\newcount\lpsize%the global internal variable... ugly but \foreach forces us to use a global variable
\newif\iflinkpatterninverted%mirror symmetry
\newif\iflinkpatterntikzstarted%in principle, shouldn't be needed -- tests itself
\newif\iflinkpatternboxed%draw a box around link patterns/tangles
\newif\iflinkpatternaxis%draw a line at starting/endpoints
\newif\iflinkpatternstraightlines%draw pipedream style
\newif\iflinkpatternnumbered%write labels
\newif\iflinkpatternalias%the vertices can be accessed via their label rather than their normal indexation. turn off for complicated labels
\newif\iflinkpatternnode%creates a node surrounding the linkpattern
\newif\iflinkpatterncentered%automatically centers the picture if outside \tikz{}
\newcount\linkpatternfused%for fused link patterns, level of fusion
\pgfkeys{/linkpattern/.cd,centered/.is if=linkpatterncentered,inverted/.is if=linkpatterninverted,numbered/.is if=linkpatternnumbered,tikzstarted/.is if=linkpatterntikzstarted,straight lines/.is if=linkpatternstraightlines,boxed/.is if=linkpatternboxed,axis/.is if=linkpatternaxis,vertexcolor/.store in=\linkpatternvertexcolor,edgecolor/.store in=\linkpatternedgecolor,boxcolor/.code={\linkpatternboxedtrue\def\linkpatternboxcolor{#1}},tikzoptions/.style={every linkpattern/.append style={#1}},size/.code={\linkpatternsize=#1},numbering0/.code={\def\lpnumbering{#1}\def\linkpatternnumbering{#1}},numbering/.style={numbered,numbering0={#1}},unit/.code={\linkpatternunit=#1},height/.store in=\linkpatternheight,shape/.store in=\linkpatternshape,looseness/.store in=\linkpatternlooseness,squareness/.store in=\linkpatternsquareness,extra space/.store in=\linkpatternextraspace,width/.store in=\linkpatternwidth,alias/.is if=linkpatternalias,pos/.store in=\linkpatternpos,labeloptions/.style={labeloptionslist/.append style={#1}},labeloptionslist/.style={inner sep=2pt,font=\scriptsize,execute at begin node=$,execute at end node=$,label anchor={#1+180}},nodeon/.is if=linkpatternnode,node/.style={nodeon,nodeoptionslist/.append style={#1}},nodeoptionslist/.style={},
pipedream/.style={shape=pipedream,looseness=0,straight lines,numbering0=tangle},
tangle/.style={shape=tangle,numbering0=tangle},
every linkpattern/.style={x=\linkpatternunit,y=\linkpatternunit},
fused/.code={\linkpatternfused=#1},%semi-implemented
vertex/.style={circle,thin,solid,draw=black,fill=\linkpatternvertexcolor,inner sep=1.5pt,draw opacity=1,transform shape},%note the transform shape to avoid all kinds of pbls
edge/.style={very thick,solid,draw=\linkpatternedgecolor,draw opacity=1}}
\linkpatterncenteredfalse%
\linkpatterninvertedfalse%
\linkpatternnumberedfalse%
\linkpatterntikzstartedfalse%
\linkpatternboxedfalse%
\linkpatternaxistrue%
\linkpatternaliastrue%
\linkpatternstraightlinesfalse%
\def\linkpatternlooseness{0.2}% how much arches go straight to their partner
\def\linkpatternsquareness{0.35}% how much arches go orthogonally
\def\linkpatternvertexcolor{red}%
\def\linkpatternedgecolor{blue}%
\def\linkpatternboxcolor{none}%
\def\linkpatternheight{0}
\def\linkpatternwidth{0}
\def\linkpatternshape{default}
\def\linkpatternnumbering{default}
\def\linkpatternpos{(0,0)}
\def\linkpatternextraspace{0}
\linkpatternnodefalse
\def\firstchar#1#2\empty{#1}%
\def\linkpatterndo#1#2{
\edef\param{\csname linkpattern#2\endcsname}
\edef\firstcharparam{\expandafter\firstchar\param\empty}
\expandafter\ifcat\firstcharparam a
\expandafter\ifx\csname linkpattern#1\param\endcsname\relax
\csname linkpattern#1unknown\endcsname
\else
\csname linkpattern#1\csname linkpattern#2\endcsname\endcsname
\fi
\else
\csname linkpattern#1unknown\endcsname
\fi
}%
\def\linkpatterncoordtangle{\ifnum\x>\lphalfsize\pgfmathparse{\lpsize+1-\x}\xdef\lpcoordx{\pgfmathresult}\xdef\lpcoordy{\lpheight}\xdef\lpangle{270}\else\xdef\lpcoordx{\x}\xdef\lpcoordy{-\lpheight}\xdef\lpangle{90}\fi}
\def\linkpatterncoordpipedream{\ifnum\x>\lphalfsize\pgfmathparse{\lpsize+1-\x-0.5}\xdef\lpcoordx{\pgfmathresult}\xdef\lpcoordy{0}\xdef\lpangle{270}\else\pgfmathparse{0.5-\x}\xdef\lpcoordy{\pgfmathresult}\xdef\lpcoordx{0}\xdef\lpangle{0}\fi}
\def\linkpatterncoordrectangle{
\ifnum\x>\lptqsize
\pgfmathparse{\lpsize+1-\x-0.5}\xdef\lpcoordx{\pgfmathresult}\xdef\lpcoordy{0}\xdef\lpangle{270}
\else\ifnum\x>\lphalfsize
\pgfmathparse{\x-\lptqsize-0.5}\xdef\lpcoordy{\pgfmathresult}\xdef\lpcoordx{\linkpatternwidth}\xdef\lpangle{180}
\else\ifnum\x>\linkpatternheight
\pgfmathparse{\x-\linkpatternheight-0.5}\xdef\lpcoordx{\pgfmathresult}\xdef\lpcoordy{-\linkpatternheight}\xdef\lpangle{90}
\else
\pgfmathparse{0.5-\x}\xdef\lpcoordy{\pgfmathresult}\xdef\lpcoordx{0}\xdef\lpangle{0}
\fi\fi\fi
}%
\def\linkpatternsetsizeunknown{
\global\lpsize=\linkpatternsize%this is the size as given by option size
\if\linkpatternheight0
\xdef\maxsep{0}
\foreach \x/\xx in \mylist%
{%
\edef\tempx{\withoutprime{\x}}
\edef\tempxx{\withoutprime{\xx}}
\pgfmathparse{max(\maxsep,abs(\tempx-\tempxx))}
\xdef\maxsep{\pgfmathresult}
}%
\pgfmathparse{0.25+0.8*\linkpatternsquareness*\maxsep}
\xdef\lpheight{\pgfmathresult}
\else
\xdef\lpheight{\linkpatternheight}
\fi
}
\newcount\tempsize
\def\linkpatternrightmostunknown{
\global\lpsize=0
\global\tempsize=0
\foreach\x/\labx in \linkpatternnumbering
{
\edef\tempx{\withoutprime{\x}}%not great... (prime not treated)
\ifnum\lpsize<\tempx\global\lpsize=\tempx\fi
\global\advance\tempsize by 1
}
\ifnum\tempsize>\lpsize\global\lpsize=\tempsize\fi
}%
\def\linkpatternrightmostdefault{
\global\lpsize=0
\global\tempsize=0
\foreach \x/\y in \mylist
{
\edef\tempx{\withoutprime{\x}}%not great... (prime not treated)
\ifnum\lpsize<\tempx\global\lpsize=\tempx\fi
\ifx\x\y
\global\advance\tempsize by 1
\else
\edef\tempy{\withoutprime{\y}}
\ifnum\lpsize<\tempy\global\lpsize=\tempy\fi%
\global\advance\tempsize by 2
\fi
}
\ifnum\tempsize>\lpsize\global\lpsize=\tempsize\fi
}%
\def\linkpatternrightmosttangle{
\global\lpsize=0
\global\tempsize=0
\foreach \x/\y in \mylist
{
\edef\tempx{\withoutprime{\x}}
\ifnum\lpsize<\tempx\global\lpsize=\tempx\fi
\ifx\x\y
\global\advance\tempsize by 1
\else
\edef\tempy{\withoutprime{\y}}
\ifnum\lpsize<\tempy\global\lpsize=\tempy\fi%
\global\advance\tempsize by 2
\fi
}
\global\advance\lpsize by\lpsize% not so great either: guess doubling because of 's
\ifnum\tempsize>\lpsize\global\lpsize=\tempsize\fi
}%

\newcommand\linkpattern[2][]{%universal command
{%we want the macro modifications to be local only
\pgfkeys{/linkpattern/.cd,#1}%parse the list of optional arguments
\edef\mylist{#2}%to avoid some annoying bugs
\def\primetest##1'{}%
\def\hasaprime##1{\expandafter\primetest##1''}%aha!
\def\internalwithoutprime##1'{##1}%
\def\withoutprime##1{\if\hasaprime##1 %
\expandafter\internalwithoutprime##1\else ##1\fi}%
\iflinkpatternnumbered%
\iflinkpatterninverted%the label anchoring hack forces us to distinguish cases
\tikzset{/linkpattern/lbl/.style n args={3}{label={[/linkpattern/labeloptionslist=-##1,##3] ##1:##2}}}%
\else%
\tikzset{/linkpattern/lbl/.style n args={3}{label={[/linkpattern/labeloptionslist=##1,##3] ##1:##2}}}%
\fi%
\else%
\tikzset{/linkpattern/lbl/.style={}}%
\fi%
\tikzifinpicture{\linkpatterntikzstartedtrue%
\begin{scope}[shift=\linkpatternpos,/linkpattern/every linkpattern]%should it be there? can't decide. and order of shift and every...
}{%
\linkpatterntikzstartedfalse%
\iflinkpatterncentered
\begin{tikzpicture}[baseline=(current  bounding  box.center),/linkpattern/every linkpattern]%
\else
\begin{tikzpicture}[baseline=0,/linkpattern/every linkpattern]%
\fi
}%
\begin{scope}[local bounding box=link pattern box]
\iflinkpatterninverted%
\begin{scope}[yscale=-1]%
\fi%
\linkpatterndo{setsize}{shape}%in case geometry dictates size
\ifnum\lpsize=0
\linkpatterndo{rightmost}{numbering}
\fi
\pgfmathtruncatemacro{\lphalfsize}{\lpsize/2}%should probably be transferred somewhere above
\linkpatterndo{numbering}{numbering}
\iflinkpatternboxed
\linkpatterndo{drawbox}{shape}
\else
\iflinkpatternaxis
\linkpatterndo{drawaxis}{shape}
\fi
\fi
\foreach\xx/\xlab/\opt in \lpnumbering
{
\ifx\xlab\opt\def\opt{}\fi%if no option avoid repeat of \xlab
\if\hasaprime\xx %
\pgfmathtruncatemacro{\xx}{\lpsize+1-\withoutprime{\xx}}
\fi
\ifnum\linkpatternfused>1
\pgfmathsetmacro{\x}{0.4*(0.5+\linkpatternfused*(0.5+floor((\xx-1)/\linkpatternfused)))+0.6*\xx}
\else
\def\x{\xx}
\fi
\linkpatterndo{coord}{shape}
\iflinkpatternalias\def\xlabb{\xlab}\else\def\xlabb{\xx}\fi
\path (\lpcoordx,\lpcoordy) coordinate[/linkpattern/vertex,/linkpattern/lbl={\lpangle+180}{\xlab}{\opt},alias=v\xlabb] (v\xx) ++(\lpangle:\linkpatternunit) coordinate[alias=vv\xlabb] (vv\xx); 
}
\foreach \a/\b/\c in \mylist
{
\if\hasaprime\a %
\pgfmathtruncatemacro{\a}{\lpsize+1-\withoutprime{\a}}
\fi
\ifx\b\c\def\c{}\fi%set option to empty if they're not there
\draw[/linkpattern/edge]
\ifx\a\b
(v\a)
\c
--
++(0,\lpheight);
\else
\pgfextra{
\if\hasaprime\b %
\pgfmathtruncatemacro{\b}{\lpsize+1-\withoutprime{\b}}
\fi
\gettikzxy{(v\a)}{\ax}{\ay}
\gettikzxy{(v\b)}{\bx}{\by}
\gettikzxy{(vv\a)}{\axx}{\ayy}
\gettikzxy{(vv\b)}{\bxx}{\byy}
\pgfmathsetmacro{\dist}{sqrt((\ax-\bx)*(\ax-\bx)+(\ay-\by)*(\ay-\by))}
\pgfmathsetmacro{\abx}{(\axx-\ax)*\dist*\linkpatternsquareness+(\bx-\ax)*\linkpatternlooseness)}
\pgfmathsetmacro{\aby}{(\ayy-\ay)*\dist*\linkpatternsquareness+(\by-\ay)*\linkpatternlooseness)}
\pgfmathsetmacro{\bax}{(\bxx-\bx)*\dist*\linkpatternsquareness+(\ax-\bx)*\linkpatternlooseness)}
\pgfmathsetmacro{\bay}{(\byy-\by)*\dist*\linkpatternsquareness+(\ay-\by)*\linkpatternlooseness)}
}
(v\a)
\c%option
\iflinkpatternstraightlines
\pgfextra{
\pgfmathsetmacro{\t}{((\ax-\bx)*\bay-(\ay-\by)*\bax)/(\aby*\bax-\abx*\bay)}
\pgfmathsetmacro{\abx}{\t*\abx}
\pgfmathsetmacro{\aby}{\t*\aby}
}
[rounded corners=0.2\linkpatternunit] -- ++(\abx,\aby) -- (v\b);
\else
.. controls ++(\abx,\aby) and ++(\bax,\bay) .. 
\fi
(v\b);
\fi
}
\end{scope}
\iflinkpatternnode
\node[fit=(link pattern box),/linkpattern/nodeoptionslist] {};
\fi
\iflinkpatterninverted
\end{scope}
\fi
\iflinkpatterntikzstarted
\end{scope}
\else%
\end{tikzpicture}%
\fi%
}}%
\newcommand\tanglelinkpattern[3][]{%
{%we want the macro modifications to be local only
\pgfkeys{/linkpattern/.cd,#1}
\iflinkpatterninverted
\begin{tikzpicture}[/linkpattern/every linkpattern,baseline=\linkpatternunit]%
\else
\begin{tikzpicture}[/linkpattern/every linkpattern,baseline=-\linkpatternunit]%
\fi
\linkpattern[#1,tikzstarted,numbered=false]{#3}
\pgfmathtruncatemacro{\lptempsize}{2*\linkpatternsize}
\iflinkpatterninverted
\begin{scope}[yshift=0.5*\linkpatternunit]
\else
\begin{scope}[yshift=-0.5*\linkpatternunit]
\fi
\linkpattern[tangle,#1,tikzstarted,size=\lptempsize,
numbering=halftangle,
height=0.5]{#2}
\end{scope}
\end{tikzpicture}%
}}
\newcommand\diag[4][]{%
\pgfkeys{/linkpattern/.cd,#1}%parse the list of optional arguments
\iflinkpatterntikzstarted\else%
\begin{tikzpicture}[scale=0.5]%TEMP
\fi%
\iflinkpatterninverted%
\begin{scope}[yscale=-1]%
\fi%
\draw (0,0) grid (#2,#3);
\edef\mylist{#4}%to avoid some annoying bugs
\foreach\y/\x/\z in \mylist
{
\ifx\x\z
\draw[decorate,decoration={zigzag,
amplitude=1pt,segment length=5pt}]
(\x-0.5,#3) -- (\x-0.5,\y-0.5) node[circle,fill=black,inner sep=2pt] {} -- (#2,\y-0.5);
\else
\node at (\x-0.5,\y-0.5) {$\z$};
\fi
}
\iflinkpatterninverted
\end{scope}
\fi
\iflinkpatterntikzstarted\else%
\end{tikzpicture}%
\fi%
}
\makeatletter
\tikzset{circle split part fill/.style  args={#1,#2}{%
 alias=tmp@name,
  postaction={%
    insert path={
     \pgfextra{% 
     \pgfpointdiff{\pgfpointanchor{\pgf@node@name}{center}}%
                  {\pgfpointanchor{\pgf@node@name}{east}}%            
     \pgfmathsetmacro\insiderad{\pgf@x}
      \fill[#1] (\pgf@node@name.base) ([xshift=-\pgflinewidth]\pgf@node@name.east) arc
                          (0:180:\insiderad-\pgflinewidth)--cycle;
      \fill[#2] (\pgf@node@name.base) ([xshift=\pgflinewidth]\pgf@node@name.west)  arc
                           (180:360:\insiderad-\pgflinewidth)--cycle;                    }}}}}  
 \makeatother
\tikzset{bdot/.style={circle,circle split,draw,circle split part fill={black,white},thin,inner sep=1pt}}%
\tikzset{wdot/.style={circle,circle split,draw,circle split part fill={white,black},thin,inner sep=1pt}}%
\newcommand\circlelinkpattern[2][]{%in principle, could be merged
{%we want the macro modifications to be local only
\pgfkeys{/linkpattern/.cd,#1}%parse the list of optional arguments
\iflinkpatterntikzstarted\else%
\begin{tikzpicture}[/linkpattern/every linkpattern]%
\fi%
\iflinkpatterninverted%
\begin{scope}[yscale=-1]%
\fi%
\global\lpsize=\linkpatternsize
\edef\mylist{#2}%to avoid some annoying bugs
\foreach \x/\y in \mylist
{
\ifnum\x>\lpsize\global\lpsize=\x\fi
\ifnum\y>\lpsize\global\lpsize=\y\fi
}
\iflinkpatternaxis
\draw (0,0) circle (1);
\fi
\foreach\x in {1,...,\lpsize}
{
\pgfmathparse{(0.3*floor((\x-1)/\linkpatternfused)+0.7*((\x-0.5)/\linkpatternfused-0.5))*\linkpatternfused*360/\lpsize}
\coordinate[/linkpattern/vertex] (v\x) at (\pgfmathresult:1);
}
\foreach \x/\y/\z in \mylist
{
\ifx\y\z%
\draw[/linkpattern/edge] (v\x) .. controls ($0.5*(v\x)$) and  ($0.5*(v\y)$) .. (v\y);
\else
\draw[/linkpattern/edge] \z (v\x) .. controls ($0.5*(v\x)$) and  ($0.5*(v\y)$) .. (v\y);
\fi
}
\iflinkpatternnumbered%
\pgfmathparse{\lpsize/\linkpatternfused}
\global\lpsize=\pgfmathresult
\def\linkpatternnumbering{1,...,\lpsize}
\newdimen\angle
\foreach\x/\xx/\opt in \linkpatternnumbering
{
  \pgfmathsetmacro{\angle}{360/\lpsize*(\x-1)}
\ifx\xx\opt%
  \node[outer sep=1pt,anchor=180+\angle] at (\angle:1) {$\scriptstyle\xx$}; % note the subtle anchoring: if we used for example "label", the anchoring would be restricted to the cardinal points i.e. north, north east, etc. instead we're using angle as anchor (documented in sec 39)
\else
  \node[outer sep=1pt,anchor=180+\angle,\opt] at (\angle:1) {$\scriptstyle\xx$}; % note the subtle anchoring: if we used for example "label", the anchoring would be restricted to the cardinal points i.e. north, north east, etc. instead we're using angle as anchor (documented in sec 39)
\fi
}
\fi%
\iflinkpatterninverted%
\end{scope}
\fi%
\iflinkpatterntikzstarted\else%
\end{tikzpicture}%
\fi%
}}%
\newdimen{\loopcellsize}
\tikzset{bgplaq/.style={draw=black,fill=\linkpatternboxcolor}}
\def\plaqwest{}
\def\plaqeast{}
\def\plaqnorth{}
\def\plaqsouth{}
\def\plaqa{
\draw[/linkpattern/edge,\plaqwest,\plaqnorth] (0,0.5) .. controls (0,0.2) and (-0.2,0) .. (-0.5,0);
\draw[/linkpattern/edge,\plaqeast,\plaqsouth] (0,-0.5) .. controls (0,-0.2) and (0.2,0) .. (0.5,0);
}
\def\plaqb{
\draw[/linkpattern/edge,\plaqwest,\plaqsouth] (0,-0.5) .. controls (0,-0.2) and (-0.2,0) .. (-0.5,0);
\draw[/linkpattern/edge,\plaqeast,\plaqnorth] (0,0.5) .. controls (0,0.2) and (0.2,0) .. (0.5,0);
}
\def\plaqc{
\draw[/linkpattern/edge,\plaqsouth,\plaqnorth] (0,0.5) -- (0,-0.5);
\draw[/linkpattern/edge,\plaqwest,\plaqeast] (0.5,0) -- (-0.5,0);
}
\pgfkeys{/linkpattern/.cd,west/.store in=\plaqwest,east/.store in=\plaqeast,north/.store in=\plaqnorth,south/.store in=\plaqsouth}
\def\plaqname{plaq}%default name
\newcommand\plaq[2][]{
\node[bgplaq,rectangle,draw,use as bounding box,minimum size=\loopcellsize,transform shape] (\plaqname) {};
\pgfkeys{/linkpattern/.cd,#1}%parse the list of optional arguments
\ifx#2\empty\else
\begin{scope}[x=\loopcellsize,y=\loopcellsize]
\csname plaq#2\endcsname
\end{scope}\fi
}
\newcommand\halfplaq[1][]{
\node[draw=none,fill=none,rectangle,use as bounding box,minimum size=\loopcellsize,transform shape] (\plaqname) {};
\pgfkeys{/linkpattern/.cd,#1}%parse the list of optional arguments
\begin{scope}[x=\loopcellsize,y=\loopcellsize]
\draw[bgplaq] (-0.5,-0.5) -- ++(1,1) -- ++(-1,0) -- cycle;
\draw[/linkpattern/edge,\plaqwest,\plaqnorth] (0,0.5) .. controls (0,0.2) and (-0.2,0) .. (-0.5,0);
\end{scope}
}
\tikzset{loop/.code={\def\plaqname{loop-\the\pgfmatrixcurrentrow-\the\pgfmatrixcurrentcolumn}},loop/.append style={matrix,row sep={\loopcellsize,between origins},column sep={\loopcellsize,between origins}}}%note: define loop-row-col as name for each elementary plaquette

\linkpatterninvertedtrue
\linkpatternnumberedtrue
\def\linkpatternheight{1}
\def\linkpatternboxcolor{pink!20!white}
\def\linkpatternextraspace{0.5}
\newcount\leftcount
\pgfkeys{/linkpattern/every linkpattern/.add code={\global\leftcount=1}}
\newcommand\toleft{{[rounded corners=0.2\linkpatternunit] |- (0.5,\leftcount*0.4); \global\advance\leftcount by 1}}
\newcommand\toleftx[1]{{[rounded corners] |- (0.5,#1*0.4); \global\advance\leftcount by 1}}
\pgfkeys{/linkpattern/.cd,small/.style={numbered=false,tikzoptions={scale=0.5}}}
\colorlet{cola}{green!50!black}
\colorlet{colb}{red!50!black}
\newdimen{\cellsize}
\newcommand\medboxes{\setlength{\cellsize}{14.22pt}\def\boxformat{}}%exactly 0.5cm

\medboxes
\tikzset{tableaubox/.style={draw=black,thin,sharp corners,solid,minimum size=\cellsize,inner sep=0pt}}
\tikzset{tableau/.style={matrix,name=tab,matrix anchor=tab-1-1.south west,inner sep=1pt,matrix of math nodes,cells={anchor=center,draw=black,thin,solid,arrows=-},nodes={tableaubox,execute at begin node=\boxformat},nodes in empty cells,row sep={\cellsize,between origins},column sep={\cellsize,between origins}}}
\makeatletter
\newcommand\cellextra[1]{#1\expandafter\tikz@lib@matrix@start@cell}%trick: after doing the extra stuff we restart the node reading macro. hopefully expandafter is ok
\makeatother
\def\activate#1{\begingroup
  \lccode`\~=`#1%
  \lowercase{\endgroup \let~#1}%
  \catcode`#1=13\relax}
\activate &

\tikzset{bgplaq/.style={draw=black,dotted,fill=\linkpatternboxcolor}}

\long\def\junk#1{}
\def\lozXY(#1,#2){
\begin{scope}[x={(0.866cm,-0.5cm)},y={(0.866cm,0.5cm)}]
\draw[fill=cyan!30!white] (#1,#2) -- ++(1,0) -- ++(0,1) -- ++(-1,0) -- cycle;
\end{scope}
}
\def\lozX(#1,#2){
\begin{scope}[x={(0cm,-1cm)},y={(0.866cm,0.5cm)},shift={(\a*0.866cm+0.866cm,\a*0.5cm+0.5cm)}]
\draw[fill=pink] (#1,#2) -- ++(1,0) -- ++(0,1) -- ++(-1,0) -- cycle;
\end{scope}
}
\def\lozY(#1,#2){
\begin{scope}[x={(0.866cm,-0.5cm)},y={(0cm,1cm)},shift={(\a*0.866cm+0.866cm,-\a*0.5cm-0.5cm)}]
\draw[fill=orange!50!white] (#1,#2) -- ++(1,0) -- ++(0,1) -- ++(-1,0) -- cycle;
\end{scope}
}
\newcommand{\loz}[3]{%
\begin{scope}[scale=0.6,xshift=-\a*0.866cm-\b*0.433cm-\c*0.433cm-1.732cm,yshift=\b*0.25cm-\c*0.25cm]
\def\la{#1}\def\lb{#2}\def\lc{#3}
\foreach\co in \la { \expandafter\lozXY\co }
\foreach\co in \lb { \expandafter\lozX\co }
\foreach\co in \lc { \expandafter\lozY\co }
\end{scope}
}

\begin{document}
%\allowdisplaybreaks

\newcommand{\arXivNumber}{1612.05361}

\renewcommand{\thefootnote}{}

\renewcommand{\PaperNumber}{069}

\FirstPageHeading

\ShortArticleName{Loop Models and $K$-Theory}

\ArticleName{Loop Models and $\boldsymbol{K}$-Theory\footnote{This paper is a~contribution to the Special Issue on Combinatorics of Moduli Spaces: Integrability, Cohomo\-logy, Quantisation, and Beyond. The full collection is available at \href{http://www.emis.de/journals/SIGMA/moduli-spaces-2016.html}{http://www.emis.de/journals/SIGMA/moduli-spaces-2016.html}}}

\Author{Paul ZINN-JUSTIN}

\AuthorNameForHeading{P.~Zinn-Justin}

\Address{School of Mathematics and Statistics, The University of Melbourne, Victoria 3010, Australia}
\Email{\href{mailto:pzinn@unimelb.edu.au}{pzinn@unimelb.edu.au}}
\URLaddress{\url{http://blogs.unimelb.edu.au/paul-zinn-justin/}}

\ArticleDates{Received November 28, 2017, in final form June 27, 2018; Published online July 13, 2018}

\Abstract{This is a review/announcement of results concerning the connection between certain exactly solvable two-dimensional models of statistical mechanics, namely loop models, and the equivariant $K$-theory of the cotangent bundle of the Grassmannian. We interpret various concepts from integrable systems ($R$-matrix, partition function on a finite domain) in geometric terms. As a byproduct, we provide explicit formulae for $K$-classes of various coherent sheaves, including structure and (conjecturally) square roots of canonical sheaves and canonical sheaves of conormal varieties of Schubert varieties.}

\Keywords{quantum integrability; loop models; $K$-theory}

\Classification{14M15; 82B23}

\tableofcontents

\renewcommand{\thefootnote}{\arabic{footnote}}
\setcounter{footnote}{0}

\section{Introduction}
Loop models are an important class of lattice models in two-dimensional statistical mechanics. They display a broad range of critical phenomena, and in fact many classical models (such as the 2D Ising model) are equivalent to a loop model. A subclass of loop models is particularly interesting: these are the {\em exactly solvable} ones, or equivalently the ones that display the pro\-perty of quantum integrability, under the form of the Yang--Baxter equation satisfied by their Boltzmann weights. The main reason that these exactly solvable/quantum integrable models are studied is of course because one can perform various exact calculations in them, that are in general not available. However, for our purposes, there is another reason to consider them, which is the recently discovered connection between quantum integrable systems and generalized cohomology theories (see \cite{GKV-coh} for the first hint of such a connection). In fact, loop models were among the first in which this connection was made explicit \cite{artic34, artic33}; much later, a~framework for vertex models was set up in \cite{MO-qg}; see also \cite{AO-ell, GK-K,GRTV,RTV-K}.

So far, the work on the ``geometry'' of loop models (in the sense of the connection above) has only focused on ordinary cohomology. A first step towards $K$-theory was made in \cite{artic69}; here, we continue in this direction, as well as propose new ones. We first explain the ``natural'' way to generalize from cohomology to $K$-theory the results of~\cite{artic34,artic56}; this is the same approach as in~\cite{artic69}, and though it differs from the vertex model approach in several important aspects, it still follows the same very general philosophy: we define a certain basis of the $K$-theory of the ambient space,\footnote{The $K$-theoretic basis considered here is related to the stable basis defined in~\cite{Ok-K} and further studied in \cite{SZZ-Kstable} by a~triangular matrix of maximal parabolic Kazhdan--Lusztig polynomials.} which in almost all of this paper will be the cotangent bundle of the Grassmannian, and then consider the action of the Weyl group on it. This means that, on the geometric side, we shall define certain coherent sheaves~$\sigma_\pi$, which will lead us on the integrable side to the so-called (Temperley--Lieb) {\em noncrossing loop model}. We shall then further depart from this philosophy by asking questions about {\em structure sheaves} of conormal varieties of Schubert varieties (in short, conormal Schubert varieties). This will lead us this time to a {\em crossing loop model}.

A significant part of the paper will be concerned with enlarging the dictionary between geometry and integrability. In particular, we shall discuss in some detail (extending joint work with A.~Knutson \cite{artictt}) the interpretation of the {\em partition function} of the various loop models on an arbitrary finite domain.

It should be noted since no proofs are provided in this paper, all new results should technically be considered as conjectures at this stage. We differentiate below ``Claims'', for which the idea of proof should be clear, from proper ``Conjectures''.

As motivation for what follows, we now provide two such claims, which are byproducts of the framework that is developed here.

The first one concerns certain explicit formulae in the equivariant $K$-theory of the cotangent bundle of the Grassmannian; we state it in words only here:
\begin{thm}%{thm*}
The equivariant $K$-classes of structure sheaves of conormal varieties of Schubert varieties in the Grassmannian $($resp.\ of the sheaves $\sigma_\pi$ supported on these varieties$)$ are given by partition functions of the trigonometric crossing loop model $($resp.\ noncrossing loop model$)$ on a~$k\times n$ rectangular domain with prescribed connectivity of boundary vertices.
\end{thm}%{thm*}

The details, including the construction of the sheaves $\sigma_\pi$, as well as the ``trigonometric'' weights of the models, will be given explicitly in Corollaries~\ref{cor:struc}, \ref{cor:sqrt}, along with the choice of boundary conditions (connectivity of the boundary vertices). Modulo a certain conjecture (Conjecture~\ref{conj:CM}), $\sigma_\pi$ is the square root of the canonical sheaf of its support, and $K$-classes of canonical sheaves of conormal Schubert varieties are also given by partition functions of the crossing loop model, see Corollary~\ref{cor:canon}.

The second one is an explicit description of a certain Hilbert series. Set $n=2k$, where $k$ is a positive integer. Consider the $n\times n$ matrix $J$ made of $k$ Jordan blocks of size~$2$:
\begin{gather*}
J=\begin{pmatrix}
0&1\\
0&0\\
 & &\ddots\\
 & & & 0 & 1\\
 & & & 0 & 0
\end{pmatrix}.
\end{gather*}
If $B_+$ is the group of invertible upper triangular matrices, then we can consider the $B_+$-orbit closure
\begin{gather*}
\mathcal O = \overline{\{xJx^{-1},\, x\in B_+\}}.
\end{gather*}
$\mathcal O$ is an affine variety, hence has a coordinate ring $R$. It is an easy exercise to check that $\mathcal O$ is invariant by scaling ($M\in\mathcal O\Rightarrow \lambda M\in\mathcal O$, $\lambda\in\CC$), so that $R$ is {\em graded}. We can therefore define its Hilbert series
\begin{gather*}
\chi(t)=\sum_{i=0}^\infty \dim R_i\, t^i.
\end{gather*}
\begin{thm}%{thm*}
\begin{gather*}
\chi(t)=\frac{t^{k(k+1)/4}P_k\big(t^{1/2}+t^{-1/2}\big)}{(1-t)^{k^2}},
\end{gather*}
where $P_k$ is a polynomial of degree $k(k+1)/2$ which is a weighted enumeration of {\em totally symmetric self-complementary plane partitions} of size $k-1$.
\end{thm}%{thm*}
The weighting in the enumeration will be detailed in Section~\ref{sec:TSSCPP}. In particular, $P_k(1)$ is the famous sequence
(\href{https://oeis.org/A005130}{A005130})
\begin{gather*}
P_k(1)=1,1,2,7,42,429,\ldots,\qquad k=1,2,\ldots,
\end{gather*}
which enumerates alternating sign matrices, totally symmetric self-complementary plane partitions and descending plane partitions. The argument~$1$ of~$P_k$ corresponds to the formal para\-me\-ter~$t$ being evaluated at a nontrivial cubic root of unity, which is somewhat mysterious from a~geometric standpoint.

\section{Setup}
\subsection{The Temperley--Lieb algebra}
Define the Temperley--Lieb algebra $\TL_n$ to be the $\QQ(\beta)$-algebra with generators $e_i$, $i=1,\ldots$, $n-1$, and relations
\begin{alignat}{3}
& e_i^2 =\beta e_i,\qquad && 1\le i\le n-1,& \nonumber\\
& e_ie_{i\pm1}e_i =e_i,\qquad &&1\le i,i\pm 1\le n-1,&\nonumber\\
& e_i e_j = e_j e_i,\qquad && 1\le i,j\le n-1,\quad |i-j|>1.& \label{eq:TLalgrel}
\end{alignat}

$\TL_n$ has a well-known diagrammatic representation; a basis of the algebra is given by diagrams on the strip $\RR\times [0,1]$ with $n$ vertices at its top and bottom boundaries (corresponding to extreme values of the second coordinate) made of $n$ disjoint arcs (i.e., smooth embeddings of $[0,1]$ into the strip) connecting vertices. In particular, the generators are depicted as
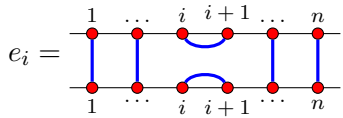
\begin{gather*}
e_i=\linkpattern[centered,tangle,height=0.6,alias=false, numbering={1/1,1'/1,2/\ldots,2'/\ldots,3/i,3'/i,4/i+1,4'/i+1,5/\ldots,5'/\ldots,6/n,6'/n}]{1/1',2/2',3/4,3'/4',5/5',6/6'}
\end{gather*}
Product corresponds to vertical concatenation of diagrams, where reading from right to left corresponds to concatenation from bottom to top. Diagrams are considered up to continuous deformation, with the extra rule that whenever a closed loop is formed, it is erased at the expense of multiplying by~$\beta$, which we describe as
\begin{gather}\label{eq:graphrulesnc}
\tikz[baseline=-3pt]{\draw[/linkpattern/edge] circle (0.15cm);} = \beta.
\end{gather}
$\beta$ is often called the ``loop weight''.

Given $0\le k\le n$, consider the $\QQ(\beta)$-submodule $\mathcal I_{k,n}$ generated by diagrams in $\TL_n$ such that there exist a pairing among the leftmost~$k$ bottom vertices, or a pairing among the rightmost $n-k$ bottom vertices. It is easy to see that this forms in fact a~left submodule (ideal) of the left regular representation of $\TL_n$, so that one can define the quotient $\TL_n$-left module
\begin{gather*}
\mathcal H_{k,n}=\TL_n/\mathcal I_{k,n}.
\end{gather*}
$\mathcal H_{k,n}$ has a canonical basis, denoted $\ket{\pi}$, made of diagrams $\pi$ with no pairings among the leftmost~$k$ bottom vertices and no pairings among the rightmost $n-k$ bottom vertices. In fact, it is convenient to simplify the diagrammatic description of this basis by noting that one only needs to keep track of the connections of the top vertices, the bottom ones being forced by the no-pairing conditions. In these truncated diagrams, some top vertices remain unpaired; it is furthermore convenient (though the information is at this stage redundant) to remember which ones used to be connected to one of the $k$ leftmost bottom vertices by making them connect to left infinity, whereas the ones connected to one of the $n-k$ rightmost bottom vertices go to downwards infinity. (The strange asymmetry in this depiction will eventually be justified in Section~\ref{sec:back}.) On an example, the diagram simplification is as follows ($k=4$, $n=8$):
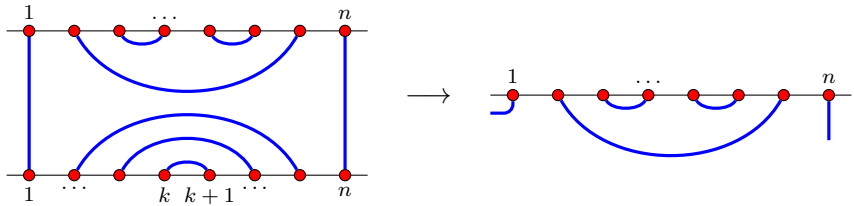
\begin{gather*}
\linkpattern[height=1.6,tangle,alias=false,numbering={1/1,2/,3/,4/\ldots,5/,6/,7/,8/n,1'/1,2'/\ldots,3'/,4'/k,5'/{k+1},6'/\ldots,7'/,8'/n}]{1/1',2/7,3/4,5/6,8/8',2'/7',3'/6',4'/5'}
\quad\longrightarrow\quad
\tikz[baseline=-3pt]{\linkpattern[alias=false,numbering={1/1,2/,3/,4/\ldots,5/,6/,7/,8/n}]{1/1/\toleft,2/7,3/4,5/6,8/8}}
\end{gather*}
Such (truncated) diagrams we call noncrossing {\em link patterns}, following the physics literature. Equivalently, we can define them as planar pairings of $n$ points on a line, with some vertices possibly left unpaired and at most $m:=\min(k,n-k)$ pairings. Their set is denoted $\LP_{k,n}$.

The rule for the $\TL_n$-action on link patterns is then as follows: $e_i$ acts in the natural way by reconnecting $i$, $i+1$ and inserting a~$(i,i+1)$ pairing if at least one of $i$, $i+1$ is paired (with a~weight of $\beta$ if they were in fact paired together), or if~$i$ is connected to left infinity and $i+1$ to bottom infinity (in the latter case, the resulting line from left to bottom infinities is erased); in all other cases, the result of the $e_i$ action is zero.

Furthermore, define for $\pi\in \LP_{k,n}$ $\cl(\pi)$ to be the subset of ``closings'' of $\pi$, that is the subset of $\{1,\ldots,n\}$ of vertices which are connected to left infinity or paired to a vertex left of them. On the example above, $\cl(\pi)=\{1,4,6,7\}$. $\cl$ is a~bijection between $\LP_{k,n}$ and $k$-subsets of $\{1,\ldots,n\}$.

Define $\rk(\pi)$ to be the number of pairings of $\pi$. There is a filtration of $\TL_n$-modules
\begin{gather*}
\mathcal H_{k,n}^{(\ge r)} = \operatorname{span} \big( \ket{\pi},\, \pi \in \LP_{k,n}^{(s)},\, s\ge r\big),
\qquad \LP_{k,n}^{(s)}=\{ \pi \in \LP_{k,n}\colon \rk(\pi)= s \},
\end{gather*}
since the $\TL_n$ action can only increase the number of pairings. In particular, the smallest such module $\mathcal H_{k,n}^{(\ge m)}$ will play a special role.

\subsection{The crossing algebra}\label{sec:crossalg}
We now present a second $\QQ(\beta)$-algebra, denoted $\CTL_n$, which contains $\TL_n$ as a subalgebra. It has generators $e_i$ and $f_i$, $i=1,\ldots,n-1$, and relations
\begin{gather}\label{eq:CTLalgrel}
e_i^2 =\beta e_i, \qquad f_i^2 =-f_i,\qquad f_ie_i=e_i f_i =-e_i,\\
e_ie_{i\pm 1}e_i =e_i,\qquad f_i f_{i+1} f_i = f_{i+1} f_i f_{i+1},\qquad f_i e_{i\pm 1} e_i = f_{i\pm 1} e_i, \qquad e_i e_{i\pm 1} f_i = e_i f_{i\pm 1},\nonumber\\
\notag e_i e_j=e_j e_i,\qquad f_i f_j = f_j f_i,\qquad e_i f_j = f_je_i,\qquad |i-j| >1.\nonumber
\end{gather}
with same range of indices as in~\eqref{eq:TLalgrel}. This algebra is somewhat similar to Brauer \cite{Brauer} and Birman--Wenzel--Murakami \cite{BW, Mura} algebras; it is also closely related to the degenerate Brauer algebra of~\cite[Section~5.5]{artic39}.

The graphical representation of $\CTL_n$, extending that of $\TL_n$, is given by diagrams on the strip $\RR\times [0,1]$ with $n$ vertices at its top and bottom boundaries made of $n$ arcs connecting vertices, {\em with regular crossings} inside the strip, such that arcs do not self-intersect and intersect each other at most once. We call these {\em reduced diagrams}. These diagrams are considered once again up to continuous deformation, and up to the move
\begin{gather}\label{eq:braid}
\vcenter{\hbox{\tikz \draw[/linkpattern/edge] (0,-0.5) -- (0,0.5) (-0.2,0.4) -- (0.5,-0.2) (-0.2,-0.4) -- (0.5,0.2); }}
=
\vcenter{\hbox{\tikz[rotate=180] \draw[/linkpattern/edge] (0,-0.5) -- (0,0.5) (-0.2,0.4) -- (0.5,-0.2) (-0.2,-0.4) -- (0.5,0.2); }}
\end{gather}
(which is necessary in order to stay within the class of diagrams with regular crossings).

Equivalence classes of reduced diagrams then form a basis of $\CTL_n$, where
the generators are given by
\begin{gather*}
e_i=\linkpattern[tangle,height=0.6,alias=false,numbering={1/1,1'/1,2/\ldots,2'/\ldots,3/i,3'/i,4/i+1,4'/i+1,5/\ldots,5'/\ldots,6/n,6'/n}]{1/1',2/2',3/4,3'/4',5/5',6/6'},
\qquad
f_i=\linkpattern[squareness=0,tangle,height=0.6,alias=false,numbering={1/1,1'/1,2/\ldots,2'/\ldots,3/i,3'/i,4/i+1,4'/i+1,5/\ldots,5'/\ldots,6/n,6'/n}]{1/1',2/2',3/4',3'/4,5/5',6/6'}
\end{gather*}

Multiplication corresponds to vertical concatenation, and the first three relations of \eqref{eq:CTLalgrel}
lead to the following additional rules:
\begin{gather}\label{eq:graphrules}
\tikz[baseline=-3pt]{\draw[/linkpattern/edge] circle (0.15cm);} = \beta
\qquad
\tikz[baseline=-3pt,rotate=90]{\draw[/linkpattern/edge,bend left=70] (-0.8,-0.3) to (0.8,-0.3); \draw[/linkpattern/edge,bend right=70] (-0.8,0.3) to (0.8,0.3);}
=
-\tikz[baseline=-3pt,rotate=90]{\draw[/linkpattern/edge] (-0.8,-0.3) to (0.8,0.3); \draw[/linkpattern/edge] (-0.8,0.3) to (0.8,-0.3);}
\qquad
\tikz[baseline=-3pt,rotate=0]{\draw[/linkpattern/edge] (-0.5,0) .. controls (0,0) and (0.5,0.5) .. (0,0.5) .. controls (-0.5,0.5) and (0,0) .. (0.5,0);}
=-\
\tikz[baseline=-3pt,rotate=0]{\draw[/linkpattern/edge] (-0.5,0) -- (0.5,0);}
\end{gather}
It is a somewhat nontrivial exercise to check that these rules are consistent, in other words, that (a) any diagram can be transformed into a reduced diagram by application of \eqref{eq:graphrules}, and (b) two reduced diagrams which can be related to each other by transformations \eqref{eq:graphrules} are equivalent, i.e., can be obtained from each other by continuous deformation and~\eqref{eq:braid}.

Given $0\le k\le n$, consider the $\QQ(\beta)$-submodule $\widetilde{\mathcal I}_{k,n}$ generated by diagrams in $\CTL_n$ such that there exist a pairing among the leftmost $k$ bottom vertices or a crossing among the arcs coming out of these vertices, or the same property for the rightmost $n-k$ bottom vertices. Once again, this forms in fact a left submodule (ideal) of the left regular representation of $\CTL_n$ (pairings among bottom vertices cannot be removed by left multiplication, and crossings can only be removed by pairing these vertices), so that one can define the quotient $\CTL_n$-left module
\begin{gather*}
\widetilde{\mathcal H}_{k,n}=\CTL_n/\widetilde{\mathcal I}_{k,n}.
\end{gather*}
$\widetilde{\mathcal H}_{k,n}$ has a canonical basis $\ket{\pi}$, made of diagrams $\pi$ with no pairings/no crossings among the leftmost~$k$ bottom vertices and among the rightmost $n-k$ bottom vertices. We can simplify the representation of these diagrams by keeping track of the connections of the top vertices only, the bottom ones being forced by the no-pairing/no-crossing conditions; however, it is now compulsory to remember which top vertices used to be connected to one of the $k$ leftmost bottom vertices (resp.\ $n-k$ rightmost bottom vertices) by connecting them to left infinity (resp.\ downwards infinity). Here is a $k=4$, $n=7$ example:
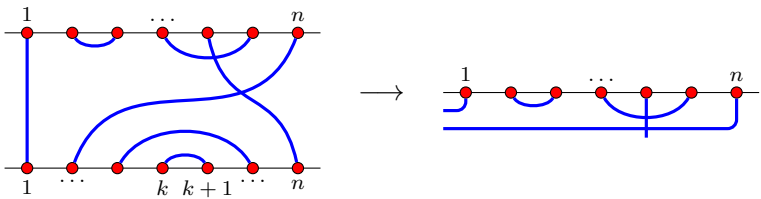
\begin{gather*}
\linkpattern[height=1.5,tangle,alias=false,numbering={1/1,2/,3/,4/\ldots,5/,6/,7/n,1'/1,2'/\ldots,3'/,4'/k,5'/{k+1},6'/\ldots,7'/n}]{1/1',2/3,4/6,7/2',3'/6',4'/5',5/7'}
%\quad\longrightarrow\quad
%\linkpattern{1/1,2/3,4/6,5,7/7}
\quad\longrightarrow\quad
\tikz[baseline=-3pt]{\linkpattern[alias=false,numbering={1/1,2/,3/,4/\ldots,5/,6/,7/n}]{1/1/\toleft,2/3,4/6,5,7/7/\toleft}}
\end{gather*}

We call these truncated diagrams crossing link patterns; their set is denoted $\CLP_{k,n}$.

The action of the algebra generators on such crossing link patterns is as follows. $e_i$ acts in the usual way if at least one of $i$, $i+1$ is paired, or if~$i$ and $i+1$ are connected to different infinities and {\em all unpaired vertices to the left of~$i$ are connected to left infinity}, {\em all unpaired vertices to the right of $i+1$ are connected to bottom infinity} (in which case the line from left to bottom is erased, possibly changing signs if it has a self-intersection); otherwise the result is zero. Similarly, $f_i$~acts in the usual way (crossing arcs coming from~$i$ and $i+1$, then applying rules~\eqref{eq:graphrules}) unless $i$ and $i+1$ are both connected to left infinity, or both to bottom infinity, in which case the result is zero.

There is again a filtration of $\CTL_n$-modules
\begin{gather*}
\widetilde{\mathcal H}_{k,n}^{(\ge r)} = \operatorname{span} \big( \ket{\pi},\, \pi \in \CLP_{k,n}^{(s)}, \, s\ge r\big),
\qquad \CLP_{k,n}^{(s)}=\big\{ \pi \in \CLP_{k,n}\colon \rk(\pi)=s \big\},
\end{gather*}
where we recall that $\rk(\pi)$ is the number of pairings of $\pi$.

Finally, denote by $\cross(\pi)$ the number of crossings of (any reduced diagram of) $\pi$.

\begin{ex}
Take $k=1$, $n=3$:
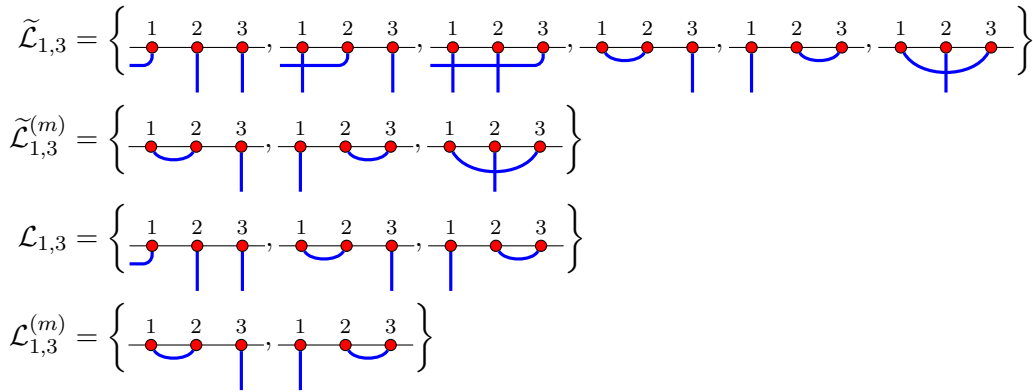
\begin{align*}
\CLP_{1,3}&=\bigg\{\linkpattern{1/1/\toleft,2/2,3/3},
\linkpattern{1/1,2/2/\toleft,3/3},
\linkpattern{1/1,2/2,3/3/\toleft},
\linkpattern{1/2,3/3},
\linkpattern{2/3,1/1},
\linkpattern{1/3,2/2}
\bigg\}
\\
\CLP^{(m)}_{1,3}&=\bigg\{
\linkpattern{1/2,3/3},
\linkpattern{2/3,1/1},
\linkpattern{1/3,2/2}
\bigg\}
\\
\LP_{1,3}&=\bigg\{\linkpattern{1/1/\toleft,2/2,3/3},
\linkpattern{1/2,3/3},
\linkpattern{2/3,1/1}
\bigg\}
\\
\LP^{(m)}_{1,3}&=\bigg\{
\linkpattern{1/2,3/3},
\linkpattern{2/3,1/1}
\bigg\}
\end{align*}
\end{ex}

\subsection[$B$-orbits, conormal Schubert varieties, and orbital varieties]{$\boldsymbol{B}$-orbits, conormal Schubert varieties, and orbital varieties}\label{sec:geomsetup}

Given nonnegative integers $k$ and $n$, $0\le k\le n$, consider the (complex) Grassmannian
\begin{gather}\label{eq:defGr}
\operatorname{Gr}_{k,n}=\big\{ V\le \CC^n\colon \dim V=k\big\}
\end{gather}
and its cotangent bundle
\begin{gather*}
T^\ast \operatorname{Gr}_{k,n}=\big\{(V,u)\in \operatorname{Gr}_{k,n}\times \operatorname{End}(\CC^n)\colon \Im u\subset V\subset \Ker u\big\}.
\end{gather*}
We define the two obvious maps $\mu\colon (V,u)\mapsto u$ and $p\colon (V,u)\mapsto V$. $\mu$ is a resolution of singularities of its image
\begin{gather*}
\mathcal N_{k,n}=\mu(T^\ast \operatorname{Gr}_{k,n})= \big\{ u\in \operatorname{End}(\CC^n)\colon u^2=0,\, \operatorname{rank}(u)\le m\big\},
\end{gather*}
where we recall that $m=\min(k,n-k)$.

The general linear group $\operatorname{GL}_n$ naturally acts on $\operatorname{Gr}_{k,n}$, $T^\ast \operatorname{Gr}_{k,n}$ and $\mathcal N_{k,n}$, making $\mu$ and $p$ equivariant. Inside $\operatorname{GL}_n$ sits the Borel subgroup $B_n$ (invertible lower triangular matrices), and the Cartan torus $T^{(0)}_n$ (invertible diagonal matrices). Define $T_n=T^{(0)}_n\times \CC^\times$, where the additional circle $\CC^\times$ acts on $T^\ast \operatorname{Gr}_{k,n}$ by scaling of the fiber. $T^{(0)}_n$-fixed points in $\operatorname{Gr}_{k,n}$, or equivalently, $T_n$-fixed points in $T^\ast \operatorname{Gr}_{k,n}$ (embedding $\operatorname{Gr}_{k,n}$ in $T^\ast \operatorname{Gr}_{k,n}$ as the image of the zero section) are coordinate subspaces $\CC^I$, where $I$ runs over $k$-subsets of $\{1,\ldots,n\}$.

Define $\mathcal O_{k,n} = \mathcal N_{k,n} \cap \mathfrak{n}_-$, where $\mathfrak{n}_-$ is the space of strict lower triangular matrices. $\mathcal O_{k,n}$ is a reducible affine scheme, known as the orbital scheme; its irreducible components are called {\em orbital varieties}. Similarly, $\mu^{-1}(\mathcal \mathcal O_{k,n})$ is a reducible scheme, whose irreducible components we call {\em conormal Schubert varieties}. In both cases, these irreducible components are Lagrangian (more on that in Section~\ref{sec:poisson}). In order to describe these, it is convenient to discuss the $B_n$-orbit decomposition of $\operatorname{Gr}_{k,n}$ and $\mathcal O_{k,n}$.

$B_n$-orbits in $\mathcal O_{k,n}$ ($B_n$ acting by conjugation) are indexed by involutions of $\{1,\ldots,n\}$ with at most~$m$ $2$-cycles (see \cite{artic39}, in particular Section~2, and references therein); given such an involution, a~representative of the corresponding orbit is the strict lower triangle of its permutation matrix. We can use link patterns to define such involutions (and therefore such matrices); namely, given a link pattern $\pi\in \CLP_{k,n}$, we can associate to it $\pi_>\in\operatorname{End}(\CC^n)$, where
\begin{gather*}
(\pi_>)_{i,j} = \begin{cases}
1 & \text{if $i$ and $j$ are paired in $\pi$, $i>j$},\\
0 & \text{otherwise}.
\end{cases}
\end{gather*}
Any $B_n$-orbit of $\mathcal O_{k,n}$ is of the form $B_n\cdot \pi_>$, $\pi\in \CLP_{k,n}$; however, note that $\pi_>$ does not distinguish between vertices that are connected to left infinity or bottom infinity by $\pi$, so that this does not provide a bijective labelling (except in the trivial cases $k=0,n$). For future use, we denote $\mathcal O_{\pi}:=\overline{B\cdot\pi_>}$. The defining equations of $\mathcal O_\pi$ are known (set-theoretically, and conjecturally, scheme-theoretically with their reduced
structure)~\cite{Rothbach-PhD}:
\begin{gather}\label{eq:orb}
\mathcal O_\pi = \big\{ u\in \operatorname{End}(\CC^n)\colon u^2=0,\, \rk(u_{i,j})_{i\ge i_0,\,j\le j_0}
\le \rk((\pi_>)_{i,j})_{i\ge i_0,\,j\le j_0},\, 1\le i_0,\,j_0\le n\big\}.\!\!\!
\end{gather}
These varieties were recently proved to be normal and to have rational singularities~\cite{BP-orb}.

$B_n$-orbits of $\operatorname{Gr}_{k,n}$ contain exactly one fixed point (coordinate subspace), and are known as Schubert cells: $S^o_I:=B_n\CC^I$, while their closures are called Schubert varieties: $S_I = \overline{S^o_I}$. Now given a link pattern $\pi\in \CLP_{k,n}$, one can consider the set of closings~$\cl(\pi)$ and then the associated Schubert variety $S_{\cl(\pi)}$. If $\pi$ runs over $\LP_{k,n}$, we get each Schubert variety exactly once as~$S_{\cl(\pi)}$. Schubert varieties are also known to be normal and to have rational singularities~\cite{Brion-flag}.

We can now describe irreducible components of $\mathcal O_{k,n}$ and $\mu^{-1}(\mathcal O_{k,n})$ as follows. Orbital varieties (irreducible components of $\mathcal O_{k,n}$) are exactly the $B_n$-orbit closures~$\mathcal O_\pi$, where $\pi\in \LP^{(m)}_{k,n}$, i.e., $\pi$ is noncrossing and has maximum number of pairings. In contrast, conormal Schubert varieties (irreducible components of $\mu^{-1}(\mathcal O_{k,n})$)
are of the form $\overline{\mu^{-1}(B_n\cdot \pi_>)}$ (being careful that the closure is taken after the preimage) where $\pi\in \LP_{k,n}$ (i.e., $\pi$ is noncrossing but has arbitrary number of pairings), in which case we denote it $\CS_{\cl(\pi)}$. Alternatively, $\CS_I$ ($I$ $k$-subset of $\{1,\ldots,n\}$) can be defined as the closure of the conormal bundle $\CS^o_I$ of $S^o_I$; or as the unique component of $\mu^{-1}(\mathcal O_{k,n})$ which satisfies $p(\CS_I)=S_I$.

We also believe the following to be true:
\begin{conj}\label{conj:CM}
$\CS_I$ is Cohen--Macaulay and normal for all $I$.
\end{conj}
Although our results do not directly depend on this conjecture, their interpretation in terms of canonical sheaves
(as mentioned in the abstract and the introduction) do.

\subsection{Cohomology theories}
We wish to study the equivariant cohomology and $K$-theory of $T^\ast \operatorname{Gr}_{k,n}$.

Let us first discuss cohomology. Since we are interested in neither ring structure nor grading, we shall simply denote by $H_{T_n}(T^\ast \operatorname{Gr}_{k,n})$ the localized $T_n$-equivariant cohomology ring of~$T^\ast \operatorname{Gr}_{k,n}$, considered as a vector space over $H_{T_n}(\cdot)$, the localized equivariant cohomology ring of a point. Explicitly, $H_{T_n}(\cdot)=\QQ(x_1,\ldots,x_n,\hbar)$, where $\hbar$ is the generator of the Lie algebra of the $\CC^\times$ scaling the fiber, and $x_1,\ldots,x_n$ are the obvious coordinates on the Lie algebra of the Cartan torus $T^{(0)}_n$ of $\operatorname{GL}_n$, and ``localized'' means we are tensoring over $\ZZ[x_1,\ldots,x_n,\hbar]$ with the fraction field $H_{T_n}(\cdot)=\QQ(x_1,\ldots,x_n,\hbar)$). $H_{T_n}(T^\ast \operatorname{Gr}_{k,n})$ is of dimension $n\choose k$, and a possible basis is given by classes of the conormal Schubert varieties $\CS_I$.

The entries $u_{i,j}$ of $u$, parameterizing the fiber of $T^\ast \operatorname{Gr}_{k,n}$, are eigenvectors of the $T_n$-action, with (additive) weights
\begin{gather}\label{eq:weiH}
\wt_H(u_{i,j})=\hbar-x_i+x_j.
\end{gather}

Similarly, in equivariant $K$-theory, we define $K_{T_n}(T^\ast \operatorname{Gr}_{k,n})$ to be the localized $K$-theory ring (or $K$-homology) of $T^\ast \operatorname{Gr}_{k,n}$, viewed as a module over $K_{T_n}(\cdot)$, which is $\QQ(z_1,\ldots,z_n,t)$, where $t$ is the coordinate on the $\CC^\times$ scaling the fiber, and $z_1,\ldots,z_N$ are coordinates on $T^{(0)}_n$. $K_{T_n}(T^\ast \operatorname{Gr}_{k,n})$ is of course also of dimension $n\choose k$. Besides the basis of classes of (structure sheaves of) conormal Schubert varieties (used in Section~\ref{sec:K2}), we shall use one more basis in what follows (see Section~\ref{sec:K1}).

We could use multiplicative notations to describe the weights of the entries $u_{i,j}$ as
\begin{gather}\label{eq:weiK}
\wt_K(u_{i,j})=t z_i^{-1}z_j.
\end{gather}

A slight subtlety will arise (in Section~\ref{sec:K1}) in that we shall sometimes wish to use the square root of $t$; to formalize this, one can introduce a double cover of $T_n$, say $T'_n$, also of the form $T^{(0)}_n\times \CC^\times$, but where the generator of the $\CC^\times$, conveniently denoted $t^{1/2}$, acts by scaling of the fiber as $\big(t^{1/2}\big)^2$. $K_{T'}(T^\ast \operatorname{Gr}_{k,n})$ is much the same as $K_{T_n}(T^\ast \operatorname{Gr}_{k,n})$, but is now a vector space over $K_{T'}(\cdot)=\QQ\big(z_1,\ldots,z_n,t^{1/2}\big)$.

In all cases, restriction $i^\ast$ to fixed points is an isomorphism (giving up to normalization the expansion in the basis of fixed points), and we shall use it for computations in examples.

Also note that we have $H_{T_n}(T^\ast \operatorname{Gr}_{k,n})\cong H(T^\ast \operatorname{Gr}_{k,n}) \otimes_\QQ \QQ(x_1,\ldots,x_n,\hbar)$, where $H(T^\ast \operatorname{Gr}_{k,n})$ is the nonequivariant localized cohomology ring, and similarly in $K$-theory.

\section[Loop models and $R$-matrices]{Loop models and $\boldsymbol{R}$-matrices}\label{sec:Rmat}
\subsection{Cohomology}\label{sec:Rmatcoh}
Write $\TL_n(\beta_0)$ for the specialization $\TL_n/\left<\beta-\beta_0\right>$ (i.e., the loop weight is fixed to be $\beta_0$), and similar notations for modules. We shall give $H(T^\ast \operatorname{Gr}_{k,n})$ the structure of the $\TL_n(2)$-right module $\mathcal H_{k,n}(2)^\ast$ by identifying $[\CS_{\cl(\pi)}]$ with the {\em dual} basis of the canonical basis $\ket{\pi}$ of $\mathcal H_{k,n}(2)$. Similarly, $H_{T_n}(T^\ast \operatorname{Gr}_{k,n})$ becomes identified with $\mathcal H_{k,n}(2)^\ast \otimes \QQ(x_1,\ldots,x_n,\hbar)$.

\subsubsection[The geometric $R$-matrix]{The geometric $\boldsymbol{R}$-matrix}
We are now ready to describe the Weyl group action. It is defined geometrically as the following right action:
\begin{gather*}
\big[w^{-1} A\big]=[A] \mathcal R_w,
\end{gather*}
where $w$ is an element of the Weyl group $W=N(T)/T$, and acts geometrically as any representative. $\mathcal R_w$ is {\em not} $H_{T_n}(\cdot)$-linear: rather,
\begin{gather*}
(f v)\mathcal R_w = (f\tau_w) (v\mathcal R_w),\qquad f\in H_{T_n}(\cdot),\qquad v\in H_{T_n}(T^\ast \operatorname{Gr}_{k,n}),
\end{gather*}
where $\tau_w$ is the automorphism of $H_{T_n}(\cdot)$ that implements the natural action of $W$ on it; explicitly, $(f\tau_w)(\hbar,x_1,\ldots,x_n)=f(\hbar,x_{w(1)},\ldots,x_{w(n)})$.

If we fix a particular basis of $H_{T_n}(T^\ast \operatorname{Gr}_{k,n})$, we can write $\mathcal R_w = \tau_w\check R_w$, where $\check R_w$ is the $H_{T_n}(\cdot)$-valued matrix of the action of $w$ in that particular basis. We choose the basis of the $[\CS_{\cl(\pi)}]$; it means that we define $\check R_w$ by
\begin{gather}\label{eq:defR}
\big[w^{-1} \CS_{\cl(\pi)}\big]= \sum_{\pi'\in \LP_{k,n}} (\check R_w)_{\pi,\pi'} [\CS_{\cl(\pi')}].
\end{gather}

The first result of this section, slightly generalizing one of the results of \cite{artic39}, is:
\begin{thm}\label{thm:coh}
In the case of the elementary transposition $(i,i+1)$, the following expression for $\check R_i \equiv \check R_{(i,i+1)}$ holds:
\begin{gather*}
\check R_i=1+\frac{x_i-x_{i+1}}{\hbar+x_{i+1}-x_i}e_i,
\end{gather*}
where $e_i$ is the matrix of the generator of $\TL_n(2)$ acting on $\mathcal H_{k,n}(2)^\ast$ equipped with its canonical basis.
\end{thm}

Since $w\mapsto \mathcal R_w$ forms a representation of the Weyl group $W\cong\mathcal S_n$, any $\check R_w$ can be written in terms of the $\check R_i$, $i=1,\ldots,n-1$. Furthermore, by writing the various Coxeter relations of $\mathcal S_n$, we immediately find identities satisfied by the $\check R_i$;
\begin{gather*}
\check R_i(x_{i+1}-x_{i+2})\check R_{i+1}(x_i-x_{i+2})\check R_i(x_i-x_{i+1}) \\
\qquad{} =\check R_{i+1}(x_i-x_{i+1})\check R_i(x_i-x_{i+2})\check R_{i+1}(x_{i+1}-x_{i+2}),\\
\check R_i(x_{i+1}-x_i)\check R_i(x_i-x_{i+1})=1, \\
\check R_i(x_i-x_{i+1})\check R_j(x_j-x_{j+1})=\check R_j(x_j-x_{j+1})\check R_i(x_i-x_{i+1}),\qquad |i-j|>1,
\end{gather*}
where we used the notation $\check R_i(u):=1+\frac{u}{\hbar-u}e_i$ to facilitate the substitution of variables due to the action of the $\tau_w$. The first two equations are known as {\em Yang--Baxter equation} and {\em unitarity equation}, respectively.

\subsubsection{Pushforward}
The second result of this section concerns the pushforward of these classes. There are three natural choices of pushforward~-- using $\mu\colon T^\ast \operatorname{Gr}_{k,n}\to \operatorname{End}(\CC^n)$; using $\mu'\colon \mu^{-1}(\mathfrak{n}_-)\to \mathfrak{n}_-$, which takes into account that all $\CS_I$ lie in $\mu^{-1}(\mathfrak{n}_-)$; and using $\pi\colon T^\ast \operatorname{Gr}_{k,n}\to \{\cdot\}$. All the target spaces are equivariantly contractible, so their cohomology is that of a point, that is, $\QQ(x_1,\ldots,x_n,\hbar)$; and we have the simple relations
\begin{gather*}
\mu_\ast(v) =\prod_{i\ge j}(\hbar+x_i-x_j) \mu'_\ast(v), \\
\pi_\ast(v) =\prod_{i<j}(\hbar+x_i-x_j)^{-1} \mu'_\ast(v)
\end{gather*}
for all $v\in H_{T_n}\big(\mu^{-1}(\mathfrak n_-)\big)$. From these formulae it is obvious that $\mu_\ast$ has the disadvantage that it introduces a common factor, whereas~$\pi_\ast$, being non proper, introduces a~denominator. We shall therefore use~$\mu'_\ast$; however we shall see that $\mu'_\ast$ has the disadvantage of changing the normalization of the $R$-matrix.

We now consider
\begin{gather*}
\Psi_\pi:=\mu'_\ast[\CS_{\cl(\pi)}].
\end{gather*}
Because of the fact that $\mu(\CS_{\cl(\pi)})=\mathcal O_\pi=\overline{B_n\cdot \pi_>}$, we see that there are two possibilities. Either
\begin{itemize}\itemsep=0pt
\item $\pi\in\LP_{k,n}^{(m)}$, in which case $\mu$ is birational on $\CS_{\cl(\pi)}$, and $\Psi_\pi=[\mathcal O_\pi]$; or
\item $\pi\not\in\LP_{k,n}^{(m)}$, in which case $\dim\mathcal O_\pi<\dim \CS_{\cl(\pi)}$ and $\Psi_\pi=0$.
\end{itemize}

Applying $\mu'_\ast$ to \eqref{eq:defR} with $w=(i,i+1)$, and paying attention to the fact that the $W$-action on $H_{T_n}(\mathfrak n_-)$ is the natural action $\tau_w$ conjugated by multiplication by $[0]_{\mathfrak n_-}$, we find that
\begin{gather}\label{eq:exchcoh}
\Psi_\pi \tau_i = \sum_{\pi'\in \LP^{(m)}_{k,n}} (\check R'_i)_{\pi,\pi'} \Psi_{\pi'},\qquad \pi\in \LP^{(m)}_{k,n}
\end{gather}
with
\begin{gather*}
\check R'_i=\frac{\hbar+x_{i+1}-x_i}{\hbar+x_i-x_{i+1}}\check R_i =\frac{\hbar+x_{i+1}-x_i+(x_i-x_{i+1})e_i}{\hbar+x_i-x_{i+1}}.
\end{gather*}
\eqref{eq:exchcoh} is known as the {\em exchange relation}. The crucial difference with~\eqref{eq:defR} is that it involves only link patterns with maximal number of pairings. In other words, $\Psi:=\sum_{\pi\in \LP^{(m)}_{k,n}} \Psi_\pi \ket{\pi}$ lives in $\mathcal H_{k,n}^{(\ge m)}(2) \otimes \QQ(x_1,\ldots,x_n,\hbar)$.

\subsection[$K$-theory 1: noncrossing loops]{$\boldsymbol{K}$-theory 1: noncrossing loops}\label{sec:K1}
In the previous section, we have obtained a solution of the Yang--Baxter equation associated to the Temperley--Lieb algebra $\TL_n$ where the parameter $\beta$ is set to the value $2$. It is natural to try to obtain a~solution for arbitrary values of $\beta$ by extending this construction to $K$-theory. Indeed, this is possible, by following the general philosophy of~\cite{artic69} (see also~\cite{RTV-K} for a different approach, leading to a distinct, but related, basis). As mentioned above, it is convenient in this section to allow oneself to use the square root of~$t$ by considering the double cover~$T'_n$ of~$T_n$; in particular, we shall identify $K_{T'_n}(T^\ast \operatorname{Gr}_{k,n})$ with the $\TL_n\big(t^{1/2}+t^{-1/2}\big)$-right module $\mathcal H_{k,n}\big(t^{1/2}+t^{-1/2}\big)^\ast\otimes \QQ\big(z_1,\ldots,z_n,t^{1/2}\big)$.

\subsubsection{The coherent sheaves}\label{sec:coh}
We now define certain $T'_n$-equivariant coherent sheaves on $T^\ast \operatorname{Gr}_{k,n}$ as follows. Given a noncrossing link pattern $\pi\in\LP_{k,n}$, we first define a sheaf on $S_{\cl(\pi)}$: the latter being Cohen--Macaulay~\cite{Ra-CM} (see also~\cite{Brion-flag}), it possesses a dualizing/canonical sheaf $\omega_{S_{\cl(\pi)}}$; tensor it with $O(n-k)$ (pulled back from $\operatorname{Gr}_{k,n}$), defining
\begin{gather}\label{eq:defsigma}
\underline{\sigma}_\pi := \omega_{S_{\cl(\pi)}} \otimes O(n-k).
\end{gather}
We then take the inverse image of $\underline{\sigma}_\pi$ under the map $p\colon \CS_{\cl(\pi)}\to S_{\cl(\pi)}$ and its direct image to $T^\ast \operatorname{Gr}_{k,n}$, thus resulting in a certain coherent sheaf $\sigma_\pi$ supported on $\CS_{\cl(\pi)}$.

$\sigma_\pi$ is naturally $T_n$- (or $T'_n$-) equivariant; however, there is the freedom to tensor by a nonequivariantly trivial line bundle carrying a representation of~$T'_n$. We take care of this freedom by fixing the weight of $\sigma_\pi$ at the fixed point $\CC^{\cl(\pi)}$:
\begin{gather}\label{eq:weinorm}
\wt_K\big(\sigma_\pi|_{\CC^{\cl(\pi)}}\big) = m_\pi:=t^{\frac{1}{2}\sum\limits_{i\in \cl(\pi)} i}\prod_{i\in\cl(\pi)} z_i^{-\#(\overline{\cl(\pi)}\cap [1,i])} \prod_{i\not\in\cl(\pi)} z_i^{-\#(\cl(\pi)\cap [1,i])}.
\end{gather}

$\underline{\sigma}_\pi$ and $\sigma_\pi$ can be more explicitly defined (nonequivariantly) in terms of certain combinatorial data attached to $\pi$. Consider
\begin{gather*}
d(\pi):=\{i\in \{1,\ldots,n-1\}\colon i\in \cl(\pi),\, i+1\not\in\cl(\pi)\}.
\end{gather*}
Geometrically, the $\{ S_{\cl(e_i \pi)},\ i\in d(\pi) \}$ (where $\cl(e_i\pi)$ is obtained from $\cl(\pi)$ by replacing $i$ with $i+1$) are exactly the Schubert varieties of codimension $1$ inside $S_{\cl(\pi)}$, and form a basis of the (Weyl) divisor class group of~$S_{\cl(\pi)}$. We then consider the divisor
\begin{gather*}
\underline{D}(\pi):=\sum_{i\in d(\pi)} a_i(\pi)S_{\cl(e_i\pi)},\qquad a_i(\pi)=i-2\#(\cl(\pi)\cap[1,i]).
\end{gather*}
Associated to it is $\underline{\sigma}_\pi$, as the sheaf of functions on $S_{\cl(\pi)}$ with poles of order at most $a_i(\pi)$ on~$S_{\cl(e_i\pi)}$ (or zeroes of order at least $-a_i(\pi)$ if $a_i(\pi)$ is negative). In general, $\underline{\sigma}_\pi$ is reflexive but not invertible, because the canonical divisor is not Cartier~\cite{WY-SchubGroth} (more on this in Section~\ref{sec:goren}).

Similarly, $\sigma_\pi$ is related to a divisor $D(\pi)$ of a $\CS_{\cl(\pi)}$; the latter is slightly more involved to define, and will force us to delve deeper into the combinatorics of link patterns. Given $\pi\in \LP_{k,n}$, the {\em depth} of an arc of $\pi$ is the number of arcs needed to escape to left infinity (starting right outside the arc), minus the number of vertices paired to left infinity to its left. The depth is just a graphical reformulation of the coefficients $a_i(\pi)$, in the sense that if $i\in d(\pi)$, then $a_i(\pi)$ is the depth of the arc starting at~$i$. We alternatively denote by $a_\alpha(\pi)$ the depth of an arc~$\alpha$.

Two arcs are {\em neighboring} if they have same depth, and they have a border region in common (i.e., they are in the closure of the same connected component of the complement of all the arcs). Note that (neighboring or equality) is an equivalence relation.

Given two neighboring arcs $\alpha$ and $\beta$, there are exactly three ways one can reconnect the endpoints of~$\alpha$ and~$\beta$, leaving the other arcs untouched:
\begin{itemize}\itemsep=0pt
\item The original link pattern $\pi=\tikz[baseline=-3pt]{\linkpattern[alias=false,numbering={1/\ldots,2/\qquad\alpha,3/,4/\ldots,5/\qquad\beta,6/,7/\ldots}]{2/3,5/6}}$.
\item Another noncrossing link pattern, denoted $e_{\alpha,\beta}\pi:=\tikz[baseline=-3pt]{\linkpattern[alias=false,numbering={1/\ldots,2/,3/,4/\ldots,5/,6/,7/\ldots}]{2/6,3/5}}$.
\item A crossing link pattern, denoted $f_{\alpha,\beta}\pi:=\tikz[baseline=-3pt]{\linkpattern[alias=false,numbering={1/\ldots,2/,3/,4/\ldots,5/,6/,7/\ldots}]{2/5,3/6}}$.
\end{itemize}
With a bit of foresight, let us denote
\begin{gather*}
X_{f_{\alpha,\beta}\pi}:= \CS_{\cl(e_{\alpha,\beta}\pi)}\cap \CS_{\cl(\pi)}.
\end{gather*}

Then we have
\begin{gather}\label{eq:div}
D(\pi)=\sum_{\substack{\alpha,\beta\\\text{neighboring arcs of $\pi$}}} a_\alpha(\pi) X_{f_{\alpha,\beta}\pi}.
\end{gather}

\begin{ex}
Let
\begin{gather*}
\pi=\tikz[baseline=-3pt,x=\linkpatternunit,y=\linkpatternunit]{
\draw[brown,latex-latex] (2,-0.8) -- (3.7,-0.8);
\draw[brown,latex-latex,bend right] (4.6,-0.3) to (6.4,-0.3);
\draw[brown,latex-latex,bend right] (6.6,-0.3) to (8.4,-0.3);
\draw[brown,latex-latex,bend right] (4.3,-0.3) to (8.7,-0.3);
\linkpattern[looseness=0.3]{1/1/\toleftx{2},2/2/\toleftx{4},3/10,4/5,6/7,8/9}
\node at (0.5,-0.4) {$\ss 0$};
\node at (1.5,-1) {$\ss -1$};
\node at (3,-1.5) {$\ss -2$};
\node at (6.5,-1.4) {$\ss -1$};
}
\end{gather*}
where we indicated neighboring arcs and depths on the diagram. Then
\begin{gather*}
\underline{D}(\pi)=-2S_{\cl(\tikz[baseline=-3pt]{\linkpattern[small]{1/1/\toleft,2/3,4/5,6/7,8/9,10/10/\toleft}})}-S_{\cl(\tikz[baseline=-3pt]{\linkpattern[small]{1/1/\toleft,2/2/\toleft,3/10,4/7,5/6,8/9}})}
-S_{\cl(\tikz[baseline=-3pt]{\linkpattern[small]{1/1/\toleft,2/2/\toleft,3/10,4/5,6/9,7/8}})},
\\[-2mm]
D(\pi)=-2X_{\linkpattern[small]{1/1/\toleft,2/10,3/3/\toleft,4/5,6/7,8/9}}-X_{\linkpattern[small]{1/1/\toleft,2/2/\toleft,3/10,4/6,5/7,8/9}}\\[-2mm]
\hphantom{D(\pi)=}{} -X_{\linkpattern[small]{1/1/\toleft,2/2/\toleft,3/10,4/5,6/8,7/9}} -X_{\linkpattern[small]{1/1/\toleft,2/2/\toleft,3/10,4/8,5/9,6/7}}.
\end{gather*}
Note that $D(\pi)$ has one more term than $\underline{D}(\pi)$.
\end{ex}

\subsubsection[The geometric $R$-matrix]{The geometric $\boldsymbol{R}$-matrix}
We can now follow the same construction as in cohomology, defining the geometric $R$-matrix to be given by the action of the Weyl group on $K_{T'_n}(T^\ast \operatorname{Gr}_{k,n})$ in the basis of the $[\sigma_\pi]$:
\begin{gather*}
\big[w^{-1} \sigma_\pi\big]= \sum_{\pi'\in \LP_{k,n}} \big(\check R^{nc}_w\big)_{\pi,\pi'} [\sigma_{\pi'}].
\end{gather*}

In the case of $w$ the elementary transposition $(i,i+1$), we find:
\begin{thm}\label{thm:nc}
\begin{gather*}
\check R^{nc}_i=1+\frac{z_{i+1}-z_i}{t^{-1/2}z_i-t^{1/2}z_{i+1}}e_i,
\end{gather*}
where $e_i$ is the generator of $\TL_n\big(t^{-1/2}+t^{1/2}\big)$.
\end{thm}

\begin{ex}\label{ex:base} Let us consider the simplest nontrivial case, which is $k=1$, $n=2$. There are two noncrossing link patterns, \linkpattern{1/1/\toleft,2/2} and \linkpattern{1/2}, corresponding to the subsets~$\{1\}$ and $\{2\}$, respectively, to Schubert varieties which are the whole of $\PP^1$ and the point $\CC^{\{2\}}$, and to conormal Schubert varieties which are the base $\PP^1$ and the fiber at~$\CC^{\{2\}}$. We compute
$d(\linkpattern[small]{1/1/\toleft,2/2})=\{1\}$, $d(\linkpattern[small]{1/2})=\varnothing$, so that $\underline{D}(\linkpattern[small]{1/1/\toleft,2/2})=-[\text{point}]$, $\underline{D}(\linkpattern[small]{1/2})=0$. The sheaf associated to the former is nothing but $O(-1)$. The classes of the conormal varieties restricted to fixed points are given by
\begin{gather*}
\begin{matrix}
&\{1\}&\{2\}
\\
[\CS_{\linkpattern[small]{1/1/\toleft,2/2}}]=(&\!\!\!\!\!1-t z_2/z_1&1-t z_1/z_2\!\!\!\!\!&),\\
[\CS_{\linkpattern[small]{1/2}}]=(&\!\!\!\!\! 0&1-z_2/z_1\!\!\!\!\! &).
\end{matrix}
\end{gather*}
Considering $O(-1)$ on $\PP^1$, choosing its weight at $\{1\}$ to be $t^{1/2}z_2^{-1}$ (according to~\eqref{eq:weinorm}), and therefore $t^{1/2}z_1^{-1}$ at $\{2\}$, and similarly the weight $t z_2^{-1}$ for the trivial bundle over the fiber at~$\CC^{\{2\}}$, we find
\begin{gather*}
\begin{matrix}
&\{1\}&\{2\} \\
[\sigma_{\linkpattern[small]{1/1/\toleft,2/2}}]=\big( &\!\!\!\!\! t^{1/2}\big(z_2^{-1}-t z_1^{-1}\big)&t^{1/2}\big(z_1^{-1}-t z_2^{-1}\big)\!\!\!\!\!&\big),\\
[\sigma_{\linkpattern[small]{1/2}}]=\big( &\!\!\!\!\! 0&t\big(z_2^{-1}-z_1^{-1}\big)\!\!\!\!\! &\big).
\end{matrix}
\end{gather*}
$\PP^1$ and $\sigma_{\linkpattern[small]{1/1/\toleft,2/2}}$ are obviously invariant by the Weyl group action; as for $w=(1,2)$,
\begin{gather*}
[w^{-1}\sigma_{\linkpattern[small]{1/2}}] =
\big( \ t\big(z_1^{-1}-z_2^{-1}\big)\qquad 0\ \big)
%\\ \hphantom{[w^{-1}\sigma_{\linkpattern[small]{1/2}}]}{}
 =-\frac{z_2-z_1}{t^{1/2}z_2-t^{-1/2}z_1}[\sigma_{\linkpattern[small]{1/1/\toleft,2/2}}]
+\frac{t z_1-z_2}{t z_2-z_1} [\sigma_{\linkpattern[small]{1/2}}].
\end{gather*}
This implies that the $R$-matrix has the form
\begin{gather*}
\check R^{nc}_1 =
\begin{pmatrix}
1&0\\
\dfrac{z_2-z_1}{t^{-1/2}z_1-t^{1/2}z_2}&\dfrac{z_2-t z_1}{z_1-tz_2},
\end{pmatrix}
\end{gather*}
which matches with the expression of Claim~\ref{thm:nc}, with $e_1=\begin{pmatrix} 0 & 0 \\ 1 & t^{1/2}+t^{-1/2}
\end{pmatrix}$.
\end{ex}

\subsubsection{Pushforward}
Again one can now pushforward the $[\sigma_\pi]$ to a point, or equivalently to $\mathfrak n_-$. The highly nontrivial facts are that
\begin{itemize}\itemsep=0pt
\item As in cohomology, $\mu'_\ast [\sigma_\pi]\ne 0$ iff $\pi\in\LP_{k,n}^{(m)}$ (i.e., iff $\pi$ has maximal number of pairings).
\item For $k\le n/2$, the $[\sigma_\pi]$ have no higher sheaf cohomology, so that $\mu'_\ast [\sigma_\pi]$ is the class of the direct image $\mu'_\ast \sigma_\pi$, which is the space of global sections of $\sigma_\pi$ viewed as a module over the coordinate ring of $\mathfrak n_-$ (equivalently, $\pi_\ast [\sigma_\pi]$ is the character of the space of global sections).\footnote{A similar statement can be made if $k\ge n/2$ by replacing $O(n-k)$ with $O(k)$ in the definition \eqref{eq:defsigma}, and modifying the weight \eqref{eq:weinorm} appropriately; see also \cite{artic69} for a more symmetric choice.} In fact, we can describe $\mu'_\ast \sigma_\pi$ explicitly for $\pi\in\LP_{k,n}^{(m)}$: it is the module of functions on $\mathcal O_\pi$ with poles of order at most $a_{\alpha}(\pi)$ on $\mathcal O_{f_{\alpha,\beta}\pi} $ for each pair of neighboring arcs~$\alpha$,~$\beta$.
\end{itemize}

For now, we only use the first property to find as a corollary the analogue of~\eqref{eq:exchcoh}, namely, defining $\Psi^{nc}_\pi=\mu'_\ast[\sigma_\pi]$,
\begin{gather}\label{eq:exchnc}
\Psi^{nc}_\pi \tau_i = \sum_{\pi'\in \LP^{(m)}_{k,n}} (\check R'_i{}^{nc})_{\pi,\pi'} \Psi^{nc}_{\pi'},\qquad \pi\in \LP^{(m)}_{k,n}
\end{gather}
with
\begin{gather*}
\check R'_i{}^{nc}=\frac{1-t\,z_{i+1}/z_i-t^{1/2}(1-z_{i+1}/z_i)e_i}{1-t\,z_i/z_{i+1}}.
\end{gather*}

As beautiful as this construction may be, it does not answer a more geometrically natural question, which is to understand the action of the Weyl group on {\em structure sheaves} of the $\CS_I$. Naive attempts at implementing the same procedure for such sheaves are a failure, since the corresponding $R$-matrices are highly nonlocal and depend on all spectral parameters. We now provide an appropriate modification of this procedure, as advertised in the introduction.

\subsection[$K$-theory 2: crossing loops]{$\boldsymbol{K}$-theory 2: crossing loops} \label{sec:K2}
The idea is to introduce a generating set of $K_{T_n}(T^\ast \operatorname{Gr}_{k,n})$ made of classes of structure sheaves of certain subvarieties (including conormal Schubert varieties), with the hope that the Weyl group action on these subvarieties will be ``nice''. In general, these classes will be linearly dependent (over $\QQ(z_1,\ldots,z_n,t)$), so that $K_{T_n}(T^\ast \operatorname{Gr}_{k,n})$ will be identified with a {\em quotient} of~$\widetilde{\mathcal H}_{k,n}(\beta_0)^\ast \otimes \QQ(z_1,\ldots,z_n,t)$ (for a~$\beta_0$ to be defined below); or equivalently, we shall find a subspace of $\widetilde{\mathcal H}_{k,n}(\beta_0)\otimes \QQ(z_1,\ldots,z_n,t)$ that is stable under
the Weyl group action $\tau_w \check R^c_w$, for some appropriate $\check R^c_w\in \CTL_n\otimes \QQ(z_1,\ldots,z_n,t)$.

The main ingredient of this section is the choice of these subvarieties. We denote them $X_\pi$, where $\pi$ runs over $\CLP_{k,n}$; their definition is deceptively simple:
\begin{gather*}
X_\pi = p^{-1}(S_{\cl(\pi)}) \cap \overline{\mu^{-1}(B\cdot\pi_<)} = \CS_{\cl(\pi)} \cap \overline{\mu^{-1}(B\cdot\pi_<)}.
\end{gather*}
These subvarieties will be studied in \cite{prep1}. Here we give certain properties of $X_\pi$.

Firstly, $X_\pi$ is irreducible, hence a subvariety. In the special case that $\pi\in \LP_{k,n}$ ($\pi$ is noncrossing), $\CS_{\cl(\pi)}=\overline{\mu^{-1}(B\cdot\pi_<)}$, so that $X_\pi = \CS_{\cl(\pi)}$. If $\pi$ has one crossing, then it is of the form $\pi=f_{\alpha,\beta}\pi'$, where $\alpha$ and $\beta$ are neighboring arcs in the sense of Section~\ref{sec:coh}, and we have $X_\pi=\CS_{\cl(\pi')}\cap \CS_{\cl(e_{\alpha,\beta}\pi')}$, coinciding with the definition given in Section~\ref{sec:coh}. In general, one has
\begin{gather*}
\dim X_\pi = k(n-k)-\cross(\pi)
\end{gather*}
(which shows in particular that if $\pi\not\in \LP_{k,n}$, $X_\pi$ is not Lagrangian~-- only isotropic).

\subsubsection[The geometric $R$-matrix]{The geometric $\boldsymbol{R}$-matrix}
We now proceed analogously to the previous sections. We look for $K_{T_n}(\cdot)$-valued matrices $\check R^c_w$, $w\in W$, satisfying
\begin{gather}\label{eq:defRcr}
[w^{-1} X_\pi]= \sum_{\pi'\in \CLP_{k,n}} (\check R^c_w)_{\pi,\pi'} [X_{\pi'}]
\end{gather}
for every $\pi\in \CLP_{k,n}$.

Since the $[X_\pi]$ are in general linearly dependent, \eqref{eq:defRcr} is a highly overdetermined system of equations, and it is somewhat miraculous that it admits a solution; in the case that $w$ is the elementary transposition $(i,i+1)$, we find
\begin{thm}[\cite{prep1}]\label{thm:cr}
\begin{gather*}
\check R^{c}_i=1+\frac{1-z_i/z_{i+1}}{1-t z_{i+1}/z_i}e_i+(1-z_i/z_{i+1})f_i,
\end{gather*}
where the $e_i$, $f_i$ are generators of $\CTL_n(1+t)$.
\end{thm}

The $\check R^c_i$ satisfy the Yang--Baxter equation, unitarity equation and commutation far apart as operators on $K_{T_n}(T^\ast \operatorname{Gr}_{k,n})$~-- in fact, they satisfy these as abstract elements of~$\CTL_n\otimes \QQ(z_1,\ldots,z_n,t)$, as a consequence of the defining relations of the algebra. As far as the author knows, this is a new solution of the Yang--Baxter equation.

Introduce the notations $a(z)=1-t/z$ and $b(z)=1-z$. It is useful to write more explicitly~\eqref{eq:defRcr} according to the following trichotomy:
\begin{itemize}\itemsep=0pt
\item $i$ and $i+1$ are not paired together, and the arcs coming out of $i$ and $i+1$ do not cross. Then
\begin{gather}\label{eq:defRcra}
[(i,i+1)X_\pi]=[X_\pi]
\end{gather}
expressing the invariance of $X_\pi$ under $(i,i+1)$ (and in fact, under the whole of the sub\-group~$\operatorname{GL}_2^{(i)}$ of~$\operatorname{GL}_n$ which differs from the identity only in rows and columns $i$, $i+1$).

\item $i$ and $i+1$ are not paired together, and the arcs coming out of $i$ and $i+1$ cross. Then there exists a unique $\rho\ne\pi$ such that $f_i\rho=\pi$ (obtained by ``uncrossing'' the arcs coming out of $i$ and $i+1$), and we have
\begin{gather}\label{eq:defRcrb}
[(i,i+1)X_\pi]=(1-b(z_i/z_{i+1}))[X_\pi]+b(z_i/z_{i+1})[X_\rho]
\end{gather}
or
\begin{gather*}
D_i [X_\pi] =[X_\rho],
\end{gather*}
where $D_i$ is the divided difference Demazure operator $D_i=\frac{z_{i+1} (i,i+1)-z_i}{z_{i+1}-z_i}$, expressing the fact that the map from $\operatorname{GL}_2^{(i)} \times_{B_2^{(i)}} X_\pi$ to $X_\rho$ ($B^{(i)}_2 = B_n \cap \operatorname{GL}_2^{(i)}$) is generically one-to-one.

\item $i$ and $i+1$ form a pairing. Then
\begin{gather}\label{eq:defRcrc}
a(z_i/z_{i+1})[(i,i+1)X_\pi]=a(z_{i+1}/z_i)[X_\pi]+b(z_i/z_{i+1})\sum_{\substack{\rho\ne\pi\\ e_i\rho=\pi}} (-1)^{\cross(\rho)-\cross(\pi)} [X_\rho]
\end{gather}
or
\begin{gather*}
D_i(a(z_{i+1}/z_i)[X_\pi]) = a(z_{i+1}/z_i)[X_\pi] + \sum_{\substack{\rho\ne\pi\\ e_i\rho=\pi}} (-1)^{\cross(\rho)-\cross(\pi)} [X_\rho].
\end{gather*}
This equation is more complicated to explain, and we shall not go into the details. Let us simply note that the l.h.s.\ corresponds to geometrically to first ``cutting'' with the equation $u_{i+1,i}=0$, then ``sweeping'' with $\operatorname{GL}_2^{(i)}$ as in the previous case. The r.h.s.\ corresponds to a careful analysis of the resulting scheme and its components (see \cite[Section~5.4]{artic39} for a~similar discussion, but in cohomology only).
\end{itemize}

\begin{rmk}%{rmk*}
One could take the limit from $K$-theory to cohomology of the results of this section and obtain a ``rational'' $R$-matrix for the varieties $X_\pi$, which generalizes the one obtained in Section~\ref{sec:Rmatcoh} to the crossing case. We would recover this way the results of \cite[Section~5.5]{artic39}. Since the focus of the present article is on $K$-theory and for the sake of compactness, we do not study such a limit here.
\end{rmk}%{rmk*}

\subsubsection{Pushforward}
\looseness=-1 Next, we discuss the pushforward to $\mathfrak n_-$. Define $\Psi^c_\pi=\mu'_\ast[X_\pi]$. The $[X_\pi]$ being classes of structure sheaves, it is of course not true that some of them are sent to zero by $\mu'_\ast$, as in the previous sections. However, note the following property: first, \eqref{eq:defRcra} and \eqref{eq:defRcrb} only involve link patterns with the same number of pairings. Only remains~\eqref{eq:defRcrc}. We rewrite it slightly, by noting that proper preimages of $\pi$ under $e_i$ always come in pairs, related to each other by~$f_i$, so that we have
\begin{gather}
a(z_i/z_{i+1})[(i,i+1)X_\pi]\nonumber\\
\qquad{} =a(z_{i+1}/z_i)[X_\pi]+b(z_i/z_{i+1})\sum_{\substack{\rho\ne\pi\\ e_i\rho=\pi\\f_i\rho\ne -\rho}} (-1)^{\cross(\rho)-\cross(\pi)} ([X_\rho]-[X_{f_i\rho}]),\label{eq:defRcrcc}
\end{gather}
where the summation is now only over the $\rho$ such that the arcs coming from $i$, $i+1$ do not cross. The $e_i$ action may change the number of pairings only in one case, and that is when neither~$i$ not $i+1$ are paired. In this case $\rho$ (resp.~$f_i\rho$) is identical to $\pi$ except $i$ is connected to left infinity (resp.\ bottom infinity), and $i+1$ is connected to bottom infinity (resp.\ left infinity).
The result is that $\rho_>=(f_i\rho)_>$ (since these matrices only care about pairings), and therefore $\mu(X_\rho)=\mu(X_{f_i\rho})$. This means that if we apply $\mu'_\ast$
to \eqref{eq:defRcrcc}, their contribution compensates. Finally, we are led to the following statement:
\begin{gather}\label{eq:exchcr}
\Psi^c_\pi \tau_i = \sum_{\pi'\in \CLP^{(r)}_{k,n}} (\check R_i'{}^c)_{\pi,\pi'} \Psi^c_{\pi'},\qquad \pi\in \CLP^{(r)}_{k,n},\qquad 0\le r\le m,
\end{gather}
where
\begin{gather*}
\check R_i'{}^c =
\frac{a(z_i/z_{i+1})+b(z_i/z_{i+1})e_i+a(z_i/z_{i+1})b(z_i/z_{i+1})f_i}{a(z_{i+1}/z_i)}.
\end{gather*}
The difference with \eqref{eq:defRcr} is that the summation in~\eqref{eq:exchcr} is over link patterns with fixed number of pairings.

\begin{ex}\label{ex:n4k2}
$n=4$, $k=2$, \eqref{eq:defRcrc} with $(i,i+1)=(1,2)$ and $\pi=\linkpattern{1/2,3/4}$. Proper preimages of $\pi$ under $e_1$ are \linkpattern{1/4,2/3}, \linkpattern{1/3,2/4}, \linkpattern{1/1/\toleft,2/2,3/4}, \linkpattern{1/1,2/2/\toleft,3/4}.

The corresponding varieties being smooth, we can easily compute the restrictions of their classes to fixed points; writing $a_{ij}=a(z_i/z_j)$, $b_{ij}=b(z_i/z_j)$, we have
\begin{gather*}
\begin{matrix}
&\!\!\!\!\!\{1,2\}&\{1,3\}&\{1,4\}&\{2,3\}&\{2,4\}&\{3,4\}\\
[X_{\linkpattern[small]{1/2,3/4}}]=
(
&\!\!\!\!\!0
&0
&0
&0
&a_{23}b_{21}b_{41}b_{43}
&a_{32}b_{31}b_{41}b_{42}
\!\!\!\!\!&),
\\
[X_{\linkpattern[small]{1/4,2/3}}]=
(
&\!\!\!\!\!0
&0
&0
&0
&0
&b_{32}b_{31}b_{41}b_{42}
\!\!\!\!\!&),
\\
[X_{\linkpattern[small]{1/3,2/4}}]=
(
&\!\!\!\!\!0
&0
&0
&0
&0
&a_{32}b_{32}b_{31}b_{41}b_{42}
\!\!\!\!\! &),
\\
[X_{\linkpattern[small]{1/1/\toleft,2/2,3/4}}]=
(
&\!\!\!\!\!0
&0
&a_{12}a_{13}b_{42}b_{43}\ \
&0
&a_{21}a_{23}b_{41}b_{43}
&a_{31}a_{32}b_{41}b_{42}
\!\!\!\!\! &),
\\
[X_{\linkpattern[small]{1/1,2/2/\toleft,3/4}}]=
(
&\!\!\!\!\!0
&0
&0
&0
&a_{21}b_{21}a_{23}b_{41}b_{43}\ \
&a_{31}a_{32}b_{31}b_{41}b_{42}
\!\!\!\!\! &),
\end{matrix}
\end{gather*}
and one can check that \eqref{eq:defRcrc} holds. Pushing forward using $\mu'$ gives
\begin{gather*}
\mu'(X_{\linkpattern[small]{1/2,3/4}})=\left\{
\begin{aligned}&\begin{pmatrix}\\u_{2,1}\\u_{3,1}&0\\\star&u_{4,2}&u_{4,3}&\end{pmatrix}\\&u_{4,2}u_{2,1}+u_{4,3}u_{3,1}=0\end{aligned}
\right\},\\
\mu'_\ast [X_{\linkpattern[small]{1/2,3/4}}] =(1-t z_2/z_3)\big(1-t^2z_1/z_4\big),\\
\mu'(X_{\linkpattern[small]{1/4,2/3}}) =\left\{\begin{pmatrix}\\0\\\star&\star\\\star&\star&0&\end{pmatrix}\right\},\\
\mu'_\ast [X_{\linkpattern[small]{1/4,2/3}}]=(1-tz_1/z_2)(1-t z_3/z_4),\\
\mu'(X_{\linkpattern[small]{1/3,2/4}})=\left\{\begin{pmatrix}\\0\\\star&0\\\star&\star&0&\end{pmatrix}\right\},\\
\mu'_\ast [X_{\linkpattern[small]{1/3,2/4}}]=(1-tz_1/z_2)(1-tz_2/z_3)(1-tz_3/z_4),\\
\mu'(X_{\linkpattern[small]{1/1/\toleft,2/2,3/4}}) =\left\{\begin{pmatrix}\\0\\0&0\\\star&\star&\star&\end{pmatrix}\right\},\\
\mu'_\ast [X_{\linkpattern[small]{1/1/\toleft,2/2,3/4}}]=(1-tz_1/z_2)(1-tz_1/z_3)(1-tz_2/z_3),\\
\mu'(X_{\linkpattern[small]{1/1,2/2/\toleft,3/4}})=\left\{\begin{pmatrix}\\0\\0&0\\\star&\star&\star&\end{pmatrix}\right\},\\
\mu'_\ast [X_{\linkpattern[small]{1/1,2/2/\toleft,3/4}}]=(1-tz_1/z_2)(1-tz_1/z_3)(1-tz_2/z_3)
\end{gather*}
(where $\star$ indicates a free entry) and the last two compensate, resulting in the simpler identity (cf.~\eqref{eq:exchcr})
\begin{gather*}
(1-t\,z_1/z_2)(D_1-1)\mu'_\ast [X_{\linkpattern[small]{1/2,3/4}}]=\mu'_\ast [X_{\linkpattern[small]{1/4,2/3}}]-
\mu'_\ast [X_{\linkpattern[small]{1/3,2/4}}].
\end{gather*}
\end{ex}

\section{The partition function of finite domains}\label{sec:fd}
This section is largely based on \cite{artictt}.

\subsection{Full crossing link patterns and their poset}\label{sec:poset}
Let $N=2K$ be an even integer,\footnote{We intentionally use upper case $K$ and $N$ in this section and the next, even though $K$ and $N$ will play similar roles as $k$ and $n$ in the rest of the paper; the reason will become clear in Section~\ref{sec:rect}.} and $\rho$ be an element of $\CLP^{(K)}_{K,N}$. Note that every vertex is paired by $\rho$, so we can equivalently think of $\rho$ as a {\em fixed-point free involution} of $\{1,\ldots,N\}$. Accordingly, for any element say $i\in \ZZ/N\ZZ$, we shall denote by $\rho(i)$ the vertex which is paired with $i$ in $\rho$.

\def\linkpatternedgecolor{colb}
In the whole of this section, it will be sometimes convenient to adopt a slightly different point of view,
by identifying $\{1,\ldots,N\}$ with $\ZZ/N\ZZ$, and accordingly, by redrawing $\rho$
as a crossing link pattern {\em on a circle} (with vertices ordered clockwise):
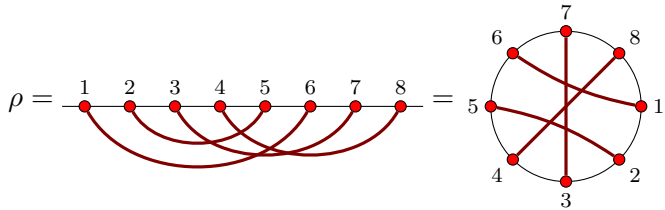
\begin{gather*}
\rho = \linkpattern{1/6,2/5,3/7,4/8} = \linkpattern[unit=1cm,shape=circle,squareness=0.2]{1/6,2/5,3/7,4/8}
\end{gather*}

As in \cite{artic33}, we write $\cyc{i_1\le\cdots\le i_k}$ iff there exist representatives of the $i_a\in\ZZ/N\ZZ$ in $\ZZ$ such that $i_1\le \cdots\le i_k<i_1+N$; any inequality $i_a\le i_{a+1}$ can be substituted with a $i_a<i_{a+1}$, with the obvious additional implication $i_a\ne i_{a+1}$.

There is a partial order on $\CLP_{K,N}^{(K)}$ which is defined as follows:\footnote{Actually, this order already appeared in \cite{artic39} in a non-cyclic formulation: it is exactly the order defined in \cite[Section~2]{artic39} {\em restricted to full link patterns}.} we say that $\rho\le\pi$ iff a~reduced diagram of $\pi$ can be obtained from a reduced diagram of $\rho$ via a sequence of moves
\begin{gather*}
\tikz[rotate=45,baseline=-3pt]{\plaqc} \quad \longrightarrow \quad \tikz[rotate=45,baseline=-3pt]{\plaqa}
\end{gather*}
(Note that there is freedom of rotating vertices, so that two different moves can be applied to any given vertex.)

We shall associate to $\rho$ an affine scheme $X_\rho$ living in an affine space $V_\rho$ of dimension $2 \cross(\rho)$ endowed with a symplectic structure, which is Lagrangian and invariant under a torus of dimension $K+1$, a subtorus of dimension $K$ of which preserves the symplectic structure.

\subsection{The ambient space}\label{sec:amb}
We first start with the complex vector space of square matrices with indices taking values in $\ZZ/N\ZZ$, which we denote $\Mat{\ZZ/N\ZZ}$. There is a natural projection to the vector space $W_\rho$, of dimension $4\cross(\rho)$, which corresponds to keeping only the following coordinates. Let $i,j\in\ZZ/N\ZZ$, such that $\cyc{i<j<\rho(i)<\rho(j)}$, i.e., the arcs $(i,\rho(i))$ and $(j,\rho(j))$ cross:
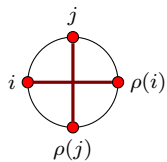
\begin{center}
\linkpattern[shape=circle,alias=false,numbering={1/{\rho(i)},2/{\rho(j)},3/i,4/j}]{1/3,2/4}
\end{center}
To each such crossing are associated {\em four} variables of $W_\rho$, namely, $M_{i,j}$, $M_{j,\rho(i)}$, $M_{\rho(i),\rho(j)}$, $M_{\rho(j),i}$. We call $\Omega$ the set of such coordinates.

Now define the subspace $V_\rho$ of $W_\rho$ by the following
equations: for any $\cyc{i<j<\rho(i)<\rho(j)}$, impose
\begin{gather}\label{eq:subst}
M_{i,j}+M_{\rho(i),\rho(j)}+\sum_{k:\cyc{j<k<\rho(i)},\cyc{\rho(j)<\rho(k)<i}} M_{i,k}M_{k,\rho(j)}=0.
\end{gather}

Clearly, \eqref{eq:subst} suggests that one should be able to express one of the two variables $M_{i,j}$, $M_{\rho(i),\rho(j)}$ in terms
of the other, and therefore to reduce to $2$ variables per crossing. This is the object of the following lemma:

\begin{lem}\label{lem:fourtotwo} Let $\Xi$ be a subset of cardinality $2\cross(\rho)$ of $\Omega$, such that $M_{i,j}\in \Xi$ iff $M_{\rho(i),\rho(j)}\not\in \Xi$. Then the natural projection $p_\Xi$ from $V_\rho$ to $\CC^{2\cross(\rho)}$ defined by keeping the variables in $\Xi$ is an isomorphism.
\end{lem}

\subsubsection{The torus action}\label{sec:torus}
We define a linear action of various tori on $\Mat{\ZZ/N\ZZ}$ by giving the additive weights of the variables $M_{i,j}$. We start with $\hat T_N=(\CC^\times)^N\times \CC^\times$ acting by
\begin{gather}\label{eq:T0-act}
\wt_H(M_{i,j})=\hbar+\sum_{\cyc{j<a\le i}} w_a,
\end{gather}
where $w_1,\ldots,w_N$ are coordinates on the Lie algebra of $(\CC^\times)^N$, and $\hbar$ is the coordinate on the Lie algebra of the last $\CC^\times$.

This linear action descends to one on $W_\rho$. $\hat T_N$ does not leave $V_\rho$ invariant; in fact, it is easy to see that to make~\eqref{eq:subst} homogeneous, the following relations must be imposed:
\begin{gather}\label{eq:restr}
\hbar+\sum_{\cyc{i<a\le \rho(i)}} w_a= 0,\qquad i=1,\ldots,N.
\end{gather}
This corresponds to a certain subgroup $T_\rho\subset \hat T_N$ of dimension $K$ which naturally acts on $V_\rho$.

In practice, the parameterization above of $T_\rho$ is inconvenient; by breaking the cyclic symmetry, one can introduce
\begin{gather*}
x_i=\sum_{a=1}^i w_a,\qquad i=1,\ldots,N,
\end{gather*}
with the relations $x_{\rho(i)}=x_i+\hbar$ if $1\le i<\rho(i)\le N$, and write instead
\begin{gather}\label{eq:torus}
\wt_H(M_{i,j})=\operatorname{sign}(j-i) \hbar+x_i-x_j, \qquad 1\le i,j\le N,\qquad M_{i,j}\in\Omega.
\end{gather}
The parameterization is redundant, in the sense that there are $K+1$ weights $\hbar,x_1,\ldots,x_K$, and only differences of $x_i$'s appear in~\eqref{eq:torus}; however, this subtlety can be safely ignored in what follows.

\subsubsection[$V_\rho$ as a slice of a nilpotent orbit]{$\boldsymbol{V_\rho}$ as a slice of a nilpotent orbit}\label{sec:slice}
We now connect $V_\rho$ with the geometry of Section~\ref{sec:geomsetup}, showing that $V_\rho$ is a ``slice'' of the nilpotent orbit closure~$\mathcal N_{K,N}$.

Unfortunately, due to irreconcilable conventions, we must relate $M$ to the {\em transpose} of the matrix of $u\in \mathcal N_{K,N}$, i.e., we shall map entries $M_{i,j}$ to entries $u_{j,i}$.

The slice $A_\rho$, which depends on the fixed-point-free involution $\rho$, is defined by a series of (inhomogeneous) linear equations. These equations are best described in terms of $2\times 2$ blocks of $M$ with rows $(i,\rho(i))$ and columns $(j,\rho(j))$, where $1\le i<\rho(i)\le N$, $1\le j<\rho(j)\le N$. Here are the possible relative configurations, described as sub-link patterns of $\rho$: (recall that $\star$ means that the entry satisfies no equation)
\begin{alignat}{3}\label{eq:stars}
&\linkpattern[alias=false,unit=1.2cm,numbering={1/i=j,2/\rho(i)=\rho(j)}]{1/2}\colon
\quad
\begin{pmatrix}0&\color{red}1\\\star&\star\end{pmatrix}; \qquad &&&
\\[-4mm] \notag
&\linkpattern[alias=false,numbering={1/i,2/\rho(i),3/j,4/\rho(j)}]{1/2,3/4}\colon
\quad
\begin{pmatrix}\color{red}0&\color{red}0\\\color{red}\star&\color{red}\star\end{pmatrix};
\qquad &&
\linkpattern[alias=false,numbering={3/i,4/\rho(i),1/j,2/\rho(j)}]{1/2,3/4}\colon
\quad
\begin{pmatrix}0&0\\\star&\star\end{pmatrix}; & \\\notag
&\linkpattern[alias=false,numbering={1/j,2/i,3/\rho(i),4/\rho(j)}]{1/4,2/3}\colon
\quad
\begin{pmatrix}0&\color{red}0\\\star&\color{red}\star\end{pmatrix};\qquad
&&
\linkpattern[alias=false,numbering={1/i,2/j,3/\rho(j),4/\rho(i)}]{1/4,2/3}\colon
\quad
\begin{pmatrix}\color{red}\star&\color{red}0\\\star&0\end{pmatrix}; &
\\\notag
&\linkpattern[alias=false,numbering={1/i,2/j,3/\rho(i),4/\rho(j)}]{1/3,2/4}\colon
\quad
\begin{pmatrix}\color{red}\star&\color{red}0\\\star&\color{red}\star\end{pmatrix};\qquad
&&
\linkpattern[alias=false,numbering={1/j,2/i,3/\rho(j),4/\rho(i)}]{1/3,2/4}\colon
\quad
\begin{pmatrix}0&\color{red}\star\\\star&\star\end{pmatrix}.&
\end{alignat}
For example, the first line of \eqref{eq:stars} means that for each pair $i<\rho(i)$ of $\rho$, one must impose the equation $M_{i,i}=0$ and $M_{i,\rho(i)}=1$. Similarly, the next diagram means that for each pair of non-crossing pairs $i<\rho(i)<j<\rho(j)$ one must impose $M_{i,j}=M_{i,\rho(j)}=0$; and so on.

The number of equations is $2\#\{\text{chords}\}+4\#\{\text{non-crossing pairs}\}+2\#\{\text{crossing pairs}\}= 2(K^2-\cross(\rho))$. We have marked in red entries that are in the strict upper triangle of~$M$; this allows to count the number of equations in the strict upper triangle, which is similarly found to be $K^2-\cross(\rho)$.

\begin{lem}\label{lem:slice} Denote $\Upsilon=\{M_{i,j}\colon i<j\text{ and the chords from }i, \ j\text{ cross in }\rho\}$. Consider the natural projection $p_\Upsilon$ from $V_\rho$ to the upper triangle which only keeps the variables $M_{i,j}$ in $\Upsilon$. Consider the similar projection $p'_\Upsilon$ from $\mathcal N_{K,N}\cap A_\rho$ which only keeps variables in $\Upsilon$, i.e., such that their chords in $\rho$ cross. Then we have the isomorphisms
\begin{gather*}
\mathcal N_{K,N} \cap A_\rho \overset{p'_\Upsilon}\cong p'_\Upsilon(\mathcal N_{K,N}\cap A_\rho) = p_\Upsilon(V_\rho) \overset{p_\Upsilon}\cong V_\rho.
\end{gather*}
\end{lem}

The action of $T_N=(\CC^\times)^N\times \CC^\times$ on $\mathcal N_{K,N}$, as defined in Section~\ref{sec:geomsetup}, corresponds to conjugation by diagonal matrices and scaling of $M$, i.e., rewriting the weights~\eqref{eq:weiH} in terms of $M=u^T$,
\begin{gather}\label{eq:boringtorus}
\wt_H(M_{i,j})=\hbar+x_i-x_j,\qquad i,j=1,\ldots,N.
\end{gather}
Only a subgroup of it preserves the slice $A_\rho$; in order for the equation $M_{i,\rho(i)}=1$ to be homogeneous, we must have the relation $x_{\rho(i)}=\hbar+x_i$ if $1\le i<\rho(i)\le N$. This naturally identifies this subgroup with the torus $T_\rho$ defined in Section~\ref{sec:torus}. The isomorphism above then commutes with the action of $T_\rho$, since the weights \eqref{eq:torus} and \eqref{eq:boringtorus} agree for $M_{i,j}\in\Upsilon$ (i.e., $i<j$).

\subsubsection{Poisson structure}\label{sec:poisson}
We view the space of $N\times N$ matrices, $\Mat{N}$, as $\mathfrak{gl}_N^\ast$; the natural Poisson bracket $\{ a,b\}=[a,b]$, where $a,b\in\mathfrak{gl}_N$ are viewed as coordinates on $\Mat{N}$, restricts to the orbit closure $\mathcal{N}_{K,N}$ (i.e., the equations $M^2=0$ form a Poisson ideal), and the orbit itself, namely $\mathcal N_{K,N}\cap \{\rk(M)=K\}$, is a symplectic leaf.

Conjugation by diagonal matrices preserves the symplectic form on $\Mat{N}$ (but not scaling); this forms a subgroup given in terms of the weights~\eqref{eq:torus} by $\hbar=0$.

Once we restrict to the slice $\mathcal N_{K,N} \cap A_\rho$, only the smaller torus $T_\rho$ acts, and setting $\hbar=0$ corresponds to restricting further to a codimension~1 subtorus which preserves the symplectic structure on $\mathcal N_{K,N} \cap A_\rho$.

The corresponding Poisson bracket on $V_\rho$ is
\begin{gather}\label{eq:poi}
\{ M_{i,j}, M_{k,\ell} \} = \delta_{j,k}
 \begin{cases}
1,& i=\rho(\ell),\\
M_{i,\ell},& \cyc{i<\ell<\rho(i)<\rho(\ell)},\\
0,&\text{otherwise},
\end{cases}
\\\notag
\hphantom{\{ M_{i,j}, M_{k,\ell} \}=}{} - \delta_{\ell,i}
 \begin{cases}
1,& k=\rho(j),\\
M_{k,j},& \cyc{k<j<\rho(k)<\rho(j)},\\
0,&\text{otherwise},
\end{cases}
\qquad
M_{i,j},M_{k,\ell}\in \Omega.
\end{gather}

\subsection{The Lagrangian subvarieties}\label{sec:lag}
We shall now define an affine scheme $L_\rho$ inside $V_\rho$.

We first consider the subscheme $\hat L_\rho$ of the space of matrices $\Mat{\ZZ/N\ZZ}$ (square matrices indexed by $\ZZ/N\ZZ$) given by the following equations:
\begin{itemize}\itemsep=0pt
\item $M_{i,\rho(i)}=1$ for all $i\in\ZZ/N\ZZ$.
\item $M_{i,j}=0$ if $\cyc{i<\rho(i)<j}$ or $\cyc{i<\rho(j)<j}$.
\item $M_{i,i}=0$ (these equations are optional, since they are implied set-theoretically by the equations that follow).
\item And quadratic equations:
\begin{gather}\label{eq:quad}
\sum_{\cyc{i\le j\le k}} M_{i,j} M_{j,k} = 0,\qquad i,k\in \ZZ/N\ZZ
\end{gather}
\end{itemize}
(note the similarity with the Brauer loop scheme \cite{artic33}).

As usual, we denote $p_\Omega\colon \Mat{\ZZ/N\ZZ} \to W_\rho$ the natural projection which consists in keeping only the variables in $\Omega$.

We then define $L_\rho:=p_\Omega(\hat L_\rho)$. It is not hard to check the following:
\begin{lem}$p_\Omega$ is an isomorphism from $\hat L_\rho$ to $L_\rho$. Furthermore, $L_\rho \subset V_\rho$.
\end{lem}

The quadratic equations~\eqref{eq:quad} are invariant under the $\hat T_N$-action~\eqref{eq:T0-act}, and so are the equations of the form $M_{i,j}=0$. However the equations $M_{i,\rho(i)}=1$ force the restriction~\eqref{eq:restr} on the weights, which implies that~$L_\rho$, just like $V_\rho$, is only invariant under the subtorus~$T_\rho$.

\subsubsection{Relation to the orbital scheme}
In the same way that $V_\rho$ can be viewed as a slice of the nilpotent orbit $\mathcal N_{K,N}$, $L_\rho$ can be viewed as a slice of the orbital scheme $\mathcal O_{K,N}=\mathcal N_{K,N}\cap \mathfrak n_-$ (being careful that $u\in \mathfrak n_-\ \Leftrightarrow\ M=u^T$ {\em upper} triangular):

\begin{lem}\label{lem:slice2} With the same setup as in Lemma~{\rm \ref{lem:slice}}, we have the isomorphisms
\begin{gather*}
\mathcal O_{K,N} \cap A_\rho \overset{p'_\Upsilon}\cong p'_\Upsilon(\mathcal O_{K,N}\cap A_\rho) = p_\Upsilon(L_\rho) \overset{p_\Upsilon}\cong L_\rho.
\end{gather*}
\end{lem}

\subsubsection{Irreducible components}\label{sec:irrcomp}
Consider for $\pi\in\CLP^{(K)}_{K,N}$, $\pi\ge\rho$, the variety
\begin{gather*}
L_{\rho,\pi}= p_\Upsilon^{-1}p'_\Upsilon(\mathcal O_\pi\cap A_\rho).
\end{gather*}
The decomposition of $\mathcal O_{K,N}$ into irreducible components $\mathcal O_\pi$, $\pi\in \LP^{(K)}_{K,N}$, leads to a similar decomposition
\begin{gather}\label{eq:irrdecomp}
L_\rho = \bigcup_{\pi\in \LP^{(K)}_{K,N}} L_{\rho,\pi}.
\end{gather}
We claim the following fact:
\begin{thm}[\cite{prep1}] Given $\pi\in\CLP_{K,N}$, $\mathcal O_\pi\cap A_\rho\ne \varnothing$ iff $\rho\le\pi$. Assume now that $\rho\le\pi$. Then~$L_{\rho,\pi}$ is irreducible; and the intersection $\mathcal O_\pi\cap A_\rho$ is transverse, so that $\dim L_{\rho,\pi}=\cross(\rho)-\cross(\pi)$ $($and~$L_{\rho,\pi}$, being isotropic, is Lagrangian if~$\pi\in \LP^{(K)}_{K,N})$.
\end{thm}
In particular, the irreducible components of $L_\rho$ are the $L_{\rho,\pi}$ for $\pi\ge\rho$ noncrossing.

The equations of each $L_{\rho,\pi}$ can be written in principle as follows:
\begin{itemize}\itemsep=0pt
\item Start from the equations \eqref{eq:orb} for $\mathcal O_\pi$.
\item Add the equations of $A_\rho$; equivalently, this means that one puts a $1$ at position $(i,j)$ in the matrix $M$ each time $\rho(i)=j$, $i<j$, and then one fills with $0$'s the row above and column to the right of each such~$1$.
\item {\em Eliminate} all the variables $M_{i,j}$ that are not in $\Upsilon$ (that this is possible is guaranteed by Lemma~\ref{lem:slice}).
\end{itemize}

\subsection{Main results}
As before we fix $\rho\in \CLP^{(K)}_{K,N}$, but now we also fix a reduced diagram $\mathcal D$ of $\rho$ (recall that $\rho$ may have several reduced diagrams, which are related by move~\eqref{eq:braid}).

A {\em noncrossing loop configuration} of $\mathcal D$ is the replacement of {\em each} crossing of $\mathcal D$ with
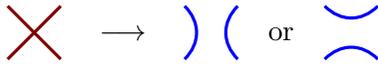
\begin{gather*}
\tikz[rotate=45,baseline=-3pt]{\plaqc}
\quad
\longrightarrow
\quad
\def\linkpatternedgecolor{blue}%
\tikz[rotate=45,baseline=-3pt]{\plaqa}
\quad
\text{or}
\quad
\tikz[rotate=45,baseline=-3pt]{\plaqb}
\end{gather*}
Note that this coincides with the moves defining the order relation in Section~\ref{sec:poset}. Similarly, a~{\em crossing loop configuration} of $\mathcal D$ is the replacement of {\em some} of the crossings of $\mathcal D$ with the same pictures; equivalently, we would rather think of it as the replacement of each crossing of $\mathcal D$ with
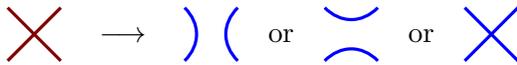
\begin{gather*}
\tikz[rotate=45,baseline=-3pt]{\plaqc}
\quad
\longrightarrow
\quad
\def\linkpatternedgecolor{blue}%
\tikz[rotate=45,baseline=-3pt]{\plaqa}
\quad
\text{or}
\quad
\tikz[rotate=45,baseline=-3pt]{\plaqb}
\quad
\text{or}
\quad
\tikz[rotate=45,baseline=-3pt]{\plaqc}
\end{gather*}

It is sometimes convenient to draw the {\em dual planar map} of $\mathcal D$; since the latter has only regular crossings, this dual planar map is nothing but a quadrangulation of a domain of the plane:
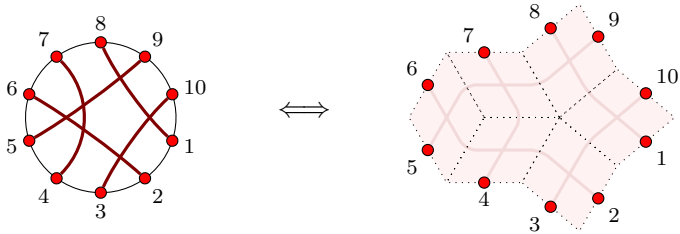
\begin{gather*}
\linkpattern[tikzoptions={rotate=-18},labeloptions={label anchor=#1+180-18},unit=1cm,shape=circle,squareness=0.1]{1/8,2/6,3/10,5/9,4/4/{[bend left=45] to (v7);}}
\qquad
\Longleftrightarrow
\qquad
\begin{tikzpicture}[baseline=0,label distance=-1mm]
\draw[bgplaq] (0,0) coordinate (a) -- (1,0) coordinate (b) -- +(120:1) coordinate (c) -- (120:1) coordinate (d) -- cycle
 (a) -- (b) -- +(-120:1) coordinate (e) -- (-120:1) coordinate (f) -- cycle
 (a) -- (d) -- +(-120:1) coordinate (g) -- (f) -- cycle
 (b) -- (c) -- ++(40:1) coordinate (h) -- ++(300:1) coordinate (i) -- cycle
 (b) -- (e) -- ++(-40:1) coordinate (j) -- ++(-300:1) coordinate (k) -- cycle
 (b) -- (i) -- ++(-40:1) coordinate (l) -- (k) -- cycle;
\draw[/linkpattern/edge,rounded corners,draw opacity=0.1]
($(g)!.5!(d)$) node[/linkpattern/vertex,label={above left:$\ss 6$}] {} -- ($(a)!.5!(f)$) -- ($(b)!.5!(e)$) -- ($(j)!.5!(k)$) node[/linkpattern/vertex,label={below right:$\ss 2$}] {}
($(g)!.5!(f)$) node[/linkpattern/vertex,label={below left:$\ss 5$}] {} -- ($(a)!.5!(d)$) -- ($(c)!.5!(b)$) -- ($(h)!.5!(i)$) node[/linkpattern/vertex,label={above right:$\ss 9$}] {}
($(f)!.5!(e)$) node[/linkpattern/vertex,label={below:$\ss 4$}] {} -- ($(a)!.5!(b)$) -- ($(c)!.5!(d)$) node[/linkpattern/vertex,label={above left:$\ss 7$}] {}
($(e)!.5!(j)$) node[/linkpattern/vertex,label={below left:$\ss 3$}] {} -- ($(b)!.5!(k)$) -- ($(i)!.5!(l)$) node[/linkpattern/vertex,label={above right:$\ss 10$}] {}
($(c)!.5!(h)$) node[/linkpattern/vertex,label={above left:$\ss 8$}] {} -- ($(b)!.5!(i)$) -- ($(k)!.5!(l)$) node[/linkpattern/vertex,label={below right:$\ss 1$}] {};
\end{tikzpicture}
\end{gather*}
\def\linkpatternedgecolor{blue} In this picture, a loop configuration of $\mathcal D$ is an assignment of so-called {\em plaquettes}, i.e., $\tikz[baseline=-3pt]{\plaq{a}}$, $\tikz[baseline=-3pt]{\plaq{b}}$ (or $\tikz[baseline=-3pt]{\plaq{c}}$ in the crossing case) to each face of the dual planar map of $D$, e.g.,
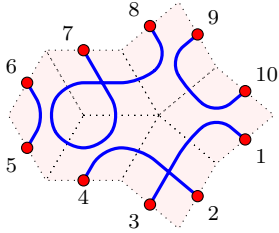
\begin{gather*}
\begin{tikzpicture}[baseline=0,label distance=-1mm]
\draw[bgplaq] (0,0) coordinate (a) -- (1,0) coordinate (b) -- +(120:1) coordinate (c) -- (120:1) coordinate (d) -- cycle
 (a) -- (b) -- +(-120:1) coordinate (e) -- (-120:1) coordinate (f) -- cycle
 (a) -- (d) -- +(-120:1) coordinate (g) -- (f) -- cycle
 (b) -- (c) -- ++(40:1) coordinate (h) -- ++(300:1) coordinate (i) -- cycle
 (b) -- (e) -- ++(-40:1) coordinate (j) -- ++(-300:1) coordinate (k) -- cycle
 (b) -- (i) -- ++(-40:1) coordinate (l) -- (k) -- cycle;
\draw[/linkpattern/edge,rounded corners=2.5mm]
($(g)!.5!(d)$) node[/linkpattern/vertex,label={above left:$\ss 6$}] {} -- ($(a)!.5!(g)$) -- ($(g)!.5!(f)$) node[/linkpattern/vertex,label={below left:$\ss 5$}] {}
($(j)!.5!(k)$) node[/linkpattern/vertex,label={below right:$\ss 2$}] {} -- ($(b)!.5!(e)$) -- ($(a)!.5!(e)$) -- ($(f)!.5!(e)$) node[/linkpattern/vertex,label={below:$\ss 4$}] {}
($(h)!.5!(i)$) node[/linkpattern/vertex,label={above right:$\ss 9$}] {} -- ($(b)!.5!(h)$) -- ($(b)!.5!(i)$) -- ($(b)!.5!(l)$) -- ($(i)!.5!(l)$) node[/linkpattern/vertex,label={above right:$\ss 10$}] {}
($(c)!.5!(d)$) node[/linkpattern/vertex,label={above left:$\ss 7$}] {} -- ($(a)!.5!(b)$) -- ($(a)!.5!(e)$) -- ($(a)!.5!(f)$) -- ($(a)!.5!(g)$) -- ($(a)!.5!(d)$) -- ($(b)!.5!(c)$) -- ($(b)!.5!(h)$) -- ($(c)!.5!(h)$) node[/linkpattern/vertex,label={above left:$\ss 8$}] {}
($(e)!.5!(j)$) node[/linkpattern/vertex,label={below left:$\ss 3$}] {} -- ($(b)!.5!(k)$) -- ($(b)!.5!(l)$) -- ($(k)!.5!(l)$) node[/linkpattern/vertex,label={below right:$\ss 1$}] {};
\end{tikzpicture}
\end{gather*}

To each loop configuration we shall now associate a {\em weight}.\footnote{These (Boltzmann) weights (a terminology borrowed from statistical mechanics) should not be confused with the weights of torus actions.} It is comprised of two parts:
\begin{itemize}\itemsep=0pt
\item A product of {\em local weights}, of the form
\begin{gather*}
\prod_{\text{crossings of $\mathcal D$}} w(\text{plaquette},\text{labels of lines crossing}),
\end{gather*}
where the function $w$ depends on which model we are considering and will be discussed in the next paragraph. The labels of lines follow the diagram of $\rho$, and should not be confused with the lines drawn on the plaquettes of the particular configuration~-- to help with the distinction, lines of $\rho$ are always drawn in red, whereas configuration lines are drawn in blue.
\item A nonlocal part, which is obtained by
applying to the loop configuration the moves of~\eqref{eq:graphrulesnc} (or~\eqref{eq:graphrules} in the crossing case) to produce a reduced diagram. This results in a factor of
\begin{gather*}
(-1)^{|\text{removed crossings}|} \beta^{|\text{removed loops}|},
\end{gather*}
which multiplies the local weight above.

The link pattern of the resulting reduced diagram is called the {\em connectivity} of the loop configuration. It is an element of $\LP^{(K)}_{K,N}$ (resp.\ $\CLP^{(K)}_{K,N}$ in the crossing case).
\end{itemize}

In the example above, the nonlocal weight is $-1$, and the connectivity is \linkpattern[shape=circle]{5/6,7/8,9/10,1/3,2/4}.

The {\em partition function} of the loop model on $\mathcal D$ is by definition the sum of weights of all loop configurations of $\mathcal D$ with some prescribed connectivity.

\subsubsection{Cohomology}
First, in order to distinguish the two types of noncrossing plaquettes at each vertex, we shall orient every line corresponding to a pairing from $i$ to $\rho(i)$ where $i<\rho(i)$, and associate to it a~formal parameter~$x_i$. The local weights at each vertex are then given by
\begin{gather*}
w\left(
\begin{tikzpicture}[x=\loopcellsize,y=\loopcellsize,baseline=-3pt]
\draw[bgplaq] (-0.5,-0.5) rectangle ++(1,1);
\draw[/linkpattern/edge,colb,->] (-0.7,0) node[left,black] {$x_j$} -- (0.9,0);
\draw[/linkpattern/edge,colb,->] (0,-0.7) node[below,black] {$x_i$} -- (0,0.9);
\end{tikzpicture}
\right)
=
\begin{cases}
\hbar-x_j+x_i
&\tikz[baseline=-3pt]{\plaq{a}}\\[2mm]
x_j-x_i
&\tikz[baseline=-3pt]{\plaq{b}}
\end{cases}
\end{gather*}

We can now state our first main result:
\begin{thm}[\cite{artictt}]\label{thm:fdcoh}
The class of $L_{\rho,\pi}$ in $H_{T_\rho}(V_\rho)\cong H_{T_\rho}(\cdot)$ is given by the partition function
of the noncrossing loop model on any diagram $\mathcal D$ of $\rho$, with loop weight $\beta=2$ and connectivity $\pi$, that is
\begin{gather*}
[L_{\rho,\pi}]=\sum_{\substack{\text{\rm loop configurations on $\mathcal D$}\\\text{\rm with connectivity $\pi$}}}
2^{|\text{\rm loops}|}
\prod_{\substack{i<j<\rho(i)<\rho(j)\\\text{\rm crossing of $\rho$}}}
\begin{cases}
\hbar-x_j+x_i
&\tikz[baseline=-3pt]{\plaq{a}}\\[2mm]
x_j-x_i
&\tikz[baseline=-3pt]{\plaq{b}}
\end{cases}
\end{gather*}
\end{thm}
An idea of the proof will be given in Section~\ref{sec:proof}; for now, we simply point out that the weights of the claim are proportional to the coefficients of the $R$-matrix of Claim~\ref{thm:coh}, where the first plaquette plays the role of identity operator and the second one, the role of Temperley--Lieb generator. The claim then says that the cohomology classes of the $L_{\rho,\pi}$ for given $\rho$ and varying~$\pi$, are proportional to a product of such $R$-matrices.

\subsubsection[$K$-theory 1: noncrossing loops]{$\boldsymbol{K}$-theory 1: noncrossing loops}\label{sec:findomnc}
A similar statement can be made in $K$-theory. For $\pi\in \LP^{(K)}_{K,N}$, define the coherent sheaves $\sigma_\pi$ as in Section~\ref{sec:K1}
with $n=N$, $k=K$. Consider then $\sigma_{\rho,\pi}$, which is by definition the restriction of~$\mu_\ast \sigma_\pi$ to~$A_\rho$. Via the isomorphisms of Lemma~\ref{lem:slice2}, we view $\sigma_{\rho,\pi}$ as sheaves on~$V_\rho$, with support~$L_{\rho,\pi}$. $\sigma_{\rho,\pi}$ is naturally $T'_\rho$-equivariant, where $T'_\rho$ is the preimage of $T_\rho$ under the double cover $T'_N\to T_N$.

Compared to the case of cohomology, we replace variables $\hbar,x_1,\ldots,x_N$ with variables $t,z_1$, $\ldots,z_N$ in the obvious way. In particular, we attach variables $z_i$ to each pairing $i<\rho(i)$, having eliminated other $z$'s thanks to $z_{\rho(i)}=t\,z_i$. We then have the generalization of the previous result:
\begin{thm}[\cite{artictt}]\label{thm:fdnc}
The class of $\sigma_{\rho,\pi}$ in $K_{T'_\rho}(V_\rho)\cong K_{T'_\rho}(\cdot)$ is given $($up to a monomial$)$ by the partition function of the noncrossing loop model on any diagram $\mathcal D$ of~$\rho$, with loop weight $\beta=t^{1/2}+t^{-1/2}$ and connectivity~$\pi$, that is
\begin{gather*}
[\sigma_{\rho,\pi}]=m_\rho\sum_{\substack{\text{\rm loop configurations on $\mathcal D$}\\\text{\rm with connectivity $\pi$}}}
\big(t^{1/2}+t^{-1/2}\big)^{|\text{\rm loops}|}
\prod_{\substack{i<j<\rho(i)<\rho(j)\\\text{\rm crossing of $\rho$}}}
\begin{cases}
1-t\,z_i/z_j
&\tikz[baseline=-3pt]{\plaq{a}}\\[2mm]
t^{1/2}(z_i/z_j-1)
&\tikz[baseline=-3pt]{\plaq{b}}
\end{cases}
\end{gather*}
where $m_\rho$ is the monomial defined in \eqref{eq:weinorm}.
\end{thm}
Again, note that the local weights are nothing but the $R$-matrix of Claim~\ref{thm:nc}, up to norma\-li\-zation.

We shall not describe more explicitly the sheaves $\sigma_{\rho,\pi}$ in full generality, though two examples will be treated in Sections~\ref{sec:hexa} and~\ref{sec:back}.

\subsubsection[$K$-theory 2: crossing loops]{$\boldsymbol{K}$-theory 2: crossing loops}
Finally, we connect $K$-classes structure sheaves of these varieties $L_{\rho,\pi}$ to the crossing loop model:
\begin{thm}[\cite{artictt}]\label{thm:fdcr}\allowdisplaybreaks
The class of $($the structure sheaf of$)$ $L_{\rho,\pi}$ in $K _{T_\rho}(V_\rho)\cong K_{T_\rho}(\cdot)$ is given by the partition function of the crossing loop model on any diagram $\mathcal D$ of $\rho$, with loop weight $\beta=1+t$ and connectivity~$\pi$, that is
\begin{gather*}
[L_{\rho,\pi}]=
\sum_{\substack{\text{\rm loop configurations on $\mathcal D$}\\\text{\rm with connectivity $\pi$}}}
(-1)^{|\text{\rm crossings}|}(1+t)^{|\text{\rm loops}|}\\
\hphantom{[L_{\rho,\pi}]=}{}\times
\prod_{\substack{i<j<\rho(i)<\rho(j)\\\text{\rm crossing of $\rho$}}}
\begin{cases}
1-t\,z_i/z_j
&\tikz[baseline=-3pt]{\plaq{a}}
\\[2mm]
1-z_j/z_i
&\tikz[baseline=-3pt]{\plaq{b}}
\\[2mm]
(1-t\,z_i/z_j)(1-z_j/z_i)
&\tikz[baseline=-3pt]{\plaq{c}}
\end{cases}
\end{gather*}
$($where $|\text{crossings}|$ means as before the number of crossings removed via~\eqref{eq:graphrules} to produce a reduced diagram$)$.
\end{thm}
Once again, the local weights reproduce the $R$-matrix of Claim~\ref{thm:cr} up to normalization.

\subsubsection{The hexagon}\label{sec:hexa}
Since the whole construction of this section may seem somewhat abstract, we provide a full description of one important special case, namely the maximally crossing link pattern for $N=2K=6$.

\begin{ex}\label{ex:hexa}
We choose $\rho$ to be the $\ZZ/6\ZZ$ involution $\rho(i)=i+3$; it has two possible diagrams among which we choose one:
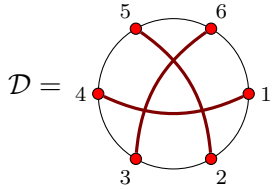
\begin{gather*}
\mathcal D=
\linkpattern[edgecolor=colb,unit=1cm,shape=circle,squareness=0.1,size=6]{1/1/{[bend right=25] to (v4);},2/2/{[bend left=25] to (v5);},3/3/{[bend right=25] to (v6);}}
\end{gather*}
We select as one possible set of coordinates on $V_\rho$ $\{M_{1,2},M_{2,4},M_{3,4},M_{4,6},M_{5,6},M_{6,2}\}$; note that
they satisfy canonical Poisson brackets (all other brackets are zero):
\begin{gather*}
\{M_{1,2},M_{2,4}\} = \{M_{3,4},M_{4,6}\} = \{M_{5,6},M_{6,2}\}=1.
\end{gather*}
The equations of $L_\rho$, once all other variables are eliminated, are
\begin{gather*}
L_\rho=\big\{ M_{1,2}M_{2,4}=M_{3,4}M_{4,6}=M_{5,6}M_{6,2},\\
\hphantom{L_\rho=\big\{}{} (M_{1,2}+M_{4,6}M_{6,2})M_{3,4}
=(M_{3,4}+M_{6,2}M_{2,4})M_{5,6}=(M_{5,6}+M_{2,4}M_{4,6})M_{1,2}=0\big\}.
\end{gather*}
Using for example Macaulay~2 \cite{M2}, one can decompose this affine scheme into its 5 irreducible components, namely
\begin{gather*}
L_{\rho,\linkpattern[shape=circle,small]{1/2,3/4,5/6}}=\{M_{1,2}+M_{4,6}M_{6,2}=M_{3,4}+M_{6,2}M_{2,4}=M_{5,6}+M_{2,4}M_{4,6}=0\},\\
L_{\rho,\linkpattern[shape=circle,small]{1/6,2/3,4/5}}=\{M_{1,2}=M_{3,4}=M_{5,6}=0\},\\
L_{\rho,\linkpattern[shape=circle,small]{1/2,3/6,4/5}}=\{M_{2,4}=M_{3,4}=M_{5,6}=0\},\\
L_{\rho,\linkpattern[shape=circle,small]{1/4,2/3,5/6}}=\{M_{1,2}=M_{3,4}=M_{6,2}=0\},\\
L_{\rho,\linkpattern[shape=circle,small]{1/6,2/5,3/4}}=\{M_{1,2}=M_{4,6}=M_{5,6}=0\},
\end{gather*}
and compute their $K$-classes, comparing them with crossing partition functions:
\tikzset{hexa/.pic={
\draw[bgplaq] (90:1) -- node[midway,/linkpattern/vertex] (6) {} (150:1) -- node[midway,/linkpattern/vertex] (1) {} (210:1) -- node[midway,/linkpattern/vertex] (2) {} (270:1) -- node[midway,/linkpattern/vertex] (3) {} (330:1) -- node[midway,/linkpattern/vertex] (4) {} (30:1) -- node[midway,/linkpattern/vertex] (5) {} (90:1) (0,0) -- coordinate[midway] (a) (-90:1) (0,0) -- coordinate[midway] (b) (150:1) (0,0) -- coordinate[midway] (c) (30:1);
}}
\begin{gather}
[L_{\rho,\linkpattern[shape=circle,small]{1/2,3/4,5/6}}] =(1-t z_1/z_2)(1-z_3/z_1)(1-t z_2/z_3) =
\tikz[baseline=-3pt,scale=0.7]{\pic[transform shape]{hexa}; \draw[/linkpattern/edge,bend left] (1) to (2) (3) to (4) [bend right] (5) to (c) [bend left] to (a) to (b) [bend right] to (6);}
+
\tikz[baseline=-3pt,scale=0.7]{\pic[transform shape]{hexa}; \draw[/linkpattern/edge,bend right] (1) to (b) [bend left] to (c) to (a) [bend right] to (2) [bend left] (3) to (4) (5) to (6);}
+
\tikz[baseline=-3pt,scale=0.7]{\pic[transform shape]{hexa}; \draw[/linkpattern/edge,bend left] (1) to (2) [bend right] (3) to (a) [bend left] to (b) to (c) [bend right] to (4) [bend left] (5) to (6);} \nonumber\\
\label{eq:hexa1} \hphantom{[L_{\rho,\linkpattern[shape=circle,small]{1/2,3/4,5/6}}]=}{}
 +(1+t)\,\tikz[baseline=-3pt,scale=0.7]{\pic[transform shape]{hexa}; \draw[/linkpattern/edge,bend left] (1) to (2) (3) to (4) (5) to (6) (a) to (b) to (c) to (a);}
 - \tikz[baseline=-3pt,scale=0.7]{\pic[transform shape]{hexa}; \draw[/linkpattern/edge,bend left] (1) to (2) (3) to (4) [bend left=0] (5) to (b) [bend right=30] to (a) to (c) [bend left=0] to (6);}
-
\tikz[baseline=-3pt,scale=0.7]{\pic[transform shape]{hexa}; \draw[/linkpattern/edge] (1) to (a) [bend right] to (c) to (b) [bend left=0] to (2) [bend left=30] (3) to (4) (5) to (6);}
-
\tikz[baseline=-3pt,scale=0.7]{\pic[transform shape]{hexa}; \draw[/linkpattern/edge,bend left] (1) to (2) [bend left=0] (3) to (c) [bend right=30] to (b) to (a) [bend right=0] to (4) [bend left=30] (5) to (6);}\,,\\
\label{eq:hexa2}
[L_{\rho,\linkpattern[shape=circle,small]{1/6,2/3,4/5}}] =(1-t z_1/z_2)(1-z_3/z_1)(1-t z_2/z_3)
 = \tikz[baseline=-3pt,scale=0.7]{\pic[transform shape]{hexa}; \draw[/linkpattern/edge,bend right] (1) to (b) to (6) (5) to (c) to (4) (3) to (a) to (2);}\,,\\
\label{eq:hexa3}
[L_{\rho,\linkpattern[shape=circle,small]{1/2,3/6,4/5}}] =(1-z_2/z_1)(1-z_3/z_1)(1-t z_2/z_3)
 = \tikz[baseline=-3pt,scale=0.7]{\pic[transform shape]{hexa}; \draw[/linkpattern/edge,bend right] (1) to (b) to (6) (5) to (c) [bend left] to (a) [bend right] to (2) (4) to (3);}\,,\\
\label{eq:hexa4}
[L_{\rho,\linkpattern[shape=circle,small]{1/4,2/3,5/6}}]=(1-tz_1/z_2)(1-z_3/z_1)(1-z_3/z_2)
=\tikz[baseline=-3pt,scale=0.7]{\pic[transform shape]{hexa}; \draw[/linkpattern/edge,bend right] (1) to (b) [bend left] to (c) [bend right] to (4) [bend left] (5) to (6) [bend right] (3) to (a) to (2);}\,, \\
\label{eq:hexa5}
[L_{\rho,\linkpattern[shape=circle,small]{1/6,2/5,3/4}}] =(1-t z_1/z_2)(1-t z_1/z_3)(1-t z_2/z_3)
= \tikz[baseline=-3pt,scale=0.7]{\pic[transform shape]{hexa}; \draw[/linkpattern/edge,bend left] (1) to (2) [bend right] (3) to (a) [bend left] to (b) [bend right] to (6) (5) to (c) to (4);}\,,
\end{gather}
where the thumb rule to compute $K$-classes is that each $M_{i,j}=0$ equation contributes a factor $1-t z_i/z_j$, with the additional replacement $z_{\rho(i)}=t z_i$ with $i<\rho(i)$, and with the convention that each loop model configuration represents its product of local weights.

The isotropic varieties corresponding to crossing link patterns can be treated similarly.

Of course one could have used the other diagram of $\rho$ to perform the computation of the $K$-classes; the fact that the result is identical is nothing but the Yang--Baxter equation. We shall reinvestigate this point in Section~\ref{sec:degen}.

One can equally easily compute the $K$-classes of the sheaves $\sigma_{\rho,\pi}$. For all but $\pi=\linkpattern[shape=circle,small]{1/6,2/3,4/5}$, $\sigma_\pi$~is nonequivariantly trivial (structure sheaf), and therefore so is~$\sigma_{\rho,\pi}$. $\sigma_{\linkpattern[shape=circle,small]{1/6,2/3,4/5}}$ is the $O(1)$ sheaf; once sliced to the linear subspace
$L_{\rho,\linkpattern[shape=circle,small]{1/6,2/3,4/5}}$, it becomes a trivial sheaf but with an extra weight of $z_3 z_4^{-1}=t^{-1}z_1^{-1}z_3$. This way, we get the following weights for $\sigma_{\rho,\pi}$ (from which we have factored out $m_\rho$ in view of Claim~\ref{thm:fdnc}):
\begin{gather*}
m_{\linkpattern[shape=circle,small]{1/2,3/4,5/6}} =t z_1z_3^{-1}m_\rho,\\
t^{-1}z_1^{-1}z_3 m_{\linkpattern[shape=circle,small]{1/6,2/3,4/5}} =t^{1/2}z_1 z_3^{-1}m_\rho,\\
m_{\linkpattern[shape=circle,small]{1/2,3/6,4/5}} =t z_1^2z_2^{-1}z_3^{-1}m_\rho,\\
m_{\linkpattern[shape=circle,small]{1/4,2/3,5/6}} =t z_1z_2z_3^{-2}m_\rho,\\
m_{\linkpattern[shape=circle,small]{1/6,2/5,3/4}} =m_\rho.
\end{gather*}
We get the classes $[\sigma_{\rho,\pi}]$ by multiplying $[L_{\rho,\pi}]$ with these weights. We can then compare them with partition functions of the noncrossing loop model:{\allowdisplaybreaks
\begin{gather*}
[\sigma_{\rho,\linkpattern[shape=circle,small]{1/2,3/4,5/6}}] =m_\rho t(1-t z_1/z_2)(z_1/z_3-1)(1-t z_2/z_3)\\
\hphantom{[\sigma_{\rho,\linkpattern[shape=circle,small]{1/2,3/4,5/6}}]}{}
= m_\rho\Big(\tikz[baseline=-3pt,scale=0.7]{\pic[transform shape]{hexa}; \draw[/linkpattern/edge,bend left] (1) to (2) (3) to (4) [bend right] (5) to (c) [bend left] to (a) to (b) [bend right] to (6);}
+
\tikz[baseline=-3pt,scale=0.7]{\pic[transform shape]{hexa}; \draw[/linkpattern/edge,bend right] (1) to (b) [bend left] to (c) to (a) [bend right] to (2) [bend left] (3) to (4) (5) to (6);}
+
\tikz[baseline=-3pt,scale=0.7]{\pic[transform shape]{hexa}; \draw[/linkpattern/edge,bend left] (1) to (2) [bend right] (3) to (a) [bend left] to (b) to (c) [bend right] to (4) [bend left] (5) to (6);}
+(t^{1/2}+t^{-1/2})\,\tikz[baseline=-3pt,scale=0.7]{\pic[transform shape]{hexa}; \draw[/linkpattern/edge,bend left] (1) to (2) (3) to (4) (5) to (6) (a) to (b) to (c) to (a);}\Big),
\\
[\sigma_{\rho,\linkpattern[shape=circle,small]{1/6,2/3,4/5}}] =m_\rho t^{1/2}(1-t z_1/z_2)(z_1/z_3-1)(1-t z_2/z_3)
 = m_\rho \tikz[baseline=-3pt,scale=0.7]{\pic[transform shape]{hexa}; \draw[/linkpattern/edge,bend right] (1) to (b) to (6) (5) to (c) to (4) (3) to (a) to (2);}\,,
\\
[\sigma_{\rho,\linkpattern[shape=circle,small]{1/2,3/6,4/5}}] =m_\rho t(z_1/z_2-1)(z_1/z_3-1)(1-t z_2/z_3)
 = m_\rho \tikz[baseline=-3pt,scale=0.7]{\pic[transform shape]{hexa}; \draw[/linkpattern/edge,bend right] (1) to (b) to (6) (5) to (c) [bend left] to (a) [bend right] to (2) (4) to (3);}\,,
\\
[\sigma_{\rho,\linkpattern[shape=circle,small]{1/4,2/3,5/6}}] =m_\rho t(1-t\,z_1/z_2)(z_1/z_3-1)(z_2/z_3-1)
 = m_\rho\,\tikz[baseline=-3pt,scale=0.7]{\pic[transform shape]{hexa}; \draw[/linkpattern/edge,bend right] (1) to (b) [bend left] to (c) [bend right] to (4) [bend left] (5) to (6) [bend right] (3) to (a) to (2);}\,,
\\
[\sigma_{\rho,\linkpattern[shape=circle,small]{1/6,2/5,3/4}}] =m_\rho(1-t\,z_1/z_2)(1-t\,z_1/z_3)(1-t z_2/z_3)
 = m_\rho\,\tikz[baseline=-3pt,scale=0.7]{\pic[transform shape]{hexa}; \draw[/linkpattern/edge,bend left] (1) to (2) [bend right] (3) to (a) [bend left] to (b) [bend right] to (6) (5) to (c) to (4);}\,.
\end{gather*}
The first equality is particularly remarkable, because it shows that up to an overall monomial, we could have calculated $[L_{\rho,\linkpattern[shape=circle,small]{1/2,3/4,5/6}}]$ by summing only {\em noncrossing} loop configurations only! A~similar phenomenon will be discussed in Section~\ref{sec:goren}.}
\end{ex}

\subsubsection{Idea of proof}\label{sec:proof}
Even though proofs are not provided in this paper, the results above are sufficiently important, and the essence of their proof sufficiently simple, that it is worth briefly mentioning here. We derive Claim~\ref{thm:fdcr}; the other results can be obtained similarly.

Fix $\pi\in\CLP^{(K)}_{K,N}$. The first statement is that the intersection $\mathcal O_\pi\cap A_\rho$ is {\em transverse} in $\mathfrak n_-$,
so that the $T_\rho$-equivariant $K$-class of $L_{\rho,\pi}$ is equal to that of $\mathcal O_\pi$ up to some prefactors:
\begin{gather*}
\frac{[L_{\rho,\pi}]_{T_\rho,V_\rho}}{[0]_{T_\rho,V_\rho}} =
\frac{[\mathcal O_\pi\cap A_\rho]_{T_\rho,\mathfrak n_-\cap A_\rho}}{[0]_{T_\rho,\mathfrak n_-\cap A_\rho}}
=\frac{[\mathcal O_\pi]_{T_\rho,\mathfrak n_-}}{[0]_{T_\rho,\mathfrak n_-\cap A_\rho}},
\end{gather*}
where the subscripts specify the choice of torus, as well as the different embedding spaces. Writing $\Psi_\pi=[\mathcal O_\pi]_{T_N,\mathfrak n_-}$ (as in Section~\ref{sec:K2}, except we drop the superscript $c$ for convenience), we have the simple identity resulting from the embedding $T_\rho\subset T_N$:
\begin{gather*}
[\mathcal O_\pi]_{T_\rho,\mathfrak n_-}= \Psi_\pi|_{z_{\rho(i)}=t z_i,\, i<\rho(i)}.
\end{gather*}

Combining these two equalities and computing $[0]_{T_\rho,\mathfrak n_-\cap A_\rho}$ from the definition of $A_\rho$ (taking the product of weights of ``red stars'' in \eqref{eq:stars}), we conclude that
\begin{gather}
[L_{\rho,\pi}]_{T_\rho,V_\rho} = \prod_{i<j<\rho(i)<\rho(j)} (1-t\,z_i/z_j)^{-1} \prod_{i<\rho(i)<j<\rho(j)} (1-t z_i/z_j)^{-1}\big(1-t^2 z_i/z_j\big)^{-1}\nonumber\\
\hphantom{[L_{\rho,\pi}]_{T_\rho,V_\rho} =}{} \times \prod_{i<j<\rho(j)<\rho(i)} (1-t z_i/z_j)^{-1}(1-t z_j/z_i)^{-1}
 \Psi_\pi|_{z_{\rho(i)}=tz_i,\, i<\rho(i)}.\label{eq:spec}
\end{gather}

So we are led to the calculation of an appropriate specialization of $\Psi_\pi$. The latter is based on two fundamental facts:
\begin{itemize}\itemsep=0pt
\item The {\em exchange relation}~\eqref{eq:exchcr}, which we rewrite here in vector notation:
\begin{gather}\label{eq:exchcr2}
\ket{\Psi} \tau_i = \check R'_i \ket{\Psi},
\end{gather}
where $\ket{\Psi}=\sum_{\pi\in\CLP_{K,N}^{(K)}} \Psi_\pi \ket{\pi}$,
\begin{gather*}
\check R'_i = \frac{a(z_i/z_{i+1})+b(z_i/z_{i+1})e_i+a(z_i/z_{i+1})b(z_i/z_{i+1})f_i}{a(z_{i+1}/z_{i})}
\end{gather*}
and $a(z)=1-t/z$, $b(z)=1-z$.

\item Let $\mathcal D$ be a diagram of $\rho$. If $\mathcal D$ has connected components $\mathcal D_1,\ldots,\mathcal D_\ell$, e.g.,
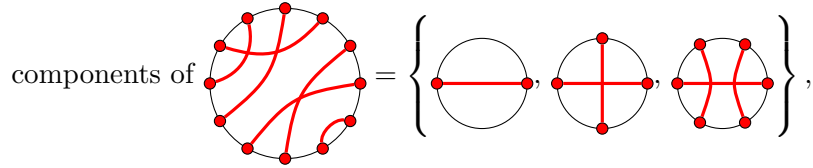
\begin{gather*}
\text{components of }\linkpattern[shape=circle,numbered=false,unit=1cm,edgecolor=red,looseness=0.15]{1/5,2/3,4/12,6/10,7/9,8/11}=
\left\{
\linkpattern[shape=circle,numbered=false,edgecolor=red]{1/2},\,
\linkpattern[shape=circle,numbered=false,edgecolor=red]{1/3,2/4},\,
\linkpattern[shape=circle,numbered=false,edgecolor=red]{1/4,2/6,3/5}
\right\},
\end{gather*}
then $V_\rho\cong V_{\rho_1}\times\cdots\times V_{\rho_\ell}$ and $L_\rho\cong L_{\rho_1}\times\cdots\times L_{\rho_\ell}$, where~$\rho_i$ is the crossing link pattern with diagram $\mathcal D_i$ (up to relabelling of vertices). This is essentially obvious from the definitions of $V_\rho$ and $L_\rho$, cf. Sections~\ref{sec:amb} and~\ref{sec:lag}.

Furthermore, any $\pi\ge\rho$ has a decomposition into connected components that refines that of~$\rho$, which means it is of the form $\pi=\pi_1\sqcup\cdots\sqcup\pi_\ell$ where~$\pi_i$ connects the same vertices as~$\rho_i$, and the diagrams of the $\pi_i$ are disjoint.
This means that irreducible components of~$L_\rho$ factor analogously:
\begin{gather}\label{eq:decomp}
L_{\rho,\pi} = L_{\rho_1,\pi_1}\times\cdots\times L_{\rho_\ell,\pi_\ell}.
\end{gather}

The same equation then holds for $K$-classes: $[L_{\rho,\pi}] = [L_{\rho_1,\pi_1}]\cdots [L_{\rho_\ell,\pi_\ell}]$. Note that this result is compatible with the formula of Claim~\ref{thm:fdcr}, since partition functions trivially factorize in the same manner, thereby allowing us to restrict ourselves to the case that~$\mathcal D$ is connected.
\end{itemize}

The strategy is now clear: since the various specializations of $\ket{\Psi}$ corresponding to different~$\rho$ are related by swapping variables $z_i$, we shall apply repeatedly the exchange relation~\eqref{eq:exchcr2} to remove a crossing from~$\rho$, each time adding an extra plaquette to the partition function; at the end of the day, the domain of the partition function will reproduce a diagram of~$\rho$.

More explicitly, we proceed inductively on the number of crossings $\cross(\rho)$ of $\rho$.
\begin{itemize}\itemsep=0pt
\item If $\cross(\rho)=0$, then $L_\rho\cong V_\rho\cong \{\cdot\}$, and
\begin{gather}\label{eq:trivfd}
[L_{\rho,\pi}]=\delta_{\rho,\pi},\qquad \rho\in \LP_{K,N}^{(K)},
\end{gather}
which matches trivially the formula of Claim~\ref{thm:fdcr}.
\item Now assume $\cross(\rho)>0$. As was explained above using the decomposition of $\rho$ into connected components, one may assume that the reduced diagram $\mathcal D$ of $\rho$ under consideration is connected; since $\mathcal D$ has at least one crossing, it is not hard to conclude that there must be an $i$ such that $i<i+1<\rho(i)<\rho(i+1)$, i.e., we have the following decomposition
\begin{gather}\label{eq:bdrycross}
\begin{tikzpicture}[baseline=-3pt]
\fill[bgplaq] (0,0) circle [radius=1cm];
\linkpattern[unit=1cm,shape=circle,size=10,numbering={1/,2/,3/,4/,5/,6/,7/,8/i,9/i+1,10/}]{};
\node {$\mathcal D$};
\end{tikzpicture}
=
\begin{tikzpicture}[baseline=-3pt]
\fill[bgplaq] (0,0) circle [radius=0.6cm];
\linkpattern[unit=1cm,shape=circle,size=10,numbering={1/,2/,3/,4/,5/,6/,7/,8/i,9/i+1,10/}]{};
\foreach\x in {1,...,10} \node also [alias=w\x] (v\x);
\linkpattern[unit=0.6cm,shape=circle,size=10,numbered=false]{};
\draw[/linkpattern/edge,colb,bend left=15,arrow=0.4] (w8) to (v9);
\draw[/linkpattern/edge,colb,bend right=15,arrow=0.4] (w9) to (v8);
\foreach\x in {1,...,7,10} \draw[/linkpattern/edge,colb] (w\x) -- (v\x);
\node {$\mathcal D'$};
\end{tikzpicture}
,\qquad \rho=f_i \rho',
\end{gather}
where $\rho'$ is the link pattern with diagram $\mathcal D'$. We then apply exchange relation \eqref{eq:exchcr2} to relate the two specializations corresponding to $\rho$ and $\rho'$:
\begin{gather}\label{eq:exchcr3}
\ket{\Psi}|_{z_{\rho(j)}=t z_j,\, j<\rho(j)} = ((\ket{\Psi}\tau_i)|_{z_{\rho'(j)}=t z_j,\, j<\rho'(j)})\tau_i
=(\check R'_i \ket{\Psi}|_{z_{\rho(j)}=t z_j,\, j<\rho'(j)})\tau_i.
\end{gather}
Carefully taking care of the prefactors in \eqref{eq:spec}, and introducing the vector notation $\ket{L_\rho}=\sum_{\pi\in\CLP_{K,N}^{(K)}}
[L_{\rho,\pi}]_{T_\rho,V_\rho} \ket{\pi}$, we conclude that
\begin{gather}\label{eq:exchcr4}
\ket{L_\rho} = (\check r_i \ket{L_{\rho'}})\tau_i,
\end{gather}
where $\check r_i$ is yet another normalization of the $R$-matrix, namely
\begin{gather*}
\check r_i = a(z_i/z_{i+1})+b(z_i/z_{i+1})e_i+a(z_i/z_{i+1})b(z_i/z_{i+1})f_i.
\end{gather*}
(Note that $\check r_i$ does not satisfy the unitarity equation.)

Now we apply the induction hypothesis to $\rho'$, which has one less crossing than $\rho$: the components of~$\ket{L_{\rho'}}$ form the partition function at fixed connectivity on the domain~$\mathcal D'$. The relation \eqref{eq:exchcr4} is nothing but the addition of the linear combination of plaquettes $ a(z_i/z_{i+1})\tikz[baseline=-3pt]{\plaq{a}}+b(z_i/z_{i+1})\tikz[baseline=-3pt]{\plaq{b}}+ a(z_i/z_{i+1})b(z_i/z_{i+1})\tikz[baseline=-3pt]{\plaq{c}}$ to the boundary of~$\mathcal D'$, followed by the appropriate relabelling of the boundary vertices, producing exactly $\mathcal D$ (cf.~\eqref{eq:bdrycross}). This shows the induction hypothesis for $\mathcal D$ and $\rho$.
\end{itemize}

\subsection{Degeneration}\label{sec:degen}
The equalities of Claims~\ref{thm:fdcoh}, \ref{thm:fdcr} are fairly suggestive, in the sense that {\em each term} in the summation of the r.h.s.\ can be interpreted as the class of a certain subvariety of $V_\rho$. In fact, ignoring the global weight of $2^{|\text{loops}|}$, we recognize the class of a {\em linear} subvariety given by the equations $M_{i,j}=0$ or $M_{j,\rho(i)}=0$, since the local weights of plaquettes in the r.h.s.\ match the factors associated to these equations (related to their additive or multiplicative weights~-- in the sense of torus action~-- for cohomology and $K$-theory respectively). In particular, the crossing plaquette corresponds to the intersection $M_{i,j}=M_{j,\rho(i)}=0$ of the two noncrossing plaquettes.

Such a situation is fairly typical, see, e.g.,~\cite{KM-Schubert} in a closely related context. It is tempting to speculate that there exists a~(Gr\"obner degeneration) of our scheme $L_\rho$ into a so-called Stanley--Reisner scheme, that is a reduced union of coordinate subspaces, such that each term in the sum of Claim~\ref{thm:fdcoh} (i.e., each noncrossing loop configuration) corresponds to one such coordinate subspace (and a similar statement for Claim~\ref{thm:fdcr} in terms of the simplicial complex associated to them, see, e.g.,~\cite{MS-book}).

Interestingly enough, the situation is more subtle: one finds that for each diagram $\mathcal D$ of $\rho$, there exists a~(possibly nonunique) degeneration of $L_\rho$ into a~{\em nonreduced} union of coordinate subspaces, the nonreducedness being responsible for the loop weight. This is the content of the following
\begin{conj}\label{conj:degen} There exists a $($partial$)$ Gr\"obner $T_\rho$-equivariant degeneration of~$L_\rho$ into an $($in general unreduced$)$ scheme whose components are indexed by noncrossing loop configurations of~$\mathcal D$, such that the components coming from $L_{\rho,\pi}$ are indexed by loop configurations with connectivity~$\pi$; each geometric component is given by the equations
\begin{gather*}
\begin{cases}
M_{i,j}=0
&\tikz[baseline=-3pt]{\plaq{a}}\\[2mm]
M_{j,\rho(i)}=0
&\tikz[baseline=-3pt]{\plaq{b}}
\end{cases},
\qquad
i<j<\rho(i)<\rho(j)\ \text{crossing of $\rho$},
\end{gather*}
and has {\em multiplicity} $2^{|\text{\rm loops}|}$.
\end{conj}
(``Partial'' degeneration means here that it does not always produce a monomial ideal.) In particular, note that the {\em radical} of the ideal of the degeneration of the whole of $L_\rho$ is simply given by the equations $M_{i,j}M_{j,\rho(i)}=0$ at each crossing, corresponding to the two choices of noncrossing plaquettes at that crossing.

Conjecture~\ref{conj:degen} directly implies Claim~\ref{thm:fdcoh} (though, typically, a proof of Conjecture~\ref{conj:degen} would rely on Claim~\ref{thm:fdcoh}); it is a bit more subtle to justify the loop weight $1+t$ in Claim~\ref{thm:fdcr} (this requires some further trickery to reduce to the case of a monomial ideal, whose simplicial complex can then be studied using standard methods, leading to crossing loop configurations), and we shall only do so in examples. Of course, Claim~\ref{thm:fdnc} should also follow by degenerating the sheaves $\sigma_{\rho,\pi}$ together with their support $L_{\rho,\pi}$.

There should also be a way to interpolate from the symplectic structure on $V_\rho$ to the one on~$T^\ast \CC^{\cross(\rho)}$ (where $M_{i,j}$ and $M_{j,\rho(i)}$ are canonically conjugate), corresponding to the special fiber, such that each fiber of the degeneration is Lagrangian; in particular, it is clear from their form (choice of either $M_{i,j}=0$ or $M_{j,\rho(i)}=0$ at each crossing) that each geometric component of the special fiber is a Lagrangian coordinate subspace.

The degeneration can be defined inductively, in a similar fashion as the partition function itself was built in Section~\ref{sec:proof}. Pick an outer crossing $(i,i+1)$ of $\rho$ as in \eqref{eq:bdrycross}; then choose as variables~$M_{i,i+1}$ (the variable ``facing outwards'', as opposed to $M_{\rho(i),\rho(i+1)}$ which is ``facing inwards'') and~$M_{i+1,\rho(i)}$ (or~$M_{\rho(i+1),i}$, it does not matter since they are opposite of each other), and an arbitrary set of variables at other crossings. The rule is then ``revlex $M_{i,i+1}$, lex $M_{i+1,\rho(i)}$, keeping their product fixed''; in other words, substitute $M_{i,i+1}\to \epsilon\, M_{i,i+1}$, $M_{i+1,\rho(i)}\to \epsilon^{-1} M_{i+1,\rho(i)}$, and take the leading behavior of the ideal of equations of $V_\rho$ as $\epsilon\to0$.

At present, no non-inductive definition of the degeneration is known, which makes computations difficult, except in some special cases, see in particular Section~\ref{sec:rectdegen}.

%\addtocounter{ex}{-1}
\begin{ex}[cont'd]
In Example~\ref{ex:hexa}, we have conveniently chosen all variables facing outwards at crossings of $\mathcal D$. Let us pick one, say $M_{1,2}$, as well as its conjugate variable $M_{2,4}$, and perform the substitution $M_{1,2}\to \epsilon\,M_{1,2}$, $M_{2,4}\to \epsilon^{-1} M_{2,4}$ (as we shall see, a one-step degeneration is enough in this case). Only the first component is affected by the degeneration, the other components being linear. Note that the given equations of $L_{\rho,\linkpattern[shape=circle,small]{1/2,3/4,5/6}}$ are {\em not} Gr\"obner, which means that their naive $\epsilon=0$ limit $\{M_{4,6}M_{6,2}=M_{6,2}M_{2,4}=M_{2,4}M_{4,6}=0\}$ is not enough to generate the ideal at the special fiber. Instead, we find
\begin{gather*}
L^{\epsilon=0}_{\rho,\linkpattern[shape=circle,small]{1/2,3/4,5/6}}=
\{M_{4,6}M_{6,2}=M_{6,2}M_{2,4}=M_{2,4}M_{4,6}=0,\ M_{3,4}M_{4,6}=M_{5,6}M_{6,2}=M_{1,2}M_{2,4} \}\\
\hphantom{L^{\epsilon=0}_{\rho,\linkpattern[shape=circle,small]{1/2,3/4,5/6}}}{} =
\{M_{2,4}=M_{4,6}=M_{5,6}=0\} \cup \{M_{2,4}=M_{3,4}=M_{6,2}=0\} \\
\hphantom{L^{\epsilon=0}_{\rho,\linkpattern[shape=circle,small]{1/2,3/4,5/6}}=}{}
\cup \{M_{1,2}=M_{4,6}=M_{6,2}=0\} \cup \{M_{2,4}=M_{4,6}=M_{6,2}=0\},
\end{gather*}
where the decomposition into irreducible components matches the first four loop configurations in~\eqref{eq:hexa1}. The fourth one appears with multiplicity~$2$. In $K$-theory, in order to see the coefficient $1+t$ appear in front of the fourth component, more work is needed: one must degenerate the whole of $L_\rho$ as before, then introduce the variable $\Phi=M_{1,2}M_{2,4}$, and revlex~$\Phi$ as well. We do not write the details here, and simply state the result: all primary components will have the extra equation $\Phi=0$ except the one corresponding to the configuration with a closed loop, which will only have $\Phi^2=0$, thus contributing an extra $1+t$ (taking into account $\wt_K(\Phi)=t$).\footnote{See Section~\ref{sec:rectdegen} for a more general example, from which it should be clear that there is exactly one such variable $\Phi$ (``flux'') satisfying $\Phi^2=0$ but not $\Phi=0$ per closed loop.} Finally, the three crossing loop configurations in~\eqref{eq:hexa1} are in one-to-one correspondence with the codimension $1$ intersections of these components.

Now assume that we had considered instead the other diagram of $\rho$:
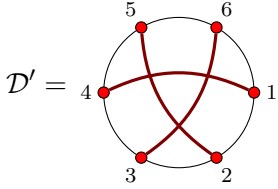
\begin{gather*}
\mathcal D'
=\linkpattern[edgecolor=colb,unit=1cm,shape=circle,squareness=0.1,size=6]{1/1/{[bend left=25] to (v4);},2/2/{[bend right=25] to (v5);},3/3/{[bend left=25] to (v6);}}%a bit clumsy
\end{gather*}
The degeneration using any of the crossings would force us to take another set of variables, say
$\{M_{2,3},M_{3,5},M_{4,5},M_{5,1},M_{6,1},M_{1,3}\}$. The irreducible components of $L_\rho$ would look different in these variables:
it is now the second component that has nontrivial equations
\begin{gather*}
L_{\rho,\linkpattern[shape=circle,small]{1/6,2/3,4/5}}=\{M_{4,5}+M_{1,3}M_{3,5}=M_{6,1}+M_{3,5}M_{5,1}=M_{2,3}+M_{5,1}M_{1,3}=0\},
\end{gather*}
while all other components are linear. This is the only component that has a nontrival de\-ge\-neration, corresponding to the multiple loop configurations of~$\mathcal D'$ with this connectivity:
\tikzset{hexa/.pic={
\draw[bgplaq] (90:1) -- node[midway,/linkpattern/vertex] (6) {} (150:1) -- node[midway,/linkpattern/vertex] (1) {} (210:1) -- node[midway,/linkpattern/vertex] (2) {} (270:1) -- node[midway,/linkpattern/vertex] (3) {} (330:1) -- node[midway,/linkpattern/vertex] (4) {} (30:1) -- node[midway,/linkpattern/vertex] (5) {} (90:1);
\draw (0,0) -- coordinate[midway] (a) (90:1) (0,0) -- coordinate[midway] (b) (-150:1) (0,0) -- coordinate[midway] (c) (-30:1);
}}
\begin{gather}
[L_{\rho,\linkpattern[shape=circle,small]{1/6,2/3,4/5}}]=(1-t z_1/z_2)(1-z_3/z_1)(1-t z_2/z_3)
 =
\tikz[baseline=-3pt,scale=0.7]{\pic[transform shape]{hexa}; \draw[/linkpattern/edge,bend left] (4) to (5) (6) to (1) [bend right] (2) to (b) [bend left] to (a) to (c) [bend right] to (3);}
+
\tikz[baseline=-3pt,scale=0.7]{\pic[transform shape]{hexa}; \draw[/linkpattern/edge,bend right] (4) to (c) [bend left] to (b) to (a) [bend right] to (5) [bend left] (6) to (1) (2) to (3);}
+
\tikz[baseline=-3pt,scale=0.7]{\pic[transform shape]{hexa}; \draw[/linkpattern/edge,bend left] (4) to (5) [bend right] (6) to (a) [bend left] to (c) to (b) [bend right] to (1) [bend left] (2) to (3);}
\nonumber\\
\hphantom{[L_{\rho,\linkpattern[shape=circle,small]{1/6,2/3,4/5}}]=}{}
+(1+t) \tikz[baseline=-3pt,scale=0.7]{\pic[transform shape]{hexa}; \draw[/linkpattern/edge,bend left] (4) to (5) (6) to (1) (2) to (3) (a) to (c) to (b) to (a);}
 - \tikz[baseline=-3pt,scale=0.7]{\pic[transform shape]{hexa}; \draw[/linkpattern/edge,bend left] (4) to (5) (6) to (1) [bend left=0] (2) to (c) [bend right=30] to (a) to (b) [bend left=0] to (3);}
-
\tikz[baseline=-3pt,scale=0.7]{\pic[transform shape]{hexa}; \draw[/linkpattern/edge] (4) to (a) [bend right] to (b) to (c) [bend left=0] to (5) [bend left=30] (6) to (1) (2) to (3);}
-
\tikz[baseline=-3pt,scale=0.7]{\pic[transform shape]{hexa}; \draw[/linkpattern/edge,bend left] (4) to (5) [bend left=0] (6) to (b) [bend right=30] to (c) to (a) [bend right=0] to (1) [bend left=30] (2) to (3);}\, .\label{eq:hexa2b}
\end{gather}

The equality of the r.h.s.\ of \eqref{eq:hexa2} and \eqref{eq:hexa2b} (and similarly for the other components) is nothing but a more explicit form of the Yang--Baxter equation in terms of loop model configurations. Geometrically, it arises as a consequence of the invariance of $K$-classes along flat families.
\end{ex}

\section{Rectangular domains and conormal matrix Schubert varieties}\label{sec:rect}
In this section, we use the same notations as in Section~\ref{sec:fd}. One may find many interesting varieties among the $L_{\rho,\pi}$ defined there. We now study in more detail one particular family of examples, which corresponds to partition functions on {\em rectangular domains}. We shall then reconnect to the results of Section~\ref{sec:Rmat}.

\subsection{Rectangular domains}\label{sec:rect0}
We assume that $K=k+n$, $N=2(k+n)$, where $k$ and $n$ are two positive integers, and choose~$\rho$ to be of the form:
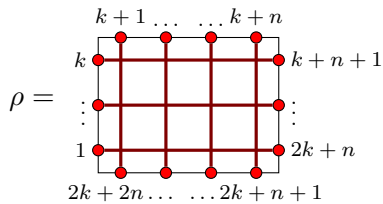
\begin{gather*}
\rho=\linkpattern[centered,edgecolor=colb,shape=rectangle,height=3,width=4,alias=false,numbering={1/1,2/\vdots,3/k,4/k+1,5/\ldots,6/\ldots,7/k+n,8/k+n+1,9/\vdots,10/2k+n,11/{\ldots 2k+n+1},12/,13/,14/{2k+2n\ldots}}]{1/10,2/9,3/8,4/14,5/13,6/12,7/11}
\end{gather*}
where we deformed the circle into a rectangle to make the structure of $\rho$ more obvious. Of course, the labelling will be momentarily redefined to something more sensible, namely matrix row/column.

The link pattern $\rho$ belongs to an important class, namely, ``triangle-free'' link patterns, i.e., for which one can never apply move~\eqref{eq:braid}, so that they possess a unique diagram. The description of $V_\rho$ and $L_\rho$ is in that case simpler, as we shall discuss now.

First, around each vertex, variables corresponding to opposite angles are opposite of each other, i.e., $M_{i,j}=-M_{\rho(i),\rho(j)}$ for all $\cyc{i<j<\rho(i)<\rho(j)}$. This means that up to a~choice of sign, $V_\rho$ possesses a canonical choice of coordinates, say the union over all crossings of $(M_{i,j},M_{j,\rho(i)})$. Introduce modified labelling for these coordinates, which in the present case are
\begin{gather}\label{eq:relabel}
M_{i,j}=Q_{k+1-i,j-k},\qquad M_{j,\rho(i)}=P_{k+1-i,j-k}, \qquad 1\le i\le k,\qquad k+1\le j\le k+n
\end{gather}
for pairs of $k\times n$ matrices $Q$ and $P$, making $V_\rho$ isomorphic to $T^\ast \Mat{k,n}$ (where $Q$ parameterizes the base and $P$ the fiber\footnote{It may seem more natural to consider pairs of matrices whose shape is transpose of each other, say $\big(Q^T,P\big)$ with our notations. Giving $P$ and $Q$ the same shape makes sense when a torus (here $T_\rho$) acts, fixing privileged coordinate subspaces; it should not obscure the fact that $T_\rho$ acts {\em differently} on $P$ and $Q$, as will be discussed in the next section.}). These satisfy canonical Poisson brackets
\begin{gather*}
\{ Q_{i,j}, P_{i',j'} \} =\delta_{i,i'}\delta_{j,j'},\qquad \{ Q_{i,j}, Q_{i',j'} \} = \{ P_{i,j}, P_{i',j'} \} =0.
\end{gather*}

Secondly, the equations of $L_\rho$ are quite easy to write explicitly. One can start directly from the definition at the start of Section~\ref{sec:lag}, but here we prefer to use the alternative characterization of Lemma~\ref{lem:slice2} as a slice of an orbital scheme. Starting from the slice \eqref{eq:stars} and keeping only the upper triangle (red entries), we obtain the shape of~$M$:
\begin{gather}\label{eq:Mrect}
\renewcommand{\kbldelim}{(}
\renewcommand{\kbrdelim}{)}
M=
\kbordermatrix{
&k\text{\tiny{ (left)}}&n\text{\tiny{ (top)}}&k\text{\tiny{ (right)}}&n\text{\tiny{ (bottom)}}\\
k\text{\tiny{ (left)}}& -/{}\mu_>{}/ & /{}Q& /&0\\
n\text{\tiny{ (top)}}& & \nu_< & P^T & /\\
k\text{\tiny{ (right)}}& & & -\mu_< & -Q{}/\\
n\text{\tiny{ (bottom)}}& & & & /{}\nu_>{}/ \\
},
\end{gather}
where $/$ is the shorthand notation for the antidiagonal matrices with $1$'s on the antidiagonal. $\mu_<$, $\mu_>$, $\nu_<$, $\nu_>$ are as yet unknown entries, which we know according to Lemma~\ref{lem:slice} must be expressible in terms of $P$ and $Q$; the minus signs and $/$s in their definition are for convenience. Note that $\mu_<$ and $\nu_<$ (resp.\ $\mu_>$ and $\nu_>$) are strict upper (resp.\ lower) triangular.

We now impose $M^2=0$; we find the equations
\begin{gather}
 \mu_<^2=\mu_>^2=\nu_<^2=\nu_>^2=0, \label{eq:321pre1} \\
 Q \nu_<-\mu_> Q = 0, \label{eq:321pre2} \\
 Q \nu_> - \mu_< Q =0, \label{eq:321pre3} \\
 P^T \mu_<-\nu_<P^T=0, \label{eq:321pre4} \\
 \mu_<+\mu_> = Q P^T, \label{eq:321pre5} \\
 \nu_<+\nu_> = P^T Q. \label{eq:321pre6}
\end{gather}
Denote $\mu=QP^T$, $\nu=P^TQ$;\footnote{These matrices will reappear in Section~\ref{sec:back} as moment maps for $\operatorname{GL}_k$ and $\operatorname{GL}_n$, respectively.} then the last two equations say that $\mu_<$ is the strict lower part of~$\mu$, and similarly for $\mu_>$, $\nu_<$, $\nu_>$; and that the diagonal parts of $\mu$ and $\nu$ vanish. Substituting this back into the other equations, and noting that~\eqref{eq:321pre2} and~\eqref{eq:321pre3} are equivalent, and that~\eqref{eq:321pre1} are consequences of say~\eqref{eq:321pre2} and~\eqref{eq:321pre4}, we finally obtain pairs of {\em cubic} equations satisfied by~$P$ and~$Q$ (plus quadratic equations for the diagonals of~$\mu$ and~$\nu$):
\begin{gather}\label{eq:321a}
\left(\sum_{i'=1}^{i-1}\sum_{j'=j}^{n}-\sum_{i'=i}^k\sum_{j'=1}^{j-1}\right) P_{i',j'}Q_{i',j}Q_{i,j'} =0,
\qquad i=1,\ldots,k,\qquad j=1,\ldots,n,
\\\label{eq:321b}
\left(\sum_{i'=1}^{i-1}\sum_{j'=1}^{j-1}-\sum_{i'=i}^k\sum_{j'=j}^{n}\right) P_{i',j}P_{i,j'}Q_{i',j'} =0,
\qquad i=1,\ldots,k,\qquad j=1,\ldots,n,
\\\label{eq:321c}
\sum_{i=1}^k P_{i,j} Q_{i,j} =0,\qquad j=1,\ldots,n,
\\\label{eq:321d}
\sum_{j=1}^n P_{i,j} Q_{i,j} =0,\qquad i=1,\ldots,k.
\end{gather}

As explained in Section~\ref{sec:irrcomp}, we can even obtain the equations for each $L_{\rho,\pi}$ by further imposing rank equations on~$M$. In particular, since the $L_{\rho,\pi}$, for $\pi\in\LP_{K,N}^{(K)}$, are Lagrangian and conical in the fiber (i.e., invariant under the torus that corresponds to the specialization $z_i=0$), they are conormal varieties of certain varieties inside $\Mat{k,n}$ whose defining equations we can determine:
\begin{itemize}\itemsep=0pt
\item {\em North-west/south-east rank conditions:} For each top vertex connected to its right nearest neighbor and each left vertex connected to its nearest neighbor below (resp.\ for each bottom vertex connected to its left nearest neighbor, and each right vertex connected to its nearest neighbor above), the rank of the north-west (resp.\ south-east) submatrix of $Q$ corresponding to the rectangle they delimit is less or equal to the number of arcs in that rectangle.
\item {\em Horizontal/vertical strip rank conditions:} For each pairing of neighbors on the left side and each pairing of neighbors on the right side (resp.\ for each pairing of neighbors on the top side and each pairing of neighbors on the bottom side), compute one half of (the Manhattan distance~-- horizontal plus vertical distance -- of their midpoints plus the number of arcs crossing the line joining them). If this number is negative, the variety is empty. Otherwise, and if the left midpoint is strictly below the right midpoint (resp.\ the top midpoint is strictly to the right of the bottom midpoint), the rank of $Q$ in the strip they delimit is less or equal to that number.
\end{itemize}
\begin{ex}Here are three rank conditions on submatrices of $Q$ for a particular link pattern:
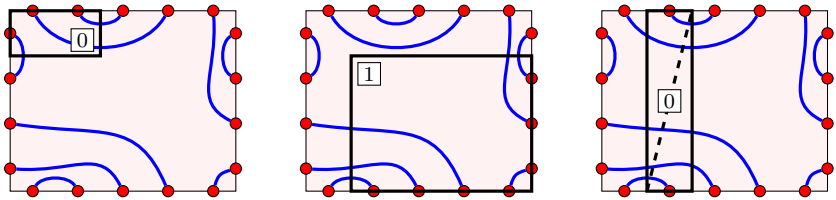
\begin{gather*}
\begin{tikzpicture}[/linkpattern/every linkpattern]
\linkpattern[shape=rectangle,boxed,numbered=false,inverted=false,height=4,width=5]{1/2,3/8,4/7,5/6,9/9',8'/5',7'/6',4'/1',3'/2'}
\draw[very thick] (0,-1) rectangle (2,0);
\node[rectangle,fill=white,draw=black,minimum size=0.3cm,inner sep=0pt] at (1.8-0.2,3.2-4+0.15) {$\ss0$};
\end{tikzpicture}
\qquad
\begin{tikzpicture}[/linkpattern/every linkpattern]
\linkpattern[shape=rectangle,boxed,numbered=false,inverted=false,height=4,width=5]{1/2,3/8,4/7,5/6,9/9',8'/5',7'/6',4'/1',3'/2'}
\draw[very thick] (1,-4) rectangle (5,-1);
\node[rectangle,fill=white,draw=black,minimum size=0.3cm,inner sep=0pt] at (1.2+0.2,2.8-4-0.2) {$\ss1$};
\end{tikzpicture}
\qquad
\begin{tikzpicture}[/linkpattern/every linkpattern]
\linkpattern[shape=rectangle,boxed,numbered=false,inverted=false,height=4,width=5]{1/2,3/8,4/7,5/6,9/9',8'/5',7'/6',4'/1',3'/2'}
\draw[very thick] (1,-4) rectangle (2,0); \draw[very thick,dashed] (1,-4) -- (2,0);
\node[rectangle,fill=white,draw=black,minimum size=0.3cm,inner sep=0pt] at (1.5,2-4) {$\ss0$};
\end{tikzpicture}
\end{gather*}
\end{ex}

Among the varieties obtained this way, and whose conormal varieties are therefore certain~$L_{\rho,\pi}$, we can recognize quite a few, including all {\em matrix Schubert varieties of 321-avoiding permutations}, certain Fomin--Zelevinsky double Bruhat cells, etc.

\subsection{Back to conormal Schubert varieties}\label{sec:back}
We now try to reconnect {\em in a different way than in Section~{\rm \ref{sec:fd}}} what precedes to the cotangent bundle of the Grassmannian. Recall that $T^\ast \operatorname{Gr}_{k,n}$, viewed as the Nakajima variety associated to the quiver $A_1$, is a symplectic quotient of $T^\ast \Mat{k,n}$ by the action of $\operatorname{GL}_k$; more precisely, using identical notations as in the previous section, and taking into account the stability condition, we have
\begin{gather*}
T^\ast \operatorname{Gr}_{k,n} \cong \big\{(Q,P)\in T^\ast \Mat{k,n}\colon \rk(Q)=k,\, PQ^T=0 \big\}/\operatorname{GL}_k.
\end{gather*}
The correspondence with the previous parameterization \eqref{eq:defGr} is $V=\Im Q^T$, $u=Q^T P$.

There are now various tori acting on our spaces:
\begin{itemize}\itemsep=0pt
\item On $T^\ast \operatorname{Gr}_{k,n}$ acts $T_n$, which is the Cartan torus of $\operatorname{GL}_n$ times $\CC^\times$. Recall that the weights of the entries of $u$ are given by~\eqref{eq:weiK}, which we rewrite here:
\begin{gather*}
\wt_K(u_{i,j})=t z_i^{-1}z_j,\qquad i,j=1,\ldots,n.
\end{gather*}

\item On $T^\ast\Mat{k,n}$, we have in the previous section the torus $T_\rho$ acting, a subgroup of $T_N$, the Cartan torus of $\operatorname{GL}_N$ times $\CC^\times$. If we use upper case for the weights of $\operatorname{GL}_N$, i.e., $Z_1,\ldots,Z_N$ (instead of $z_1,\ldots,z_N$ as before), then recall that the weights of $T_\rho$ are given by specializing $Z_{\rho(i)}=t Z_i$, $i<\rho(i)$; and the weights of the entries of $Q$, $P$ are, taking into account~\eqref{eq:torus} and the relabelling~\eqref{eq:relabel}:
\begin{gather*}
\wt_K(Q_{i,j}) =t Z_i Z_{j+k}^{-1}, \qquad \wt_K(P_{i,j}) =Z_i^{-1}Z_{j+k}.
\end{gather*}
\end{itemize}
If we try to identify the two actions via $u=Q^TP$, we find agreement on condition that
\begin{gather*}
z_j=Z_{j+k},\qquad j=1,\ldots,n.
\end{gather*}
$Z_1,\ldots,Z_k$ remain free at this stage, but it will be convenient to relabel them as well
\begin{gather*}
y_i=t\,Z_{k+1-i}=Z_{k+n+i},\qquad i=1,\ldots,k.
\end{gather*}
The weights of the entries of $Q$ and $P$ can now be rewritten as
\begin{gather*}
\wt_K(Q_{i,j}) =y_i z_j^{-1}, \qquad \wt_K(P_{i,j}) =t y_i^{-1}z_j,
\end{gather*}
which are just the weights of the natural $\operatorname{GL}_k\times \operatorname{GL}_n\times\CC^\times$ action on $T^\ast \Mat{k,n}$.

Now we discuss the various equivariant cohomology theories. In order to avoid too much repetition, we shall go straight to $K$-theory. We have
\begin{gather*}
K_{T_n}(T^\ast \operatorname{Gr}_{k,n}) \cong K_{\operatorname{GL}_k\times T_n} \big\{(Q,P)\in T^\ast \Mat{k,n}\colon
\, \rk(Q)=k,\, PQ^T=0\big\} \\
\hphantom{K_{T_n}(T^\ast \operatorname{Gr}_{k,n})}{} \twoheadleftarrow K_{\operatorname{GL}_k\times T_n} \big\{(Q,P)\in T^\ast \Mat{k,n}\colon
 PQ^T=0\big\},
\end{gather*}
where the $\twoheadleftarrow$ map is the pullback of the embedding. The last space is equivariantly contractible, and therefore its localized
equivariant $K$-theory is
\begin{gather*}
K_{\operatorname{GL}_k\times T_n}(\cdot) = K_{T_k\times T_n}(\cdot)^{\mathcal S_k} \cong \QQ(y_1,\ldots,y_k,z_1,\ldots,z_n,t)^{\mathcal S_k},
\end{gather*}
where $\mathcal S_k$ acts by permuting the $y_i$. This gives a presentation of $K_{T_n}(T^\ast \operatorname{Gr}_{k,n})$ as the quotient of the ring above by the common kernel of all restriction maps $|_I$ to fixed points, namely
\begin{align*}
|_I\colon \ \QQ(y_1,\ldots,y_k,z_1,\ldots,z_n,t)^{\mathcal S_k} &\to \QQ(z_1,\ldots,z_n,t),\\
f(y_1,\ldots,y_k,z_1,\ldots,z_n,t) &\mapsto f(z_{I_1},\ldots,z_{I_k},z_1,\ldots,z_n,t)
\end{align*}
for every $k$-subset $I=\{I_1,\ldots,I_k\}$ of $\{1,\ldots,n\}$.

Also recall that the union of conormal varieties $\CS_I$ is defined by the vanishing of the upper triangle of~$u$, or equivalently
\begin{gather*}
\bigcup_I \CS_I \cong \big\{(Q,P)\in T^\ast \Mat{k,n}\colon \rk(Q)=k,\, PQ^T=0,\, \big(Q^TP\big)_\le =0 \big\}/\operatorname{GL}_k.
\end{gather*}
Now the remarkable fact is that the equations $PQ^T=0$ and $\big(Q^T P\big)_\le=0$ imply \eqref{eq:321pre1}--\eqref{eq:321pre4} (and therefore \eqref{eq:321a}--\eqref{eq:321d}), as can be checked directly, recalling that $PQ^T=\mu^T$ and $u=Q^TP=\nu^T$.

This means that $\bigcup_I \CS_I$ is a subscheme of $(L_\rho-\{\rk(Q)<k\})/\operatorname{GL}_k$, for $\rho$ corresponding to a rectangular $k\times n$ domain. Since both schemes are equidimensional of the same dimension, it means each $\CS_I$ is a certain $(L_{\rho,\pi}-\{\rk(Q)<k\})/\operatorname{GL}_k$.

We now identify which irreducible components $L_{\rho,\pi}$ are related to each $\CS_I$ via this correspondence.
\begin{itemize}\itemsep=0pt
\item The matrix equation $P Q^T = 0$ (moment map condition for $\operatorname{GL}_k$) consists of three equations: the diagonal part (moment map condition for its Cartan torus) is nothing but \eqref{eq:321c}; and its upper triangle (resp.\ lower triangle) (moment map conditions for unipotent subgroups) implies that there cannot be any pairings between vertices on the left side (resp.\ right side) in $\pi$ (by using the linear equations among the rank equations of \eqref{eq:orb}, and applying them to~\eqref{eq:Mrect}).
\item Similarly, the equation $\big(Q^T P\big)_\le =0$ consists of: the diagonal part (moment map condition for the Cartan torus of~$\operatorname{GL}_n$) which is nothing but~\eqref{eq:321d}; and its upper triangle (moment map condition for the unipotent subgroup of $\operatorname{GL}_n$) which implies that there cannot be any pairing between vertices on the bottom side in~$\pi$.
\item Finally, the rank condition $\rk(Q)=k$ means that we should remove components satisfying $\rk(Q)<k$ (vanishing of $k\times k$ minors of~$Q$). On any component $L_{\rho,\pi}$, we have the inequality (rank equation from \eqref{eq:orb} applied to the $n\times k$ bottom rows of~\eqref{eq:Mrect})
\begin{gather*}
\rk(Q)\le \{\text{number of pairings between bottom and right sides in $\pi$}\}.
\end{gather*}
This means that we should only keep the components such that all $k$ vertices on the right side are connected to the bottom.
\end{itemize}
We now claim that such $\pi$ are in bijection with $\LP_{k,n}$: the map simply consists in keeping track only of the pairings among top vertices, the top vertices connected to the left (resp.\ bottom) boundary being marked as connected to left (resp.\ bottom) infinity. This is very similar to the truncation procedure of Section~\ref{sec:crossalg}; the only difference is that we have here $k$ additional ``spectator'' arcs connecting bottom and top. Furthermore, this finally explains, as promised, the asymmetry between left infinity and bottom infinity. We denote by $\phi$ the inverse injective map from $\LP_{k,n}$ to $\LP_{K,N}^{(K)}$; and we say that a loop configuration has ``top-connectivity'' $\pi$ when its connectivity is $\phi(\pi)$.

Finally, $L_{\rho,\phi(\pi)}$ is the conormal variety of a variety whose equations were described in Section~\ref{sec:rect0}; only north-west rank equations appear because of the particular form of $\phi(\pi)$, and we recognize the equations for {\em matrix Schubert varieties} \cite[Definition~1.3.2]{KM-Schubert} (in the case of the Grassmannian permutation associated to $\cl(\pi)$). We immediately conclude that
\begin{gather}\label{eq:identif0}
\CS_{\cl(\pi)} \cong (L_{\rho,\phi(\pi)}-\{\rk(Q)<k\})/\operatorname{GL}_k.
\end{gather}

In fact, more generally, we have
\begin{gather}\label{eq:identif}
X_\pi \cong (L_{\rho,\phi(\pi)}-\{\rk(Q)<k\})/\operatorname{GL}_k,\qquad \pi\in\CLP_{k,n},
\end{gather}
where we extend $\phi\colon \CLP_{k,n}\to\CLP_{K,N}^{(K)}$ to crossing link patterns in the obvious way: it reintroduces the various missing pairings
without introducing any extra crossings, e.g.,
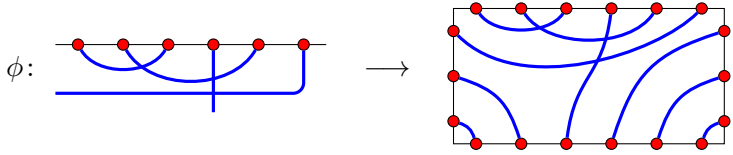
\begin{gather*}
\phi\colon \
\linkpattern[centered,height=1.5,alias=false,numbered=false%numbering={1/1,2/,3/\ldots,4/,5/,6/n}
]{1/3,2/5,4/4,6/6/\toleftx{2.7}}
\quad\longrightarrow\quad
%\tikz[baseline=(current bounding box.center)]{\linkpattern[numbered=false,shape=rectangle,width=6,height=3,squareness=0.2,looseness=0.2]{4/6,5/8,9/3,2/17,1/18,12/13,11/14,10/15,7/16}}
\linkpattern[centered,numbered=false,shape=rectangle,width=6,height=3,squareness=0.3]{4/6,5/8,9/9/{.. controls +(-1,1) and +(1.5,1.5) .. (v3);},2/17,1/18,12/13,11/14,10/15,7/16}
\end{gather*}

At this stage, we have formed a loop (pun intended) in our reasoning: starting from the cotangent bundle of the Grassmannian $T^\ast \operatorname{Gr}_{K,N}$, we have pushed forwarded our various classes to its affinization, namely the nilpotent cone $\mathcal N_{K,N}$. But then out of this cone we have sliced a~space~$V_\rho$, which up to taking an open set and dividing out by the action of $\operatorname{GL}_k$, is nothing but a smaller cotangent bundle $T^\ast \operatorname{Gr}_{k,n}$.

This has various important consequences. The most interesting one is that it provides a~formula for the classes of the~$X_\pi$ (and in particular the conormal Schubert varieties) as a function of the $y_1,\ldots,y_k$ (which are one minus the Chern roots of $T^\ast \operatorname{Gr}_{k,n}$)
and of the equivariant parameters $z_1,\ldots,z_n,t$. Indeed, putting together the results of this section, and taking into account that $\{ (Q,P)\in T^\ast \Mat{k,n}\colon PQ^T=0\}$ is a complete intersection, so that its $K$-class is easily computed, we find:
\begin{cor}\label{cor:struc}
The class of $($the structure sheaf of$)$ $X_\pi$ in
\begin{gather*}
K_{T_n}(T^\ast \operatorname{Gr}_{k,n})\cong \QQ(y_1,\ldots,y_k,z_1,\ldots,z_n,t)/\bigcap_I \Ker |_I
\end{gather*}
is given by the partition function of the crossing loop model on a rectangular domain with top-connectivity $\pi$,
divided by the class of $\{PQ^T=0\}$:
\begin{gather*}
[X_\pi]=\prod_{i,j=1}^k (1-t y_i/y_j)^{-1}
\sum_{\substack{\text{\rm loop configurations on $k\times n$}\\\text{\rm with top-connectivity $\pi$}}}
(-1)^{|\text{\rm removed crossings}|}(1+t)^{|\text{\rm loops}|}
\\
\hphantom{[X_\pi]=}{}\times \prod_{i=1}^k \prod_{j=1}^n
\begin{cases}
1-t z_j/y_i
&\tikz[baseline=-3pt]{\plaq{a}}
\\[2mm]
1-y_i/z_j
&\tikz[baseline=-3pt]{\plaq{b}}
\\[2mm]
(1-t\,z_j/y_i)(1-y_i/z_j)
&\tikz[baseline=-3pt]{\plaq{c}}
\end{cases}
\end{gather*}
\end{cor}
This in turn implies Claim~\ref{thm:cr}, which we derived this corollary from in the first place via~\eqref{eq:exchcr}.

We can derive the corresponding cohomology statement by taking the appropriate limit, but as already mentioned we now focus on $K$-theory.

We turn to coherent sheaves. Given $\pi\in \LP_{k,n}$, we need to compare the ``sliced'' sheaves $\sigma_{\phi(\pi),\rho}$ with $\sigma_\pi$ itself. We find that not only the supports match, but also the divisors \eqref{eq:div} match via the identification~\eqref{eq:identif}. We now compare their coefficients. Two arcs are neighboring in $\pi$ iff they are in $\phi(\pi)$, and their depth is shifted by exactly~$k$, corresponding to the extra $k$ arcs starting on the left side of $\phi(\pi)$ (equivalently, one can easily compute $a_{i+k}(\phi(\pi))=k+a_i(\pi)$, $a_{n+k}(\phi(\pi))=k$).

To $\sigma_{\phi(\pi),\rho}$ is associated via \eqref{eq:identif0} a sheaf on $\CS_{\cl(\pi)}$ denoted $\tilde\sigma_{\phi(\pi),\rho}$. The reasoning above show that the divisors associated to the sheaves $\tilde\sigma_{\phi(\pi),\rho}$ and $\sigma_\pi$ differ by $k$ times the sum over neighboring $\alpha$, $\beta$ of all divisors $X_{f_{\alpha,\beta}\pi}$. Now this sum defines nothing but the line bundle $O(1)$ (the corresponding statement on the base $S_{\cl(\pi)}$ is well-known), so that we have (nonequivariantly)
\begin{gather*}
\sigma_\pi=O(-k)\otimes \tilde\sigma_{\phi(\pi),\rho}.
\end{gather*}
The $K$-class of $O(-k)$ is $\prod_{i=1}^k y_i^k$. The classes of $\sigma_{\phi(\pi),\rho}$ and $\tilde\sigma_{\phi(\pi),\rho}$ are equal up to identification of equivariant parameters $y_i$ with the corresponding elements $y_i\in K_T(T^\ast \operatorname{Gr}_{k,n})$. Furthermore, we have to take into account the equivariant structure, cf.~\eqref{eq:weinorm}. This results in an additional monomial; a careful computation leads to
\begin{cor}\label{cor:sqrt}
The class of $\sigma_\pi$ in
\begin{gather*}
K_{T'_n}(T^\ast \operatorname{Gr}_{k,n})\cong \QQ\big(y_1,\ldots,y_k,z_1,\ldots,z_n,t^{1/2}\big)/\bigcap_I \Ker |_I
\end{gather*}
is given $($up to a monomial$)$ by the partition function of the noncrossing loop model on a rectangular domain with top-connectivity $\pi$, divided by the class of $\{PQ^T=0\}$:
\begin{gather*}
[\sigma_\pi]=t^{\frac{1}{4}k(2n+k+1)}\prod_{i=1}^k y_i^{k-n}\prod_{i,j=1}^k (1-t y_i/y_j)^{-1}
\sum_{\substack{\text{\rm loop configurations on $k\times n$}\\\text{\rm with top-connectivity $\pi$}}}
\big(t^{1/2}+t^{-1/2}\big)^{|\text{\rm loops}|}\\
\hphantom{[\sigma_\pi]=}{}\times
\prod_{i=1}^k \prod_{j=1}^n
\begin{cases}
t^{-1/2}y_i/z_j-t^{1/2}
&\tikz[baseline=-3pt]{\plaq{a}}\\[2mm]
1-y_i/z_j
&\tikz[baseline=-3pt]{\plaq{b}}
\end{cases}
\end{gather*}
\end{cor}

\begin{ex}Let us redo the case $k=1$, $n=2$, which was already investigated in Example~\ref{ex:base}. There is only one loop model configuration corresponding to each link pattern:
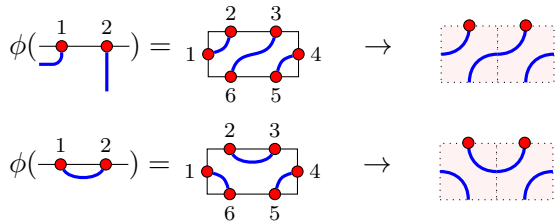
\begin{gather*}
\phi(
\linkpattern[centered]{1/1/\toleft,2/2}
) =
\linkpattern[centered,shape=rectangle,width=2,height=1]{1/2,3/6,4/5}
\quad\rightarrow\quad
\tikz[baseline=(current bounding box.center)]{
\node[loop]{\plaq{a}&\plaq{a}\\};
\node[/linkpattern/vertex] at (loop-1-1.north) {};
\node[/linkpattern/vertex] at (loop-1-2.north) {};
}
\\
\phi(
\linkpattern[centered]{1/2}
) =
\linkpattern[centered,shape=rectangle,width=2,height=1]{1/6,2/3,4/5}
\quad\rightarrow\quad
\tikz[baseline=(current bounding box.center)]{
\node[loop]{\plaq{b}&\plaq{a}\\};
\node[/linkpattern/vertex] at (loop-1-1.north) {};
\node[/linkpattern/vertex] at (loop-1-2.north) {};
}
\end{gather*}
We conclude that
\begin{gather*}
[\sigma_{\linkpattern[small]{1/1/\toleft,2/2}}] =t^{3/2}y_1^{-1}(1-t)^{-1} \big(t^{-1/2}y_1/z_1-t^{1/2}\big)\big(t^{-1/2}y_1/z_2-t^{1/2}\big),\\
[\sigma_{\linkpattern[small]{1/2}}]=t^{3/2}y_1^{-1}(1-t)^{-1} (1-y_1/z_1)\big(t^{-1/2}y_1/z_2-t^{1/2}\big).
\end{gather*}
It is not hard to check that specializing at $y_1=z_1,z_2$ reproduces the restriction to fixed points of Example~\ref{ex:base}.
\end{ex}

We can also easily deduce Anderson--Jantzen--Soergel--Billey-type formulae for restrictions of~$[X_\pi]$ and of~$[\sigma_\pi]$ to fixed points by specializing these partition functions to $y_i=z_{I_i}$, see \cite[Section~4.4]{moscowlectures} and \cite[Section~2.2]{artic69} for related computations.

We postpone to Section~\ref{sec:gorschub} further examples of application of Claims~\ref{cor:struc} and \ref{cor:sqrt}.

\subsection{Degeneration}\label{sec:rectdegen}
The conjectured degeneration of Section~\ref{sec:degen} can be made explicit in the rectangular case. Perform the substitution $Q_{i,j}\to \epsilon^{ij}Q_{i,j}$, $P_{i,j}\to \epsilon^{-ij} P_{i,j}$. We claim that the special fiber at $\epsilon=0$ is given by the equations (generating the ideal of lowest degree terms in $\epsilon$):
\begin{gather}\label{eq:degen321a}
\left(\sum_{i'=1}^{i-1} Q_{i',j}P_{i',j}-\sum_{j'=1}^{j-1} Q_{i,j'}P_{i,j'}\right) Q_{i,j} =0, \qquad i=1,\ldots,k,\qquad j=1,\ldots,n,
\\\label{eq:degen321b}
\left(\sum_{i'=i+1}^{k}Q_{i',j}P_{i',j}-\sum_{j'=1}^{j-1} Q_{i,j'}P_{i,j'}\right)P_{i,j} =0,\qquad i=1,\ldots,k,\qquad j=1,\ldots,n,
\\\label{eq:degen321c}
\sum_{i=1}^k P_{i,j} Q_{i,j} =0,\qquad j=1,\ldots,n,
\\\label{eq:degen321d}
\sum_{j=1}^n P_{i,j} Q_{i,j} =0, \qquad i=1,\ldots,k.
\end{gather}
These are the $\epsilon\to0$ limit of equations \eqref{eq:321a}--\eqref{eq:321b} up to small rearrangements (note that \eqref{eq:321c} and \eqref{eq:321d} stayed the same), implying a nontrivial Gr\"obner statement. This will be ultimately justified by the fact that the resulting degenerated scheme has the same cohomology class as the one of the original scheme as given by Claim~\ref{thm:fdcoh} (so that any further equation would ``decrease'' that class with an appropriate notion of positivity in weight space, leading to a contradiction).

Let us study these equations in more detail. Define {\em flux} variables $\Phi_e$ associated to each edge of the rectangular domain (i.e., of the dual map of $\rho$) as follows:
\begin{itemize}\itemsep=0pt
\item If $e$ is a boundary edge, then $\Phi_e=0$.
\item If $e$ is vertical, sitting between plaquettes at rows/columns $(i,j)$ and $(i,j+1)$, define
\begin{gather*}
\Phi_e=\sum_{j'\le j} Q_{i,j'} P_{i,j'} =-\sum_{j'\ge j+1} Q_{i,j'} P_{i,j'},
\end{gather*}
where all equalities are modulo \eqref{eq:degen321a}--\eqref{eq:degen321d}. One should think of it as an oriented flux from left to right across the edge.
\item
If $e$ is horizontal, sitting between plaquettes at rows/columns $(i,j)$ and $(i+1,j)$, define
\begin{gather*}
\Phi_e=\sum_{i'\le i} Q_{i',j} P_{i',j} =-\sum_{i'\ge i+1} Q_{i',j} P_{i',j}.
\end{gather*}
It is an oriented flux from bottom to top.
\end{itemize}

The denomination of flux is justified by the fact that at each plaquette $(i,j)$, from the very definition of $\Phi_e$, we have the conservation equation
\begin{gather}\label{eq:flux}
\tikz[baseline=-3pt]{
\setlength{\loopcellsize}{1.25cm}\plaq{}
\node[shape=isosceles triangle,shape border rotate=90,inner sep=2pt,fill=blue,label={above:$\Phi_{\text{N}}$}] at (plaq.north) {};
\node[shape=isosceles triangle,shape border rotate=90,inner sep=2pt,fill=blue,label={below:$\Phi_{\text{S}}$}] at (plaq.south) {};
\node[shape=isosceles triangle,shape border rotate=0,inner sep=2pt,fill=blue,label={right:$\Phi_{\text{E}}$}] at (plaq.east) {};
\node[shape=isosceles triangle,shape border rotate=0,inner sep=2pt,fill=blue,label={left:$\Phi_{\text{W}}$}] at (plaq.west) {};
}
\qquad
\Phi_W+\Phi_S=\Phi_E+\Phi_N.
\end{gather}

It follows directly from \eqref{eq:degen321a}--\eqref{eq:degen321d} that
\begin{gather*}
\Phi_e^2=0
\end{gather*}
for all edges $e$. In particular, if we are only interested in the {\em radical} of these equations, then we have $\Phi_e=0$, and therefore by taking differences, $Q_{i,j}P_{i,j}=0$ for all~$i$,~$j$. We conclude that the degenerated scheme has $2^{kn}$ geometric components, corresponding to the choice of $P_{i,j}=0$ (plaquette $\tikz[baseline=-3pt]{\plaq{a}}$) or $Q_{i,j}=0$ (plaquette $\tikz[baseline=-3pt]{\plaq{b}}$) at crossing $(i,j)$, as expected.

However, we want to go further and recover the multiplicities of the components in the unreduced degenerated scheme. For that, we rewrite \eqref{eq:degen321a}--\eqref{eq:degen321b} in terms of the fluxes:
\begin{gather*}
\tikz[baseline=-3pt]{
\setlength{\loopcellsize}{1.25cm}\plaq{}
\node at (plaq) {$\ss (i,j)$};
\node[shape=isosceles triangle,shape border rotate=90,inner sep=2pt,fill=blue,label={above:$\Phi_{\text{N}}$}] at (plaq.north) {};
\node[shape=isosceles triangle,shape border rotate=90,inner sep=2pt,fill=blue,label={below:$\Phi_{\text{S}}$}] at (plaq.south) {};
\node[shape=isosceles triangle,shape border rotate=0,inner sep=2pt,fill=blue,label={right:$\Phi_{\text{E}}$}] at (plaq.east) {};
\node[shape=isosceles triangle,shape border rotate=0,inner sep=2pt,fill=blue,label={left:$\Phi_{\text{W}}$}] at (plaq.west) {};
}
\qquad
\begin{aligned}
& (\Phi_W-\Phi_N)Q_{i,j}=(\Phi_E-\Phi_S)Q_{i,j}=0,
\\
& (\Phi_W+\Phi_S)P_{i,j}=(\Phi_E+\Phi_N)P_{i,j}=0
\end{aligned}
\end{gather*}
(the two formulations are equivalent modulo \eqref{eq:flux}).

Let us now consider the {\em primary} (i.e., unreduced) component corresponding to a given loop configuration. We can immediately simplify the equations above by using the fact that $Q_{i,j}$ is not in its radical for $\tikz[baseline=-3pt]{\plaq{a}}$, and $P_{i,j}$ is not in its radical for $\tikz[baseline=-3pt]{\plaq{b}}$:
\begin{gather*}
\tikz[baseline=-3pt]{
\setlength{\loopcellsize}{1.25cm}\plaq{a}
\node[shape=isosceles triangle,shape border rotate=90,inner sep=2pt,fill=blue,label={above:$\Phi_{\text{N}}$}] at (plaq.north) {};
\node[shape=isosceles triangle,shape border rotate=90,inner sep=2pt,fill=blue,label={below:$\Phi_{\text{S}}$}] at (plaq.south) {};
\node[shape=isosceles triangle,shape border rotate=0,inner sep=2pt,fill=blue,label={right:$\Phi_{\text{E}}$}] at (plaq.east) {};
\node[shape=isosceles triangle,shape border rotate=0,inner sep=2pt,fill=blue,label={left:$\Phi_{\text{W}}$}] at (plaq.west) {};
}
\qquad
\Phi_W-\Phi_N =\Phi_E-\Phi_S=0,
\\
\tikz[baseline=-3pt]{
\setlength{\loopcellsize}{1.25cm}\plaq{b}
\node[shape=isosceles triangle,shape border rotate=90,inner sep=2pt,fill=blue,label={above:$\Phi_{\text{N}}$}] at (plaq.north) {};
\node[shape=isosceles triangle,shape border rotate=90,inner sep=2pt,fill=blue,label={below:$\Phi_{\text{S}}$}] at (plaq.south) {};
\node[shape=isosceles triangle,shape border rotate=0,inner sep=2pt,fill=blue,label={right:$\Phi_{\text{E}}$}] at (plaq.east) {};
\node[shape=isosceles triangle,shape border rotate=0,inner sep=2pt,fill=blue,label={left:$\Phi_{\text{W}}$}] at (plaq.west) {};
}
\qquad
\Phi_W+\Phi_S =\Phi_E+\Phi_N=0.
\end{gather*}
We reach the important conclusion that fluxes are conserved {\em along loops} (i.e., blue lines on the picture). Hence fluxes along lines that connect to the boundary are zero (since the flux is zero at boundaries), whereas fluxes along closed loops remain nonzero (although their square is zero).

We can now bound from above the cohomology class of each degenerated primary component by its product of local weights times the contribution of the flux equations $\Phi_e^2=0$, that is the multiplicity $2^{|\text{loops}|}$. This is the correct result according to Claim~\ref{thm:fdcoh}, which proves the Gr\"obner statement.

Similarly, the $K$-theory formula of Claim~\ref{thm:fdcr} can be recovered by further revlex-ing the flux variables $\Phi_e$. We shall skip the details here. Presumably, something similar works for Claim~\ref{thm:fdnc}.

\subsection{Pipe dreams}
The loop configurations that are considered in this paper, especially the ones on rectangular domains, are very similar to so-called {\em pipe dreams} \cite{KM-Schubert}; the difference is that we allow three plaquettes per site~-- $\tikz[baseline=-3pt]{\plaq{a}}$, $\tikz[baseline=-3pt]{\plaq{b}}$ and $\tikz[baseline=-3pt]{\plaq{c}}$ -- whereas only $\tikz[baseline=-3pt]{\plaq{a}}$ and $\tikz[baseline=-3pt]{\plaq{c}}$ appear in pipe dreams. We now show how to recover pipe dreams as a special case of our loop configurations.

We keep the same setup as in Section~\ref{sec:rect0}, but now we set $k=n$, that is, we consider a square $n\times n$ domain. Given a permutation $w\in \mathcal S_n$, we choose the following crossing link pattern $\pi_w$ (in terms of the original labelling of vertices from $1$ to $N=4n$):
\begin{itemize}\itemsep=0pt
\item It pairs $i$ and $w_{n+1-i}+n$, $i=1,\ldots,n$.
\item It pairs $2n+i$ and $4n+1-i$, $i=1,\ldots,n$.
\end{itemize}
Using a row/column relabelling, we obtain something of the form
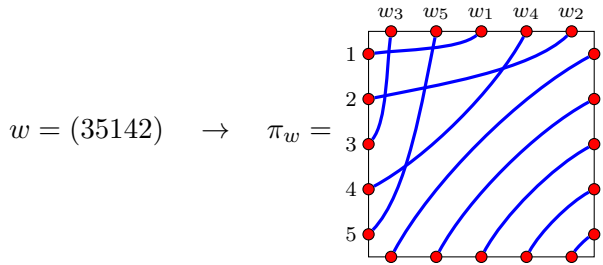
\begin{gather*}
w=(35142)\quad\to\quad
\pi_w= \linkpattern[centered,shape=rectangle,width=5,height=5,inverted=false,squareness=0.12,numbering={1/1,2/2,3/3,4/4,5/5,6/,7/,8/,9/,10/,11/,12/,13/,14/,15/,16/{w_2},17/{w_4},18/{w_1},19/{w_5},20/{w_3}}]{1/18,2/16,3/20,4/17,5/19,6/15,7/14,8/13,9/12,10/11}
\end{gather*}
where the north-west half-square reproduces the usual diagram of the permutation $w$.

The corresponding components $L_{\rho,\pi_w}$ all live in the subscheme of $L_{\rho}$ for which there are no pairings among vertices on the bottom and left or among vertices on the top and right, which is nothing but the subscheme $\{P=0\}$ (linear rank condition from \eqref{eq:orb} on the $(2k+n)\times(2k+n)$ south-west block of~\eqref{eq:Mrect}). Furthermore the north-west rank conditions for $Q\in L_{\rho,\pi_w}$ from Section~\ref{sec:rect0} are exactly the rank conditions for the matrix Schubert associated to $w$, cf.~\cite[Definition~1.3.2]{KM-Schubert}, so that $L_{\rho,\pi_w}$ is isomorphic to the latter ({\em not} its conormal variety!).

Now study the loop configurations associated to $\pi_w$. Because of the lack of pairings among left/bottom and among top/right, it is easy to see that the plaquette $\tikz[baseline=-3pt]{\plaq{b}}$ is forbidden, as expected. Furthermore, the whole south-west half of the square, diagonal included, is frozen to be $\tikz[baseline=-3pt]{\plaq{a}}$, because no crossings are allowed among bottom/right vertices. (Conventionally, this half is usually erased when drawing pipe dreams.) Furthermore, connectivity $\pi_w$ for the loop configuration is equivalent to connectivity $w$ for the pipe dream (noting that in both cases, if two lines cross multiple times, only the first crossing counts), e.g., with the same example as above,
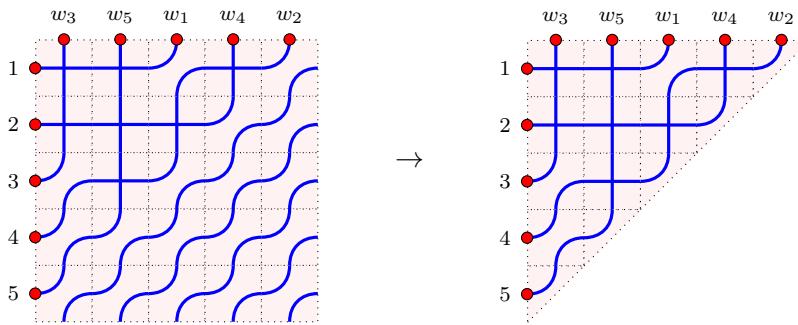
\begin{gather*}
\tikz[baseline=(current bounding box.center)]{
\node[loop]
{
\plaq{c}&\plaq{c}&\plaq{a}&\plaq{c}&\plaq{a}\\
\plaq{c}&\plaq{c}&\plaq{c}&\plaq{a}&\plaq{a}\\
\plaq{a}&\plaq{c}&\plaq{a}&\plaq{a}&\plaq{a}\\
\plaq{a}&\plaq{a}&\plaq{a}&\plaq{a}&\plaq{a}\\
\plaq{a}&\plaq{a}&\plaq{a}&\plaq{a}&\plaq{a}\\
};
\foreach\x/\y in {1/3,2/5,3/1,4/4,5/2} \node[/linkpattern/vertex,label={above:$\ss w_\y$}] at (loop-1-\x.north) {};
\foreach\x in {1,...,5} \node[/linkpattern/vertex,label={left:$\ss\x$}] at (loop-\x-1.west) {};
}
\qquad\to\qquad
\tikz[baseline=(current bounding box.center)]{
\node[loop]
{
\plaq{c}&\plaq{c}&\plaq{a}&\plaq{c}&\halfplaq\\
\plaq{c}&\plaq{c}&\plaq{c}&\halfplaq\\
\plaq{a}&\plaq{c}&\halfplaq\\
\plaq{a}&\halfplaq\\
\halfplaq\\
};
\foreach\x/\y in {1/3,2/5,3/1,4/4,5/2} \node[/linkpattern/vertex,label={above:$\ss w_\y$}] at (loop-1-\x.north) {};
\foreach\x in {1,...,5} \node[/linkpattern/vertex,label={left:$\ss\x$}] at (loop-\x-1.west) {};
}
\end{gather*}

Finally, in Claim~\ref{thm:fdcr}, applying the same relabelling of spectral parameters as in Section~\ref{sec:back}, we note that all weights have a common factor, so that
\begin{gather*}
[L_{\rho,\pi_w}]=\prod_{i,j=1}^n (1-t\,z_j/y_i)
\sum_{\substack{\text{pipe dreams}\\\text{with connectivity $w$}}}
(-1)^{|\text{removed crossings}|}
\prod_{i,j=1}^n
\begin{cases}
1
&\tikz[baseline=-3pt]{\plaq{a}}
\\[2mm]
1-y_i/z_j
&\tikz[baseline=-3pt]{\plaq{c}}
\end{cases}
\end{gather*}
The prefactor $\prod_{i,j=1}^n (1-t\,z_j/y_i)$ is nothing but the contribution of the equations $P=0$. The rest of the r.h.s.\ is the $K$-class of the matrix Schubert variety, i.e., $L_{\rho,\pi_w}$ embedded in $\Mat{n}$, and indeed we recognize the pipe dream formula for double Grothendieck polynomials (see \cite[Theorem~3.1]{FK-Groth} as well as \cite[Theorem~A]{KM-Schubert} and \cite{KM-Schubert2}).

Note that the Yang--Baxter equation already plays a prominent role in \cite{FK-Groth, FK-Schubert}; in fact, the ``algebra of projectors'' of~\cite{FK-Groth} is nothing but the subalgebra of $\CTL_n$ that is generated by the $f_i$, consistent with the fact that we are using plaquettes corresponding to the identity (\tikz[baseline=-3pt]{\plaq{a}}) and the $f_i$ (\tikz[baseline=-3pt]{\plaq{c}}), but not the one associated to the $e_i$ (\tikz[baseline=-3pt]{\plaq{b}}).

Furthermore, it is easy to see that the degeneration of Section~\ref{sec:rectdegen} reduces, in the case of link patterns $\pi_w$, to the degeneration of \cite{KM-Schubert} (in particular, the absense of loops makes the de\-ge\-neration reduced, hence a~Stanley--Reisner scheme). In that sense, the present study generalizes both integrable \cite{FK-Groth} and geometric \cite{KM-Schubert} aspects of $K$-classes of matrix Schubert varieties/Schubert varieties of the flag variety.

\section[Connection to the quantum Knizhnik--Zamolodchikov equation\\ and combinatorics]{Connection to the quantum Knizhnik--Zamolodchikov\\ equation and combinatorics}

We collect in this section various properties and interpretations of the $K$-theoretic quantities introduced in Section~\ref{sec:Rmat}, in particular the Laurent polynomials $\Psi^c_\pi$ and $\Psi^{nc}_\pi$.

\subsection[$q$KZ equation]{$\boldsymbol{q}$KZ equation}\label{sec:qKZ}

We assume in this section that $n=2k$. It is then known that the polynomials $\Psi^{nc}_\pi$ are related to the quantum Knizhnik--Zamolodchikov ($q$KZ) equation~\cite{artic34,artic69} (see also~\cite[Section~4]{hdr} for a~review, noting the correspondence of notations: $t^{1/2}=-q$). More precisely, if we redefine
\begin{gather*}
\tilde\Psi^{nc}_\pi := t^{-\frac{1}{4}k(5k-1)}\prod_{i=1}^n z_i^{i-1} \Psi_\pi^{nc}
\end{gather*}
(the power of $t$ is added for convenience), then the $\tilde\Psi^{nc}_\pi$ are polynomials (as opposed to Laurent polynomial) of degree $k(k-1)$ in the $z_i$, which collectively satisfy the system:
\begin{gather*}
\tilde\Psi^{nc}_\pi \tau_i = \sum_{\pi'\in \LP^{(m)}_{k,n}} \left(\frac{1-t\,z_{i+1}/z_i-t^{1/2}(1-z_{i+1}/z_i)e_i}{z_{i+1}/z_i-t}\right)_{\pi,\pi'}\tilde\Psi^{nc}_{\pi'},\qquad \pi\in \LP^{(m)}_{k,n},\\
\tilde\Psi^{nc}_\pi\rho = \big({-}t^{1/2}\big)^{3(k-1)} \tilde\Psi^{nc}_{r\pi},
\end{gather*}
where $\rho$ is the operator that permutes cyclically the $z_i$ according to $\rho(z_i)=z_{i+1}$ for $i<n$ and $\rho(z_n)=t^3 z_1$, and $r$ is the ``rotation'' of link patterns that amounts to relabelling the vertices $i\mapsto i+1$ in $\ZZ/n\ZZ$. This is the level $1$ $q$KZ system, which itself implies the usual $q$KZ equation~\cite{FR-qKZ}. The first equation is of course nothing but the exchange relation~\eqref{eq:exchnc}. The second equation is known as the {\em cyclicity} relation; its geometric interpretation is unknown. It is tantalizing that the $q$KZ equation also appears in~\cite{AO-ell, MO-qg}.

The $\tilde\Psi^{nc}_\pi$ are also known to satisfy a palindromy property:
\begin{gather}\label{eq:palin}
\tilde\Psi^{nc}_\pi(z_1,\ldots,z_n,t)= \prod_{i=1}^n z_i^{k-1}\tilde\Psi^{nc}_\pi(1/z_1,\ldots,1/z_n,1/t).
\end{gather}

We are particularly interested in the combinatorial interpretation of the specialization $z_i=1$, which in geometric terms means keeping only the equivariance w.r.t.\ scaling. General pro\-perties of one-variable Hilbert series imply that all such specializations have a zero of order $\codim_{\mathfrak n_-} \mathcal O_\pi = k(k-1)$ at $t=1$. We thus denote
\begin{gather*}
\Psi^{(1)}_\pi=\frac{\tilde\Psi_\pi^{nc}|_{z_i=1}}{(t^{1/2}-t^{-1/2})^{k(k-1)}}=t^{-\frac{1}{4}k(3k+1)} \frac{\Psi_\pi^{nc}|_{z_i=1}}{(1-t)^{k(k-1)}}.
\end{gather*}
According to \eqref{eq:palin}, we have $\Psi^{(1)}_\pi(1/t)=\Psi^{(1)}_\pi(t)$; it is therefore natural to use the loop weight $\beta=t^{1/2}+t^{-1/2}$ as a parameter, writing
\begin{gather*}
\Psi^{(1)}_\pi(t) = P_\pi(\beta).
\end{gather*}
$P_\pi(\beta)$ is a polynomial of a given parity in $\beta$. An important conjecture, made in \cite{artic41}, is that the coefficients of $P_\pi$ (as a polynomial in~$\beta$) are {\em positive} integers.

The positivity of these coefficients suggests that they should have some enumerative meaning. For some families of link patterns, it is known, and we shall discuss one such case in Section~\ref{sec:TSSCPP}. In general, we only have partial information:
\begin{itemize}\itemsep=0pt
\item The sum of these coefficients, i.e., the further specialization $P_\pi(\beta=1)$, is known to have a combinatorial meaning: it is the content of the famous Razumov--Stroganov conjecture~\cite{RS-conj}, proven in~\cite{CS-RS} (see also \cite[Section~1.4]{artic69} for a discussion in the present context). As mentioned in the introduction, $\beta=1$ corresponds to~$t$ nontrivial cubic root of unity, and it is unclear geometrically why such values should play a special role.
\item The sum $\sum_{\pi\in \LP^{(m)}_{k,n}} P_\pi(\beta)$ has an interpretation in terms of a weighted enumeration of totally symmetric self-complementary plane partitions~\cite{artic41}. Since the summation over $\pi$ is also unclear geometrically, we shall not develop this here.
\end{itemize}

\subsection{The Gorenstein case}\label{sec:goren}
The previous paragraph was concerned with the noncrossing loop model. No such connection with the $q$KZ equation appears in the crossing loop model. We now discuss cases where noncrossing and crossing loop models produce the same result.

\subsubsection{Gorenstein Schubert varieties}\label{sec:gorschub}
We temporarily return to general $k$ and $n$. For the purposes of this paragraph, we ignore the equivariant structure of our sheaves on $T^\ast \operatorname{Gr}_{k,n}$, which only contribute monomials in $t$ and the~$z$'s to their $K$-classes.

Assuming Conjecture~\ref{conj:CM} to hold, consider, for $\pi\in\LP_{k,n}$, the canonical sheaf of $\CS_{\cl(\pi)}$, which we denote~$\kappa_\pi$. By a simple duality argument, we can deduce from Corollary~\ref{cor:struc}
a~formula for its $K$-class:
\begin{cor}\label{cor:canon}
\begin{gather*}
[\kappa_\pi]=(-1)^{kn}\prod_{i=1}^k y_i^{-n}
\prod_{i,j=1}^k (1-t\,y_i/y_j)^{-1}
\sum_{\substack{\text{\rm loop configurations}\\\text{\rm on $k\times n$ with}\\ \text{\rm top-connectivity $\pi$}}}
(-1)^{|\text{\rm removed crossings}|}(1+t^{-1})^{|\text{\rm loops}|}
\\
\hphantom{[\kappa_\pi]=}{}\times \prod_{i=1}^k \prod_{j=1}^n
\begin{cases}
1-t^{-1}y_i/z_j &\tikz[baseline=-3pt]{\plaq{a}}
\\[2mm]
1-z_j/y_i &\tikz[baseline=-3pt]{\plaq{b}}
\\[2mm]
\big(1-t^{-1}y_i/z_j\big)(1-z_j/y_i) &\tikz[baseline=-3pt]{\plaq{c}}
\end{cases}
\end{gather*}
\end{cor}

This is {\em per se} not very interesting, but is to be compared with the next remark. It is not hard to show by explicit computation that the coherent sheaf associated to $2 D(\pi)$ (twice the divisor of $\sigma_\pi$; equivalently, it can be defined as~$(\sigma_\pi^{\otimes 2})''$), which we denote by slight abuse of notation~$\sigma_\pi^2$, is
\begin{gather*}
\sigma_\pi^2 = \kappa_\pi \otimes O(n-2k).
\end{gather*}
In particular, if $n$ is even, $\kappa_\pi$ possesses a ``square root'', given by $\sigma_\pi \otimes O(k-n/2)$, and by multiplying the formula of Corollary~\ref{cor:sqrt} by $\prod_i y_i^{n/2-k}$, we obtain a formula for the class of the square root of the canonical sheaf, as advertised in the introduction.

In general, $[\sigma_\pi^2]\ne [\sigma_\pi]^2$ because $D(\pi)$ is not Cartier. In fact, $D(\pi)$ is Cartier iff $\CS_{\cl(\pi)}$, or equivalently $S_{\cl(\pi)}$, is Gorenstein, in which case all three sheaves (structure sheaf, $\sigma_\pi$, $\kappa_\pi$) are some $O(\cdot)$ sheaves. The condition for $S_{\cl(\pi)}$ to be Gorenstein is known \cite{Svanes,WY-SchubGroth}: in the language of link patterns, it means that all neighboring arcs form a~single group, e.g., link patterns typically look like
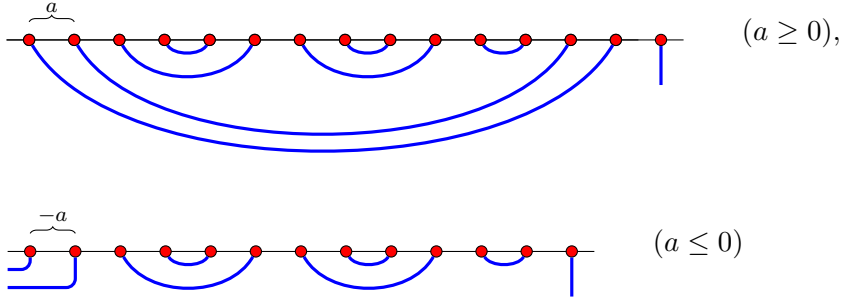
\begin{gather*}
 \tikz[baseline=0,/linkpattern/every linkpattern]{
\linkpattern[numbered=false]{3/6,4/5,7/10,8/9,11/12,15/15}
\linkpattern[numbered=false,squareness=0.25,looseness=0.15]{1/14,2/13}
\draw[decorate,decoration=brace] (1,0.4) -- node[above] {$\ss a$} (2,0.4);
} \qquad (a\ge 0),
\\
\tikz[baseline=0,/linkpattern/every linkpattern]{
\linkpattern[numbered=false]{1/1/\toleft,2/2/\toleft,3/6,4/5,7/10,8/9,11/12,13/13}
\draw[decorate,decoration=brace] (1,0.4) -- node[above] {$\ss -a$} (2,0.4);
}
\qquad (a\le 0)
\end{gather*}
for $\sigma_\pi=O(a)$. In this case, partition functions of crossing and noncrossing loop models give the same result up to a monomial; in a way, the noncrossing loop model is ``smarter'' because it produces the result with fewer configurations. We give such an example now.

\begin{ex}Let us consider the same link pattern as in Example~\ref{ex:n4k2}, namely \linkpattern{1/2,3/4}. We draw the corresponding loop configurations:
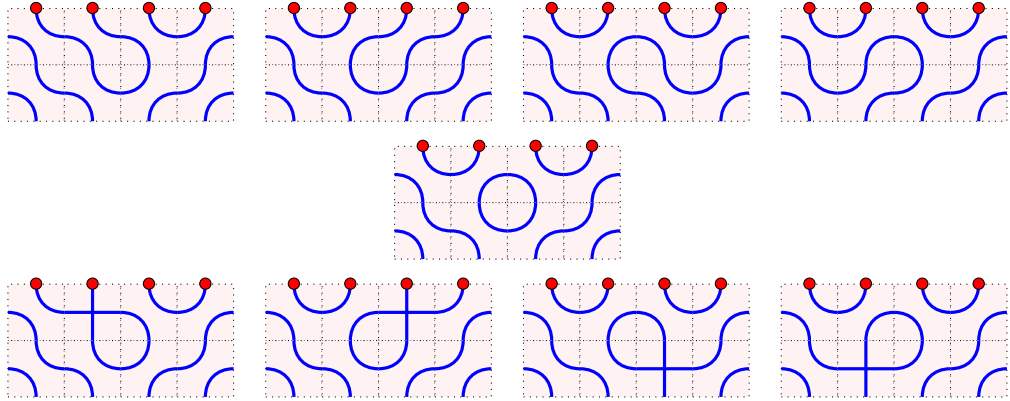
\begin{center}
\tikz{
\node[loop] {\plaq{b}&\plaq{b}&\plaq{b}&\plaq{a}\\\plaq{b}&\plaq{b}&\plaq{a}&\plaq{a}\\};
\foreach\x in {1,...,4} \node[/linkpattern/vertex] at (loop-1-\x.north) {};
}
\tikz{
\node[loop] {\plaq{b}&\plaq{a}&\plaq{a}&\plaq{a}\\\plaq{b}&\plaq{b}&\plaq{a}&\plaq{a}\\};
\foreach\x in {1,...,4} \node[/linkpattern/vertex] at (loop-1-\x.north) {};
}
\tikz{
\node[loop] {\plaq{b}&\plaq{a}&\plaq{b}&\plaq{a}\\\plaq{b}&\plaq{b}&\plaq{b}&\plaq{a}\\};
\foreach\x in {1,...,4} \node[/linkpattern/vertex] at (loop-1-\x.north) {};
}
\tikz{
\node[loop] {\plaq{b}&\plaq{a}&\plaq{b}&\plaq{a}\\\plaq{b}&\plaq{a}&\plaq{a}&\plaq{a}\\};
\foreach\x in {1,...,4} \node[/linkpattern/vertex] at (loop-1-\x.north) {};
}
\\
\tikz{
\node[loop] {\plaq{b}&\plaq{a}&\plaq{b}&\plaq{a}\\\plaq{b}&\plaq{b}&\plaq{a}&\plaq{a}\\};
\foreach\x in {1,...,4} \node[/linkpattern/vertex] at (loop-1-\x.north) {};
}
\\
\tikz{
\node[loop] {\plaq{b}&\plaq{c}&\plaq{b}&\plaq{a}\\\plaq{b}&\plaq{b}&\plaq{a}&\plaq{a}\\};
\foreach\x in {1,...,4} \node[/linkpattern/vertex] at (loop-1-\x.north) {};
}
\tikz{
\node[loop] {\plaq{b}&\plaq{a}&\plaq{c}&\plaq{a}\\\plaq{b}&\plaq{b}&\plaq{a}&\plaq{a}\\};
\foreach\x in {1,...,4} \node[/linkpattern/vertex] at (loop-1-\x.north) {};
}
\tikz{
\node[loop] {\plaq{b}&\plaq{a}&\plaq{b}&\plaq{a}\\\plaq{b}&\plaq{b}&\plaq{c}&\plaq{a}\\};
\foreach\x in {1,...,4} \node[/linkpattern/vertex] at (loop-1-\x.north) {};
}
\tikz{
\node[loop] {\plaq{b}&\plaq{a}&\plaq{b}&\plaq{a}\\\plaq{b}&\plaq{c}&\plaq{a}&\plaq{a}\\};
\foreach\x in {1,...,4} \node[/linkpattern/vertex] at (loop-1-\x.north) {};
}
\end{center}
Applying Corollary~\ref{cor:struc} (i.e., summing up the weights of these configurations of the {\em crossing} loop model), we find
\begin{gather*}
[X_{\linkpattern[small]{1/2,3/4}}]\\
\qquad{} =\big((1-t)^{2}(1-t y_1/y_2)(1-t y_2/y_1)\big)^{-1} (1-t z_4/y_1)(1-t z_4/y_2)(1-y_2/z_1)(1-y_1/z_1)\\
\qquad\quad{}\times \big({-}t^3 z_2 z_3 y_2^{-1} y_1^{-1}+2 t^2 z_2 y_1^{-1}
+2 t^2 z_3 y_1^{-1}-t^2 z_2 z_3 y_2^{-1} y_1^{-1}+2 t^2 z_2 y_2^{-1}+2 t^2 z_3 y_2^{-1}\\
 \qquad\quad{}-t^2 y_2 y_1^{-1} -t^2 y_1 y_2^{-1}
-t^2 z_2 z_3^{-1}-t^2 z_3 z_2^{-1}+2 t y_1 z_2^{-1}-t y_2 y_1 z_2^{-1} z_3^{-1}+2 t y_1 z_3^{-1}\\
\qquad\quad{}+2 t y_2 z_2^{-1}+2 t y_2 z_3^{-1} -t y_2 y_1^{-1}-t y_1 y_2^{-1}-t z_2 z_3^{-1} -t z_3 z_2^{-1}-y_2 y_1 z_2^{-1} z_3^{-1}\\
\qquad\quad{}+t^3-3 t^2-3 t+1\big).
 \end{gather*}
One can check that once specialized to $y_i=z_{I_i}$, we recover the first line of the table of restrictions to fixed points of Example~\ref{ex:n4k2}.

If we apply instead Corollary~\ref{cor:sqrt} (i.e., we only sum over the noncrossing configurations among these, giving them the weights of the noncrossing loop model), we find almost the same expression:
\begin{gather*}
[\sigma_{\linkpattern[small]{1/2,3/4}}]=t^3 z_2^{-1}z_3^{-1}z_4^{-2} [X_{\linkpattern[small]{1/2,3/4}}].
\end{gather*}
\end{ex}

\subsubsection{Gorenstein orbital varieties}
We apply $\mu'$, sending $\CS_{\cl(\pi)}$ to the orbital variety $\mathcal O_\pi$, and $\sigma_\pi$ to the module $\mu'_\ast\sigma_\pi$. Once again, note that the module associated to twice the divisor of~$\mu'_\ast\sigma_\pi$ is nothing but the canonical module of~$\mathcal O_\pi$. If (and only if) $\sigma_\pi$ is (nonequivariantly) trivial, $\mu'_\ast\sigma_\pi$ is free and its $K$-class coincides (up to a monomial, accounting for the grading shift) with that of the structure sheaf of~$\mathcal O_\pi$. In other words, we have to further distinguish a special case in the paragraph above (where we discussed the situation~$\sigma_\pi=O(a)$), which is when $a=0$.

For simplicity, we state the following result only for $n=2k$. Using in particular \cite{Stanley-Groth}, we obtain the following:
\begin{thm}\label{thm:goren}
Let $\pi\in\LP^{(k)}_{k,n}$, $n=2k$. The following conditions are equivalent:
\begin{itemize}\itemsep=0pt
\item $\mathcal O_\pi$ is Gorenstein.
\item $\mu_\ast \sigma_\pi$ is isomorphic $($up to a grading shift$)$ to both structure and canonical sheaves of $\mathcal O_\pi$.
\item $\Psi_\pi^c$ has degree range $($i.e., highest degree minus lowest degree$)$ $k(k-1)$ in the variables $z_1,\ldots,z_n$.
\item $\Psi_\pi^c$ is palindromic.
\item $\Psi_\pi^c=m_\pi^{-1} \Psi_\pi^{nc}$ where $m_\pi$ is the monomial defined in~\eqref{eq:weinorm}.
\item The leading coefficient of $P_\pi$ is $1$.
\item $\pi$ avoids $\linkpattern[numbered=false,squareness=0.25]{1/6,2/3,4/5}$ $($i.e., no restriction of $\pi$ to a subset of vertices forms this link pattern$)$.
\item $\pi$ is a series of nested arcs, i.e., of the form
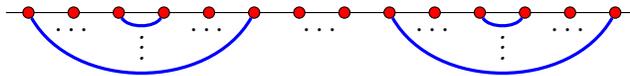
\begin{gather*}
\tikz[/linkpattern/every linkpattern]{
\linkpattern[numbered=false]{1/6,3/4,9/14,11/12}
\node at (2,-0.4) {$\cdots$};
\node at (3.5,-0.6) {$\vdots$};
\node at (5,-0.4) {$\cdots$};
\node at (7.5,-0.4) {$\cdots$};
\node at (10,-0.4) {$\cdots$};
\node at (11.5,-0.6) {$\vdots$};
\node at (13,-0.4) {$\cdots$};
}
\end{gather*}
\end{itemize}
\end{thm}
Only the statement about the leading coefficient of $P_\pi$ is worth explaining. The leading coefficient in $\beta$ corresponds to the limit $t\to 0$ where one is ``killing the fiber'' (one is considering functions that are constant on the fibers, i.e., functions on the base of the conormal variety). In other words, $\Psi^{(1)}_\pi(t\to 0)$ is simply the dimension of the space of global sections of~$\underline{\sigma}_\pi$. It is not hard to get an explicit combinatorial description of the space of global sections of any sheaf associated to an effective divisor on a Schubert variety, and check that its dimension is equal to one if and only if the divisor is zero.

\subsection{Lattice paths and TSSCPPs}\label{sec:TSSCPP}
We finally discuss explicit combinatorial expressions for certain entries $P_\pi(\beta)$ (as introduced in Section~\ref{sec:qKZ}), which happen to correspond to $S_{\cl(\pi)}$ Gorenstein (as discussed
in Section~\ref{sec:goren}).

The link pattern $\pi\in \LP^{(k)}_{k,n}$ is taken to be of the form of $a$ arcs atop a series of $k-a$ arcs connecting neighbors:
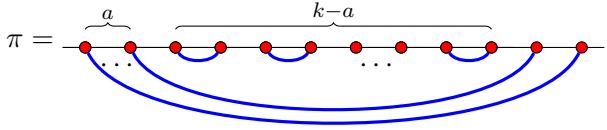
\begin{gather*}
\pi=
\tikz[baseline=0,/linkpattern/every linkpattern]{
\linkpattern[numbered=false]{3/4,5/6,9/10}
\linkpattern[numbered=false,squareness=0.2,looseness=0.1]{1/12,2/11}
\draw[decorate,decoration=brace] (1,0.4) -- node[above] {$\ss a$} (2,0.4);
\draw[decorate,decoration=brace] (3,0.4) -- node[above] {$\ss k-a$} (10,0.4);
\node at (7.5,-0.4) {$\cdots$};
\node at (1.75,-0.4) {$\cdots$};
}
\end{gather*}

Correcting some small mistakes in \cite{artic47}, we have the following expression for $P_\pi(\beta)$. It is the partition function of certain {\em non-intersecting lattice paths} on the square lattice, i.e., paths that are made of right steps and down steps and are not allowed to touch, with prescribed starting points, endpoints and weights which we describe now. There are $k-a-1$ paths; the starting points are fixed, equal to $(i,2i)$, $a+1\le i\le k-1$; and the endpoints are variable, of the form $(r_i-1,r_i)$, where $r_i$ is an integer which satisfies $r_i>2a$. Here is an example with $k=5$, $a=1$:
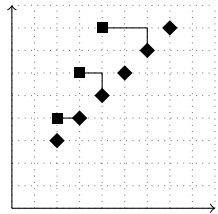
\begin{center}
\begin{tikzpicture}[scale=0.3]
\draw[gray,dotted] (0,0) grid (9,9);
\draw[->] (0,0) -- (9,0);
\draw[->] (0,0) -- (0,9);
\foreach\i in {2,...,4}
 \node[rectangle,fill,inner sep=2pt] at (\i,2*\i) {};
\foreach\i in {2,...,7}
 \node[diamond,fill,inner sep=1.5pt] at (\i,\i+1) {};
\draw (2,4) -- (3,4);
\draw (3,6) -- (4,6) -- (4,5);
\draw (4,8) -- (6,8) -- (6,7);
\end{tikzpicture}
\end{center}

We give a weight of $\beta$ to each vertical step of each path, with an additional {\em parity} contribution: if the horizontal displacement of a given path is equal to~$a$ modulo~2, then we include an extra~$\beta$. (In the example above, the weight is~$\beta^4$.) This results in the formula:
\begin{gather*}
P_\pi(\beta)=\sum_{\rm NILPs} \beta^{\sum\limits_{i=a+1}^{k-1} (2i-r_i+(r_i-i-a\mod2))}.
\end{gather*}

{\em Totally symmetric self-complementary plane partitions} (TSSCPPs) of size $k-1$ are defined as lozenge tilings of a regular hexagon of edge length $2(k-1)$ which possess all the symmetries of the hexagon.

The NILPs above are known to be in bijection with TSSCPPs of size $k-1$ with a frozen central hexagon of size $2a$ \cite{artic47}. In particular, for $a=0$, $P_\pi(\beta)$ is a certain weighted enumeration of TSSCPPs of size $k-1$; but in this case, according to Claim~\ref{thm:goren}, $\mu_\ast \sigma_\pi$ is isomorphic to the structure sheaf of~$\mathcal O_\pi$. Translating this result into elementary terms, we obtain the second claim of the introduction.

\begin{ex}
If $\pi=\linkpattern{1/2,3/4,5/6,7/8}$, we find the following 7 NILPs/TSSCPPs:
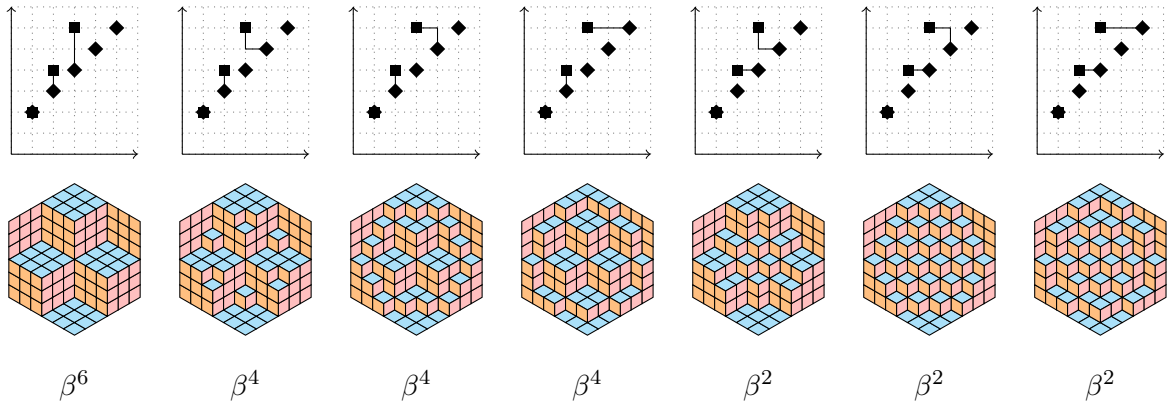
\begin{center}
\def\a{6}\def\b{6}\def\c{6}
\begin{tikzpicture}[scale=0.28]
\begin{scope}[xshift=-3cm,yshift=5cm]
\draw[gray,dotted] (0,0) grid (6,7);
\draw[->] (0,0) -- (6,0);
\draw[->] (0,0) -- (0,7);
\foreach\i in {1,...,3}
 \node[rectangle,fill,inner sep=2pt] at (\i,2*\i) {};
\foreach\i in {1,...,5}
 \node[diamond,fill,inner sep=1.5pt] at (\i,\i+1) {};
\draw (2,4) -- (2,3);
\draw (3,6) -- (3,4);
\end{scope}
\node at (0,-6) {$\beta^6$};
\loz%
{(4, 4), (5, 4), (6, 4), (10, 1), (11, 1), (12, 1), (4, 5), (5, 5), (6, 5), (10, 2), (11, 2), (12, 2), (4, 6), (5, 6), (6, 6), (10, 3), (11, 3), (12, 3), (1, 10), (2, 10), (3, 10), (7, 7), (8, 7), (9, 7), (1, 11), (2, 11), (3, 11), (7, 8), (8, 8), (9, 8), (1, 12), (2, 12), (3, 12), (7, 9), (8, 9), (9, 9)}
{(9, 4), (9, 5), (9, 6), (12, 10), (12, 11), (12, 12), (8, 4), (8, 5), (8, 6), (11, 10), (11, 11), (11, 12), (7, 4), (7, 5), (7, 6), (10, 10), (10, 11), (10, 12), (3, 1), (3, 2), (3, 3), (6, 7), (6, 8), (6, 9), (2, 1), (2, 2), (2, 3), (5, 7), (5, 8), (5, 9), (1, 1), (1, 2), (1, 3), (4, 7), (4, 8), (4, 9)}
{(1, 1), (2, 1), (3, 1), (7, 4), (8, 4), (9, 4), (1, 2), (2, 2), (3, 2), (7, 5), (8, 5), (9, 5), (1, 3), (2, 3), (3, 3), (7, 6), (8, 6), (9, 6), (4, 7), (5, 7), (6, 7), (10, 10), (11, 10), (12, 10), (4, 8), (5, 8), (6, 8), (10, 11), (11, 11), (12, 11), (4, 9), (5, 9), (6, 9), (10, 12), (11, 12), (12, 12)}
\end{tikzpicture}
\quad
\begin{tikzpicture}[scale=0.28]
\begin{scope}[xshift=-3cm,yshift=5cm]
\draw[gray,dotted] (0,0) grid (6,7);
\draw[->] (0,0) -- (6,0);
\draw[->] (0,0) -- (0,7);
\foreach\i in {1,...,3}
 \node[rectangle,fill,inner sep=2pt] at (\i,2*\i) {};
\foreach\i in {1,...,5}
 \node[diamond,fill,inner sep=1.5pt] at (\i,\i+1) {};
\draw (2,4) -- (2,3);
\draw (3,6) -- (3,5) -- (4,5);
\end{scope}
\loz%
{(4, 4), (5, 4), (7, 3), (10, 1), (11, 1), (12, 1), (4, 5), (5, 5), (6, 5), (10, 2), (11, 2), (12, 2), (3, 7), (5, 6), (6, 6), (9, 4), (11, 3), (12, 3), (1, 10), (2, 10), (4, 9), (7, 7), (8, 7), (10, 6), (1, 11), (2, 11), (3, 11), (7, 8), (8, 8), (9, 8), (1, 12), (2, 12), (3, 12), (6, 10), (8, 9), (9, 9)}
{(9, 4), (9, 5), (10, 7), (12, 10), (12, 11), (12, 12), (8, 4), (8, 5), (8, 6), (11, 10), (11, 11), (11, 12), (6, 3), (7, 5), (7, 6), (9, 9), (10, 11), (10, 12), (3, 1), (3, 2), (4, 4), (6, 7), (6, 8), (7, 10), (2, 1), (2, 2), (2, 3), (5, 7), (5, 8), (5, 9), (1, 1), (1, 2), (1, 3), (3, 6), (4, 8), (4, 9)}
{(1, 1), (2, 1), (3, 1), (6, 3), (8, 4), (9, 4), (1, 2), (2, 2), (3, 2), (7, 5), (8, 5), (9, 5), (1, 3), (2, 3), (4, 4), (7, 6), (8, 6), (10, 7), (3, 6), (5, 7), (6, 7), (9, 9), (11, 10), (12, 10), (4, 8), (5, 8), (6, 8), (10, 11), (11, 11), (12, 11), (4, 9), (5, 9), (7, 10), (10, 12), (11, 12), (12, 12)}
\node at (0,-6) {$\beta^4$};
\end{tikzpicture}
\quad
\begin{tikzpicture}[scale=0.28]
\begin{scope}[xshift=-3cm,yshift=5cm]
\draw[gray,dotted] (0,0) grid (6,7);
\draw[->] (0,0) -- (6,0);
\draw[->] (0,0) -- (0,7);
\foreach\i in {1,...,3}
 \node[rectangle,fill,inner sep=2pt] at (\i,2*\i) {};
\foreach\i in {1,...,5}
 \node[diamond,fill,inner sep=1.5pt] at (\i,\i+1) {};
\draw (2,4) -- (2,3);
\draw (3,6) -- (4,6) -- (4,5);
\end{scope}
\loz{(4, 4), (6, 3), (8, 2), (10, 1), (11, 1), (12, 1), (3, 6), (5, 5), (6, 5), (9, 3), (11, 2), (12, 2), (2, 8), (5, 6), (6, 6), (9, 4), (10, 4), (12, 3), (1, 10), (3, 9), (4, 9), (7, 7), (8, 7), (11, 5), (1, 11), (2, 11), (4, 10), (7, 8), (8, 8), (10, 7), (1, 12), (2, 12), (3, 12), (5, 11), (7, 10), (9, 9)}{(9, 4), (10, 6), (11, 8), (12, 10), (12, 11), (12, 12), (7, 3), (8, 5), (8, 6), (10, 9), (11, 11), (11, 12), (5, 2), (7, 5), (7, 6), (9, 9), (9, 10), (10, 12), (3, 1), (4, 3), (4, 4), (6, 7), (6, 8), (8, 11), (2, 1), (2, 2), (3, 4), (5, 7), (5, 8), (6, 10), (1, 1), (1, 2), (1, 3), (2, 5), (3, 7), (4, 9)}{(1, 1), (2, 1), (3, 1), (5, 2), (7, 3), (9, 4), (1, 2), (2, 2), (4, 3), (7, 5), (8, 5), (10, 6), (1, 3), (3, 4), (4, 4), (7, 6), (8, 6), (11, 8), (2, 5), (5, 7), (6, 7), (9, 9), (10, 9), (12, 10), (3, 7), (5, 8), (6, 8), (9, 10), (11, 11), (12, 11), (4, 9), (6, 10), (8, 11), (10, 12), (11, 12), (12, 12)}
\node at (0,-6) {$\beta^4$};
\end{tikzpicture}
\quad
\begin{tikzpicture}[scale=0.28]
\begin{scope}[xshift=-3cm,yshift=5cm]
\draw[gray,dotted] (0,0) grid (6,7);
\draw[->] (0,0) -- (6,0);
\draw[->] (0,0) -- (0,7);
\foreach\i in {1,...,3}
 \node[rectangle,fill,inner sep=2pt] at (\i,2*\i) {};
\foreach\i in {1,...,5}
 \node[diamond,fill,inner sep=1.5pt] at (\i,\i+1) {};
\draw (2,4) -- (2,3);
\draw (3,6) -- (5,6);
\end{scope}
\loz%
{(4, 4), (7, 2), (8, 2), (10, 1), (11, 1), (12, 1), (2, 7), (5, 5), (6, 5), (9, 3), (10, 3), (12, 2), (2, 8), (5, 6), (6, 6), (9, 4), (10, 4), (12, 3), (1, 10), (3, 9), (4, 9), (7, 7), (8, 7), (11, 5), (1, 11), (3, 10), (4, 10), (7, 8), (8, 8), (11, 6), (1, 12), (2, 12), (3, 12), (5, 11), (6, 11), (9, 9)}
{(9, 4), (11, 7), (11, 8), (12, 10), (12, 11), (12, 12), (6, 2), (8, 5), (8, 6), (10, 9), (10, 10), (11, 12), (5, 2), (7, 5), (7, 6), (9, 9), (9, 10), (10, 12), (3, 1), (4, 3), (4, 4), (6, 7), (6, 8), (8, 11), (2, 1), (3, 3), (3, 4), (5, 7), (5, 8), (7, 11), (1, 1), (1, 2), (1, 3), (2, 5), (2, 6), (4, 9)}
{(1, 1), (2, 1), (3, 1), (5, 2), (6, 2), (9, 4), (1, 2), (3, 3), (4, 3), (7, 5), (8, 5), (11, 7), (1, 3), (3, 4), (4, 4), (7, 6), (8, 6), (11, 8), (2, 5), (5, 7), (6, 7), (9, 9), (10, 9), (12, 10), (2, 6), (5, 8), (6, 8), (9, 10), (10, 10), (12, 11), (4, 9), (7, 11), (8, 11), (10, 12), (11, 12), (12, 12)}
\node at (0,-6) {$\beta^4$};
\end{tikzpicture}
\quad
\begin{tikzpicture}[scale=0.28]
\begin{scope}[xshift=-3cm,yshift=5cm]
\draw[gray,dotted] (0,0) grid (6,7);
\draw[->] (0,0) -- (6,0);
\draw[->] (0,0) -- (0,7);
\foreach\i in {1,...,3}
 \node[rectangle,fill,inner sep=2pt] at (\i,2*\i) {};
\foreach\i in {1,...,5}
 \node[diamond,fill,inner sep=1.5pt] at (\i,\i+1) {};
\draw (2,4) -- (3,4);
\draw (3,6) -- (3,5) -- (4,5);
\end{scope}
\loz%
{(4, 4), (5, 4), (7, 3), (10, 1), (11, 1), (12, 1), (4, 5), (5, 5), (7, 4), (10, 2), (11, 2), (12, 2), (3, 7), (4, 7), (6, 6), (8, 5), (11, 3), (12, 3), (1, 10), (2, 10), (5, 8), (7, 7), (9, 6), (10, 6), (1, 11), (2, 11), (3, 11), (6, 9), (8, 8), (9, 8), (1, 12), (2, 12), (3, 12), (6, 10), (8, 9), (9, 9)}
{(9, 4), (9, 5), (10, 7), (12, 10), (12, 11), (12, 12), (8, 4), (8, 5), (9, 7), (11, 10), (11, 11), (11, 12), (6, 3), (6, 4), (7, 6), (8, 8), (10, 11), (10, 12), (3, 1), (3, 2), (5, 5), (6, 7), (7, 9), (7, 10), (2, 1), (2, 2), (2, 3), (4, 6), (5, 8), (5, 9), (1, 1), (1, 2), (1, 3), (3, 6), (4, 8), (4, 9)}
{(1, 1), (2, 1), (3, 1), (6, 3), (8, 4), (9, 4), (1, 2), (2, 2), (3, 2), (6, 4), (8, 5), (9, 5), (1, 3), (2, 3), (5, 5), (7, 6), (9, 7), (10, 7), (3, 6), (4, 6), (6, 7), (8, 8), (11, 10), (12, 10), (4, 8), (5, 8), (7, 9), (10, 11), (11, 11), (12, 11), (4, 9), (5, 9), (7, 10), (10, 12), (11, 12), (12, 12)}
\node at (0,-6) {$\beta^2$};
\end{tikzpicture}
\quad
\begin{tikzpicture}[scale=0.28]
\begin{scope}[xshift=-3cm,yshift=5cm]
\draw[gray,dotted] (0,0) grid (6,7);
\draw[->] (0,0) -- (6,0);
\draw[->] (0,0) -- (0,7);
\foreach\i in {1,...,3}
 \node[rectangle,fill,inner sep=2pt] at (\i,2*\i) {};
\foreach\i in {1,...,5}
 \node[diamond,fill,inner sep=1.5pt] at (\i,\i+1) {};
\draw (2,4) -- (3,4);
\draw (3,6) -- (4,6) -- (4,5);
\end{scope}
\loz%
{(4, 4), (6, 3), (8, 2), (10, 1), (11, 1), (12, 1), (3, 6), (5, 5), (7, 4), (9, 3), (11, 2), (12, 2), (2, 8), (4, 7), (6, 6), (8, 5), (10, 4), (12, 3), (1, 10), (3, 9), (5, 8), (7, 7), (9, 6), (11, 5), (1, 11), (2, 11), (4, 10), (6, 9), (8, 8), (10, 7), (1, 12), (2, 12), (3, 12), (5, 11), (7, 10), (9, 9)}
{(9, 4), (10, 6), (11, 8), (12, 10), (12, 11), (12, 12), (7, 3), (8, 5), (9, 7), (10, 9), (11, 11), (11, 12), (5, 2), (6, 4), (7, 6), (8, 8), (9, 10), (10, 12), (3, 1), (4, 3), (5, 5), (6, 7), (7, 9), (8, 11), (2, 1), (2, 2), (3, 4), (4, 6), (5, 8), (6, 10), (1, 1), (1, 2), (1, 3), (2, 5), (3, 7), (4, 9)}
{(1, 1), (2, 1), (3, 1), (5, 2), (7, 3), (9, 4), (1, 2), (2, 2), (4, 3), (6, 4), (8, 5), (10, 6), (1, 3), (3, 4), (5, 5), (7, 6), (9, 7), (11, 8), (2, 5), (4, 6), (6, 7), (8, 8), (10, 9), (12, 10), (3, 7), (5, 8), (7, 9), (9, 10), (11, 11), (12, 11), (4, 9), (6, 10), (8, 11), (10, 12), (11, 12), (12, 12)}
\node at (0,-6) {$\beta^2$};
\end{tikzpicture}
\quad
\begin{tikzpicture}[scale=0.28]
\begin{scope}[xshift=-3cm,yshift=5cm]
\draw[gray,dotted] (0,0) grid (6,7);
\draw[->] (0,0) -- (6,0);
\draw[->] (0,0) -- (0,7);
\foreach\i in {1,...,3}
 \node[rectangle,fill,inner sep=2pt] at (\i,2*\i) {};
\foreach\i in {1,...,5}
 \node[diamond,fill,inner sep=1.5pt] at (\i,\i+1) {};
\draw (2,4) -- (3,4);
\draw (3,6) -- (5,6);
\end{scope}
\loz%
{(4, 4), (7, 2), (8, 2), (10, 1), (11, 1), (12, 1), (2, 7), (5, 5), (7, 4), (9, 3), (10, 3), (12, 2), (2, 8), (4, 7), (6, 6), (8, 5), (10, 4), (12, 3), (1, 10), (3, 9), (5, 8), (7, 7), (9, 6), (11, 5), (1, 11), (3, 10), (4, 10), (6, 9), (8, 8), (11, 6), (1, 12), (2, 12), (3, 12), (5, 11), (6, 11), (9, 9)}
{(9, 4), (11, 7), (11, 8), (12, 10), (12, 11), (12, 12), (6, 2), (8, 5), (9, 7), (10, 9), (10, 10), (11, 12), (5, 2), (6, 4), (7, 6), (8, 8), (9, 10), (10, 12), (3, 1), (4, 3), (5, 5), (6, 7), (7, 9), (8, 11), (2, 1), (3, 3), (3, 4), (4, 6), (5, 8), (7, 11), (1, 1), (1, 2), (1, 3), (2, 5), (2, 6), (4, 9)}
{(1, 1), (2, 1), (3, 1), (5, 2), (6, 2), (9, 4), (1, 2), (3, 3), (4, 3), (6, 4), (8, 5), (11, 7), (1, 3), (3, 4), (5, 5), (7, 6), (9, 7), (11, 8), (2, 5), (4, 6), (6, 7), (8, 8), (10, 9), (12, 10), (2, 6), (5, 8), (7, 9), (9, 10), (10, 10), (12, 11), (4, 9), (7, 11), (8, 11), (10, 12), (11, 12), (12, 12)}
\node at (0,-6) {$\beta^2$};
\end{tikzpicture}
\end{center}
(to recover the NILP from the TSSCPP, keep track of pink and blue lozenges in the region \tikz[scale=0.5,baseline=0]{\draw[clip] (0,-1) -- ++(30:1) -- ++(90:1) -- ++(150:1) -- ++(210:1) -- ++(270:1) -- ++(330:1);
\draw[line cap=round,red,ultra thick] (30:1) -- (0,0) -- (60:0.866); }) hence a polynomial $P_\pi(\beta)=\beta^6+3\beta^4+3\beta^2$. In particular, the degree of $\mathcal O_\pi$ is $P_\pi(\beta=2)=124$.
\end{ex}

\section{Other loop models}
We conclude by mentioning a few possible generalizations and variations of the ideas of this paper.

\subsection{Other types}
\looseness=-1 The discussion above was entirely restricted to type A, in the sense that the group acting on our variety was $({\rm P})\operatorname{GL}_n$, with corresponding Weyl group the symmetric group. One should be able to adapt our method to other types. For example, one may consider the cotangent bundle of the {\em Lagrangian} Grassmannian; it should be related to various loop models with one (integrable) boundary.

\subsection{The Brauer loop model}
One should point out that another crossing loop model was discussed in connection with (ordinary) cohomology, namely the Brauer loop model \cite{artic33,artic39}. However, several important ingredients are missing in order to make this work fit into the framework of the current paper. In particular, there, we only have the analogue of the orbital varieties~$\mathcal O_\pi$, not of the conormal Schubert varieties. This prevents us at the moment from extending it to $K$-theory.

\subsection{The dilute loop model}
There is another important noncrossing loop model, which is the {\em dilute} loop model~\cite{Nienhuis-loop}. In the same way that the (dense) noncrossing loop model discussed here is related to desingularizations of the nilpotent orbit $\{u^2=0\}$, the dilute noncrossing loop model should be related to desingularizations of $\{u^2=v^3\}$. This will be discussed elsewhere. It is not clear whether there will also be a crossing loop model associated to that geometry.

\subsection[Beyond $K$-theory]{Beyond $\boldsymbol{K}$-theory}
We conclude by saying that $K$-theory is of course not the most general complex-oriented cohomology theory. Even if we restrict ourselves to those whose formal group law is an actual group law on a curve (which seems natural on the integrable side), then we should consider elliptic cohomology, as in~\cite{AO-ell}. Since no elliptic weights are known for loop models, it is unclear how to generalize our work in this direction. However, one can extend a little less by going over to the singular elliptic cohomology considered in~\cite{LZZ-para}, where loop models should still play a role.

\vspace{-1mm}

\subsection*{Acknowledgements}

PZJ was supported by ERC grant 278124 and ARC grant FT150100232. Computerized checks of the results of this paper were performed with the help of Macaulay~2~\cite{M2}.

\vspace{-1mm}

%\pdfbookmark[1]{References}{ref}
\LastPageEnding

\end{document}